\documentclass[11pt,a4paper]{article}
\usepackage[T1]{fontenc}
\usepackage[utf8]{inputenc}
\usepackage{lmodern}
\usepackage[margin=1in]{geometry}
\usepackage{amsmath,amssymb,amsthm,mathtools,bm}
\usepackage{booktabs,tabularx,array}
\usepackage{placeins,graphicx}
\usepackage{enumitem}
\usepackage[authoryear,round]{natbib}
\usepackage{microtype}
\usepackage{xcolor}
\definecolor{LinkBlue}{RGB}{0,0,160}
\usepackage[colorlinks=true,linkcolor=LinkBlue,citecolor=LinkBlue,urlcolor=LinkBlue]{hyperref}
\providecommand{\doi}[1]{doi:\,\href{https://doi.org/#1}{\nolinkurl{#1}}}
\hypersetup{pdftitle={An F-Type Normal-Reference Calibration for High-Dimensional Covariance Testing},pdfauthor={Pengfei Wang},pdfsubject={High-dimensional covariance inference},pdfkeywords={covariance equality, normal reference, F-type test, high-dimensional inference}}
\setlist{itemsep=2pt,topsep=4pt}
\newtheorem{theorem}{Theorem}
\newtheorem{proposition}{Proposition}
\newtheorem{lemma}{Lemma}
\newtheorem{corollary}{Corollary}
\theoremstyle{definition}
\newtheorem{assumption}{Assumption}
\theoremstyle{remark}
\newtheorem{remark}{Remark}
\DeclareMathOperator{\tr}{tr}
\DeclareMathOperator{\Cov}{Cov}
\DeclareMathOperator{\Var}{Var}
\DeclareMathOperator{\svec}{svec}

\DeclareMathOperator{\diag}{diag}

\newcommand{\E}{\mathbb E}
\newcommand{\Pp}{\mathbb P}
\newcommand{\R}{\mathbb R}
\newcommand{\ind}{\mathbb I}
\newcommand{\op}{o_{\mathrm p}}
\newcommand{\Op}{O_{\mathrm p}}

\newcommand{\norm}[1]{\left\lVert#1\right\rVert}

\title{An \(F\)-Type Normal-Reference Calibration for\\
High-Dimensional Covariance Testing}
\author{Pengfei Wang\\
\small Mathematics and Mathematics Education\\
\small National Institute of Education, Nanyang Technological University}
\date{October 22, 2026}

\begin{document}
\maketitle
\begin{abstract}
An \(F\)-type normal-reference calibration is developed for testing equality
of two high-dimensional covariance matrices with unknown means. The test
statistic is a ratio whose numerator measures the Frobenius covariance
discrepancy and whose denominator estimates its null expectation. Their
difference is the complete Li--Chen statistic. Under suitable moment and
trace conditions, we establish a Gaussian-reference approximation. After
standardization, the statistic and its reference share the normal or
nonnormal limits determined by a convergent normalized spectrum. The
reference is a ratio of
independent chi-square-type mixtures. Applying Welch--Satterthwaite
matching to its two components gives an ordinary \(F\) approximation,
with degrees of freedom estimated from independent four-observation
blocks. The effect of estimating the population means is controlled by
an exact expansion, and local-power and consistency results are
established. The ordinary \(F\) approximation is asymptotically valid
under a diffuse spectrum, whereas a covariance spike can produce a
nonvanishing approximation error. Simulation studies examine empirical
size and power under normal and nonnormal sampling. Averaging the trace
estimates over fixed partitions improves their stability, but does not
remove size inflation in the strongly correlated settings.
\end{abstract}

\section{Introduction}\label{sec:intro}
Covariance matrices describe both marginal variation and dependence among
variables. Comparing them is therefore of interest even when the
population means are not the same. For high-dimensional samples, however,
the sample covariance matrices may be singular, and classical tests based
on their inverses are no longer directly applicable. We consider the
following two-sample problem.

Let
\[
 \bm{X}_{i1},\ldots,\bm{X}_{in_i}\quad (i=1,2)
\]
be two independent samples, each consisting of independent and identically
distributed observations. Write
\[
 \E \bm{X}_{ij}=\bm{\mu}_i,\qquad \Cov(\bm{X}_{ij})=\bm{\Sigma}_i,
 \qquad \bm{X}_{ij}\in\R^p.
\]
The means are unknown and may differ. We consider
\begin{equation}\label{eq:null}
 H_0:\bm{\Sigma}_1=\bm{\Sigma}_2
 \quad\text{against}\quad H_1:\bm{\Sigma}_1\ne\bm{\Sigma}_2.
\end{equation}
Under the null hypothesis, the distributions and fourth moments may differ
between groups. The dimension \(p\) may increase with the sample sizes.
Vectors and matrices are denoted by boldface letters; scalars and individual
entries are set in ordinary mathematical type.

Several tests have been developed for this problem. Under Gaussian
sampling, \citet{Schott2007} and \citet{SrivastavaYanagihara2010}
constructed trace-based statistics that avoid matrix inversion.
\citet{LiChen2012} estimated the squared Frobenius distance by a
mean-corrected \(U\)-statistic and derived its asymptotic normal
distribution under moment and trace conditions. Maximum-type procedures,
such as that of \citet{CaiLiuXia2013}, are more directly suited to sparse
covariance differences. Their calibration can be improved by multiplier
methods that allow dependence \citep{ChangZhouZhouWang2017}; a related
jackknife construction for covariance \(U\)-statistics with unknown
means was studied by \citet{ChakrabortyKarmakarKoul2026}.

For a Frobenius-type statistic, the adequacy of normal approximation
depends on the covariance spectrum. If a few directions contribute a
substantial proportion of the variation, the limiting distribution can
be nonnormal. To accommodate this possibility,
\citet{WangZhuZhang2024} treated vectorized outer products as observations
and compared the induced mean statistic with a Gaussian reference having
the same feature means and covariances. They used the three-cumulant
chi-square approximation of \citet{Zhang2005} to obtain critical values.
Normal-reference methods have also been developed for several covariance
matrices \citep{WangZhuZhang2025} and covariance general linear hypotheses
\citep{WangZhuZhang2026}. The generalized quadratic-form results of
\citet{WangXu2022} provide a related basis for Gaussian comparison.

Motivated by the two-sample mean test of \citet{ZhuWangZhang2024}, we
consider an \(F\)-type calibration of the covariance discrepancy.
We divide the squared Frobenius distance between sample covariance
matrices by an unbiased estimate of its null expectation. Thus the
numerator and denominator have the same expectation under the null
hypothesis. Their difference is exactly the complete Li--Chen statistic,
so the construction preserves the covariance discrepancy while changing
its calibration. In the Gaussian reference, the numerator and denominator
are independent weighted chi-square sums. Matching the first two moments
of each sum gives an ordinary \(F\) distribution whose degrees of
freedom depend on feature-covariance trace products.

Two points need to be addressed to justify this procedure. First,
transforming covariance equality into a mean problem connects the test
to quadratic mean statistics \citep{BaiSaranadasa1996,ChenQin2010}, but
substituting sample means makes the outer-product features dependent.
We derive an exact expansion for the resulting correction and estimate
the required trace products from independent four-observation blocks.
These calculations establish the feasible reference approximation and
its normal and nonnormal limits under explicit conditions. Second,
validity of the Gaussian reference does not ensure accuracy of its
ordinary \(F\) approximation. We give a diffuse-spectrum condition
under which the matching error vanishes and a spiked model in which it
persists. The \(F\) and corresponding chi-square critical values are
also shown to be first-order equivalent. Local-power and consistency
results complete the theoretical analysis.

The main results are presented in Section~\ref{sec:mainresults}.
Section~\ref{sec:simulation} compares empirical sizes and powers, with
particular attention to the effect of size distortion on power comparisons.
Section~\ref{sec:discussion} gives concluding remarks. Proofs and
additional numerical results are collected in the appendices.

\section{Main results}\label{sec:mainresults}
\subsection{Test statistic}\label{sec:model}
Covariance equality can be expressed as equality of means after
vectorizing the outer products of population-centered observations.
Since these matrices are symmetric, we
use \(\svec(\bm{A})\in\R^d\), \(d=p(p+1)/2\), which contains the diagonal
entries of \(\bm{A}\) and its upper-triangular entries multiplied by
\(\sqrt2\), in a fixed order. This convention preserves the Frobenius
inner product:
\begin{equation}\label{eq:isometry}
 \svec(\bm{A})^\top\svec(\bm{B})=\tr(\bm{A}\bm{B}),\qquad
 \norm{\svec(\bm{A})}^2=\norm{\bm{A}}_F^2.
\end{equation}

\subsubsection{Quadratic features}
Set \(\bm{\epsilon}_{ij}=\bm{X}_{ij}-\bm{\mu}_i\) and define the oracle features
\begin{equation}\label{eq:oraclefeatures}
 \bm{Z}_{ij}^{\,o}=\svec(\bm{\epsilon}_{ij}\bm{\epsilon}_{ij}^\top),\qquad
 \bm{\theta}_i=\svec(\bm{\Sigma}_i),\qquad
 \bm{\Omega}_i^{\,o}=\Cov(\bm{Z}_{ij}^{\,o}).
\end{equation}
The hypothesis \eqref{eq:null} is equivalent to \(\bm{\theta}_1=\bm{\theta}_2\), with squared signal
\[
 \delta_n=\norm{\bm{\theta}_1-\bm{\theta}_2}^2
          =\norm{\bm{\Sigma}_1-\bm{\Sigma}_2}_F^2.
\]
Thus the squared distance between the feature means is precisely the
squared Frobenius distance between the covariance matrices. Since
\(\bm{\Omega}_1^{\,o}\) and \(\bm{\Omega}_2^{\,o}\) may differ under \(H_0\),
the induced mean problem is heteroscedastic.

The effect of estimating the means will also involve the Gaussian
quadratic-form covariance. Denote by \(\bm{\Gamma}_i\) the covariance matrix
of \(\svec(\bm{G}_i\bm{G}_i^\top)\), where \(\bm{G}_i\sim N_p(0,\bm{\Sigma}_i)\).
For symmetric matrices \(\bm{A},\bm{B}\), this covariance is characterized by
\begin{equation}\label{eq:Gamma}
 \Cov(\bm{G}_i^\top \bm{A} \bm{G}_i,\bm{G}_i^\top \bm{B} \bm{G}_i)
       =2\tr(\bm{A}\bm{\Sigma}_i\bm{B}\bm{\Sigma}_i).
\end{equation}
In the isometric coordinates \eqref{eq:isometry},
\begin{align}
 \tr(\bm{\Gamma}_i)&=\{\tr(\bm{\Sigma}_i)\}^2+\tr(\bm{\Sigma}_i^2),\label{eq:Gamma1}\\
 \tr(\bm{\Gamma}_i^2)&=2\left[\{\tr(\bm{\Sigma}_i^2)\}^2+\tr(\bm{\Sigma}_i^4)\right].\label{eq:Gamma2}
\end{align}
For Gaussian observations, \(\bm{\Omega}_i^{\,o}=\bm{\Gamma}_i\); their quadratic
features remain non-Gaussian.

\subsubsection{Construction with unknown means}\label{sec:fullstat}
To construct the test from the observed data, let
\[
 \bm{S}_i=\frac1{n_i-1}\sum_{j=1}^{n_i}
 (\bm{X}_{ij}-\bm{\bar X}_i)(\bm{X}_{ij}-\bm{\bar X}_i)^\top.
\]
A natural measure of departure from \(H_0\) is
\(Q_{\rm full}=\norm{\bm{S}_1-\bm{S}_2}_F^2\). Its null expectation is generally nonzero because the sample covariance
matrices fluctuate around their population values. We estimate this expectation and use it to normalize the distance.
To obtain a mean-invariant estimate, for a sample of size \(n\ge4\), set
\begin{equation}\label{eq:h}
 \bm{h}_{ab}=\frac12\svec\{(\bm{X}_a-\bm{X}_b)(\bm{X}_a-\bm{X}_b)^\top\},
 \qquad a\ne b.
\end{equation}
Pairwise differences remove the unknown mean, and the average of
\(\bm{h}_{ab}\) over unordered pairs is \(\svec(\bm{S})\). Products formed from
disjoint pairs give the estimator
\begin{equation}\label{eq:D}
 D_n=\frac1{(n)_4}\sum_{\substack{a,b,c,e\\\text{all distinct}}}
       \bm{h}_{ab}^\top \bm{h}_{ce},
 \qquad (n)_4=n(n-1)(n-2)(n-3),
\end{equation}
where the sum is over ordered quadruples. The independence of the two
pairs makes \(D_n\) unbiased for \(\tr(\bm{\Sigma}^2)\). It is a fourth-order
U-statistic \citep{Hoeffding1948}. We therefore estimate the sampling
contribution to the squared distance and construct the ratio as follows:
\begin{align}
 \widehat b_i&=\norm{\bm{S}_i}_F^2-D_{n_i},&
 B_{\rm full}&=\widehat b_1+\widehat b_2,\label{eq:Bfull}\\
 U_{\rm full}&=D_{n_1}+D_{n_2}-2\tr(\bm{S}_1\bm{S}_2),&
 F_{\rm full}&=\frac{Q_{\rm full}}{B_{\rm full}}
       =1+\frac{U_{\rm full}}{B_{\rm full}}.\label{eq:Ffull}
\end{align}
We use \(F_{\rm full}\) when \(B_{\rm full}>0\).

\begin{proposition}[Exact full-sample identities]\label{prop:full}
Suppose that \(\E\norm{\bm{\epsilon}_{i1}}^4<\infty\). Then
\begin{equation}\label{eq:covS}
 \Cov\{\svec(\bm{S}_i)\}
 =\frac{\bm{\Omega}_i^{\,o}}{n_i}
  +\frac{\bm{\Gamma}_i}{n_i(n_i-1)}.
\end{equation}
Moreover,
\[
 \E D_{n_i}=\tr(\bm{\Sigma}_i^2),\qquad
 \E\widehat b_i=
 \frac{\tr(\bm{\Omega}_i^{\,o})}{n_i}
 +\frac{\tr(\bm{\Gamma}_i)}{n_i(n_i-1)}=:b_i^{\rm full},
\]
and
\[
 \E U_{\rm full}=\norm{\bm{\Sigma}_1-\bm{\Sigma}_2}_F^2,\qquad
 \E Q_{\rm full}=b_1^{\rm full}+b_2^{\rm full}+\delta_n.
\]
All quantities in \eqref{eq:D}--\eqref{eq:Ffull} are translation invariant.
\end{proposition}

Proposition~\ref{prop:full} shows that the numerator and denominator
have the same expectation under \(H_0\), as required for an \(F\)-type
construction. The denominator is unbiased for
\(b_1^{\rm full}+b_2^{\rm full}\), but need not be positive: the scalar
sample \((-1,-1,1,1)\) gives \(S=4/3\), \(D_4=8/3\), and
\(\widehat b=-8/9\).
We retain the unmodified denominator and treat \(B_{\rm full}\le0\) as a
calibration failure, resulting in nonrejection. The simulations report
these failure frequencies.

\begin{proposition}[Identity with the Li--Chen statistic]\label{prop:LiChen}
For any samples with \(n_1,n_2\ge4\), let \(A_{n_i}^{\rm LC}\) and
\(C_{n_1,n_2}^{\rm LC}\) denote the complete mean-corrected within-group
and cross-group estimators in \citet[Equations (2.1)--(2.3)]{LiChen2012}.
Then, pointwise in the observed samples,
\begin{equation}\label{eq:LiChenidentity}
 D_{n_i}=A_{n_i}^{\rm LC},\qquad
 \tr(\bm{S}_1\bm{S}_2)=C_{n_1,n_2}^{\rm LC},\qquad
 U_{\rm full}=T_{\rm LC}.
\end{equation}
\end{proposition}
Thus \(F_{\rm full}=1+T_{\rm LC}/B_{\rm full}\) gives a ratio calibration
of the complete mean-corrected Li--Chen statistic.

\subsubsection{Gram-matrix computation}\label{sec:computation}
The fourth-order statistic \(D_n\) has a simple Gram-matrix expression.
Set \(\bm{r}_j=\bm{X}_j-\bm{\bar X}\), \(G_{jk}=\bm{r}_j^\top \bm{r}_k\), \(t_G=\tr(\bm{G})\),
\(u_G=\tr(\bm{G}^2)\), and \(q_G=\sum_jG_{jj}^2\). Then
\begin{equation}\label{eq:Dgram}
 D_n=\frac{(n-1)(n-2)u_G+t_G^2-n(n-1)q_G}{(n)_4}.
\end{equation}
Thus \(D_n\) is computable in \(O(n^2p)\) operations without enumerating
quadruples. The derivation is given in Appendix~\ref{app:identities}.

\subsection{Normal-reference approximation}
\label{sec:oraclereference}
To study the null distribution of the ratio, first consider the
independent quadratic features formed using the population means.
For the features in \eqref{eq:oraclefeatures}, write
\(\bm{u}_{ij}^o=\bm{Z}_{ij}^o-\bm{\theta}_i\), and let \(\bm{\bar Z}_i^o\) and
\(\bm{\widehat\Omega}_i^o\) denote the sample mean and covariance matrix.
Define
\[
 B^o=\sum_{i=1}^2\frac{\tr(\bm{\widehat\Omega}_i^o)}{n_i},\qquad
 Q^o=\norm{\bm{\bar Z}_1^o-\bm{\bar Z}_2^o}^2,\qquad U^o=Q^o-B^o.
\]
For the distributional calculations, let
\begin{align*}
 \bm{V}^o&=\frac{\bm{\Omega}_1^o}{n_1}+\frac{\bm{\Omega}_2^o}{n_2},&
 b^o&=\tr(\bm{V}^o),&
 J^o&=\tr\{(\bm{V}^o)^2\},\\
 s_i^o&=\tr\{(\bm{\Omega}_i^o)^2\},&
 R^o&=\sum_{i=1}^2\frac{s_i^o}{n_i^2(n_i-1)},&
 (a^o)^2&=2(J^o+R^o).
\end{align*}
Also put \(c_{12}^o=\tr(\bm{\Omega}_1^o\bm{\Omega}_2^o)\).
Under \(H_0\), \(U^o\) has mean zero and variance \((a^o)^2\).

Throughout the asymptotic analysis, \(b^o,a^o>0\),
\(n_{\min}=\min(n_1,n_2)\to\infty\), and \(n_1/(n_1+n_2)\) is
bounded away from zero and one. In addition to the sampling conditions in
Section~\ref{sec:model}, assume that a constant \(K<\infty\), independent
of the dimension and sample sizes, satisfies
\[
 \E(\bm{v}^\top \bm{u}_{i1}^o)^4\le K(\bm{v}^\top\bm{\Omega}_i^o \bm{v})^2
 \quad\text{for every }\bm{v}\in\R^d,\quad i=1,2.
\]
These are the oracle-feature versions of
Assumptions~\ref{ass:sampling}--\ref{ass:fourth}. The directional moment
bound concerns centered quadratic forms and requires control of the
original-data eighth moments. Section~\ref{sec:verifiedmodels} gives
sufficient conditions on the original observations.

The normal-reference approach replaces the feature observations by
independent Gaussian vectors with the same means and covariance matrices.
Specifically, let \(\bm{Z}_{ij}^{o,*}\sim N_d(\bm{\theta}_i,\bm{\Omega}_i^o)\), and
compute \(Q^{o,*},B^{o,*},U^{o,*}\) from these vectors by the same
formulas, with \(F^{o,*}=Q^{o,*}/B^{o,*}\). Under \(H_0\), the numerator
and denominator have the representations
\[
 Q^{o,*}\overset d=\sum_r\lambda_r^o\xi_r,\qquad
 B^{o,*}\overset d=\sum_{i=1}^2\sum_r
 \frac{\lambda_{ir}^o}{n_i(n_i-1)}\eta_{ir},
\]
where \(\lambda_r^o\) and \(\lambda_{ir}^o\) are the eigenvalues of
\(\bm{V}^o\) and \(\bm{\Omega}_i^o\), respectively,
\(\xi_r\sim\chi_1^2\), \(\eta_{ir}\sim\chi_{n_i-1}^2\), and all these
variables are independent. Thus \(F^{o,*}\) is a ratio of independent
chi-square mixtures. Its numerator and denominator have common mean
\(b^o\) and variances \(2J^o\) and \(2R^o\), respectively.
Under alternatives, define the central counterpart \(U_0^{o,*}\) using
independent \(N_d(0,\bm{\Omega}_i^o)\) features, where the covariance matrices
are those of the current populations.

To obtain a convenient critical value, approximate each mixture by a
scaled chi-square variable whose first two moments agree with those of
the mixture \citep{Satterthwaite1946,Welch1947}. Since the two mixtures
have the same mean, the ratio of the matched variables follows an ordinary
\(F_{d_1^o,d_2^o}\) distribution. The degrees of freedom and the associated
matching error are
\[
 d_1^o=\frac{(b^o)^2}{J^o},\qquad
 d_2^o=\frac{(b^o)^2}{R^o},\qquad
 \varepsilon_{\rm WS}^o
   =d_{\rm K}\{\mathcal L(F^{o,*}),F_{d_1^o,d_2^o}\}.
\]
Here \(d_{\rm K}\) denotes Kolmogorov distance. Write
\(q_{1-\alpha}(\nu_1,\nu_2)\) for the lower \(1-\alpha\) quantile of
\(F_{\nu_1,\nu_2}\). Appendix~\ref{app:features} justifies the Gaussian
comparison by an \(O(n_{\min}^{-1/14})\) Kolmogorov bound for the
centered independent-feature statistics, without restricting the spectrum.
The accuracy of the subsequent \(F\) approximation depends on the
mixture weights. In particular, \(\varepsilon_{\rm WS}^o\to0\) if
\[
 \frac{\tr\{(\bm{V}^o)^3\}}{(J^o)^{3/2}}\longrightarrow0.
\]
The bound \(d_2^o\ge(n_{\min}-1)d_1^o\) also implies that the
finite-denominator \(F\) and matched chi-square critical values are
first-order equivalent; see Proposition~\ref{prop:Fchi}.
Section~\ref{sec:originalspectralcheck} gives original-covariance examples
and a spiked counterexample to vanishing matching error.

\subsection{Estimation of the approximation parameters}\label{sec:fulltransfer}
For the observed-data statistic, we must account for estimating the group
means and replace the reference parameters by consistent estimates.
Proposition~\ref{prop:full} gives the exact mean scale
\begin{equation}\label{eq:bfull}
 b^{\rm full}:=\E_0 B_{\rm full}
 =b^o+\sum_i\frac{\tr(\bm{\Gamma}_i)}{n_i(n_i-1)}.
\end{equation}
Relative to the oracle scale \(b^o\), mean estimation adds the second
term in \eqref{eq:bfull}. Theorem~\ref{thm:fulltransfer} shows that this
term is asymptotically negligible under the stated moment and trace
conditions. It also establishes \(B_{\rm full}/b^o\to_{\rm p}1\) and the
same first-order reference for the centered full-sample statistic and
\(U^o\) under \(H_0\). The exact expansion, including third-moment terms
for asymmetric populations, is given in Appendix~\ref{app:unknownmeans}.

\label{sec:fullblocktraces}
Since \(B_{\rm full}\) estimates the common mean scale, estimation of the
degrees of freedom reduces to estimating the squared and cross traces of
\(\bm{\Omega}_i^{\,o}\). We use independent observation differences to remove
\(\bm{\mu}_i\). Divide group \(i\) into \(L_i=\lfloor n_i/4\rfloor\)
disjoint blocks according to a partition fixed independently of the data.
For block \(r\), set
\begin{align}
 \bm{a}_{ir}&=(\bm{X}_{i,4r-3}-\bm{X}_{i,4r-2})/\sqrt2,&
 \bm{b}_{ir}&=(\bm{X}_{i,4r-1}-\bm{X}_{i,4r})/\sqrt2,\notag\\
 \bm{v}_{ir}&=\svec(\bm{a}_{ir}\bm{a}_{ir}^\top-\bm{b}_{ir}\bm{b}_{ir}^\top),&
 \bm{w}_{ir}&=\svec\{(\bm{a}_{ir}\bm{b}_{ir}^\top+\bm{b}_{ir}\bm{a}_{ir}^\top)/\sqrt2\}.
 \label{eq:fourblockfeatures}
\end{align}
Define the signed rank-two matrix
\begin{equation}\label{eq:fourblockK}
 \bm{K}_{ir}=\bm{v}_{ir}\bm{v}_{ir}^\top-\bm{w}_{ir}\bm{w}_{ir}^\top.
\end{equation}
The blocks are independent within and between groups. Direct calculation gives
\begin{equation}\label{eq:fourblockmean}
 \E(\bm{v}_{ir}\bm{v}_{ir}^\top)=\bm{\Omega}_i^{\,o}+\bm{\Gamma}_i,
 \qquad \E(\bm{w}_{ir}\bm{w}_{ir}^\top)=\bm{\Gamma}_i,
 \qquad \E \bm{K}_{ir}=\bm{\Omega}_i^{\,o}.
\end{equation}
The cross-contrast term removes the additional covariance \(\bm{\Gamma}_i\)
introduced by differencing. Consequently, for \(L_i\ge2\), multiplying
independent block estimates and taking traces gives
\begin{align}
 \widetilde s_i&=\frac1{L_i(L_i-1)}
                   \sum_{r\ne s}\tr(\bm{K}_{ir}\bm{K}_{is}),\label{eq:fullshat}\\
 \widetilde c_{12}&=\frac1{L_1L_2}
                   \sum_{r,s}\tr(\bm{K}_{1r}\bm{K}_{2s}).\label{eq:fullchat}
\end{align}
These estimators are unbiased for
\(s_i^o=\tr\{(\bm{\Omega}_i^{\,o})^2\}\) and
\(c_{12}^o=\tr(\bm{\Omega}_1^{\,o}\bm{\Omega}_2^{\,o})\), respectively.
Independence requires eight distinct observations for a within-group
product and four from each group for a cross-group product. Since the unbiased estimates can be
negative, define \(\widehat s_i=(\widetilde s_i)_+\),
\(\widehat c_{12}=(\widetilde c_{12})_+\), and
\begin{align}
 \widehat J_{\rm full}
 &=\frac{\widehat s_1}{n_1^2}
   +\frac{2\widehat c_{12}}{n_1n_2}
   +\frac{\widehat s_2}{n_2^2},&
 \widehat R_{\rm full}
 &=\sum_{i=1}^2\frac{\widehat s_i}{n_i^2(n_i-1)}.
 \label{eq:fullblockJR}
\end{align}
The weights in \eqref{eq:fullblockJR} use the original sample sizes
\(n_i\). An incomplete final block is omitted from trace estimation, but
all observations contribute to \(Q_{\rm full}\) and \(B_{\rm full}\).

\begin{proposition}[Feasible full-sample trace calibration]
\label{prop:fullblocktraces}
Suppose the oracle features satisfy the directional fourth-moment
condition in Section~\ref{sec:oraclereference}, the sample sizes are
balanced and tend to infinity, and \(S_o=s_1^o+s_2^o>0\). If
\begin{equation}\label{eq:blockdominance}
 \tr(\bm{\Gamma}_1^2)+\tr(\bm{\Gamma}_2^2)=O(S_o),
\end{equation}
then
\[
 \E(\widehat J_{\rm full}/J^o-1)^2=O(L_{\min}^{-1}),\qquad
 \E(\widehat R_{\rm full}/R^o-1)^2=O(L_{\min}^{-1}).
\]
In particular, both calibration quantities are relatively consistent and
their probabilities of being zero tend to zero. This result does not require
\(s_i^o\) or \(c_{12}^o\) to be bounded away from zero individually.
\end{proposition}

Under \(H_0\), \eqref{eq:blockdominance} is equivalent in order to
\eqref{eq:GammaOmegacondition}, which controls mean estimation. No sample
separation is needed between the statistic and the trace estimates.

Original-vector inner products suffice to compute the block products;
no \(d\times d\) matrix is needed
(Appendix~\ref{app:blockcomputation}). The cost per partition is
\(O((n_1+n_2)^2p)\). To reduce dependence on a single partition, we
average the raw trace estimates before taking positive parts. For each
fixed row order, apply cyclic shifts \(q=0,1,2,3\) and form consecutive
blocks of four. Every observation then enters a complete block in at
least one partition. This construction requires \(n_1,n_2\ge8\).

The test studied here combines these shifts with \(M=20\) permutation
pairs, generated independently of the data values and then held fixed.
Write \(\widetilde s_i^{(m,q)}\) and
\(\widetilde c_{12}^{(m,q)}\) for the raw estimates from permutation
pair \(m\) and shift \(q\). Set
\begin{equation}\label{eq:pooledtraces}
 \bar s_i=\frac1{4M}\sum_{m=1}^{M}\sum_{q=0}^{3}
              \widetilde s_i^{(m,q)},\qquad
 \bar c_{12}=\frac1{4M}\sum_{m=1}^{M}\sum_{q=0}^{3}
              \widetilde c_{12}^{(m,q)}.
\end{equation}
Substitute \((\bar s_i)_+\) and \((\bar c_{12})_+\) into
\eqref{eq:fullblockJR}, keeping \(Q_{\rm full}\), \(B_{\rm full}\), and
\(U_{\rm full}\) unchanged. For fixed partitions, Jensen's inequality gives
\[
 \E(\bar s_i-s_i^o)^2
 \le \frac1{4M}\sum_{m,q}
          \E(\widetilde s_i^{(m,q)}-s_i^o)^2,
\]
and likewise for the cross trace. Proposition~\ref{prop:fullblocktraces}
therefore applies to the pooled estimates despite dependence between
partitions. The bound alone does not imply a \(1/M\) variance reduction.
Section~\ref{sec:simulation} examines the resulting calibration and
row-order sensitivity.

\subsubsection{The feasible \texorpdfstring{$F$}{F} test}
\label{sec:releasealgorithm}
We now estimate the two degrees of freedom using
\(B_{\rm full}\), \(\widehat J_{\rm full}\), and
\(\widehat R_{\rm full}\). The following theorem gives their ratio
consistency and the resulting null rejection probability.
Its proof is in Appendix~\ref{app:unknownmeans}.

\begin{theorem}[Feasible full-sample calibration]\label{thm:fulltransfer}
Suppose \(H_0\) holds, \(n_1,n_2\ge8\), the sampling and moment
conditions in Section~\ref{sec:oraclereference} hold, and
\eqref{eq:blockdominance} holds. Use the full-sample statistic
\eqref{eq:Ffull} and the four-block estimates \eqref{eq:fullblockJR}, including
the prescribed partition averaging if desired. Define
\[
 \widehat d_{1,\rm full}=B_{\rm full}^2/\widehat J_{\rm full},\qquad
 \widehat d_{2,\rm full}=B_{\rm full}^2/\widehat R_{\rm full}.
\]
Then
\begin{align}
 (U_{\rm full}-U^o)/a^o&=\Op(n_{\min}^{-1/2}),\label{eq:FB1}\\
 B_{\rm full}/b^o&\to_{\rm p}1,\label{eq:FB2}\\
 \widehat d_{j,\rm full}/d_j^o&\to_{\rm p}1,\qquad j=1,2.\label{eq:fullDFconsistent}
\end{align}
If a required statistic or calibration quantity is nonfinite, or any of
\(B_{\rm full},\widehat J_{\rm full},\widehat R_{\rm full}\) is nonpositive,
record calibration failure and nonrejection. Otherwise reject if and only if
\(F_{\rm full}>q_{1-\alpha}(\widehat d_{1,\rm full},\widehat d_{2,\rm full})\);
a valid nonrejection is not a calibration failure. The failure probability
tends to zero, and
\begin{equation}\label{eq:fullsize}
 \left|\Pp_0(\mathrm{reject})-\alpha\right|
 \le\varepsilon_{\rm WS}^{\,o}+o(1).
\end{equation}
If, in addition,
\(\tr\{(\bm{V}^o)^3\}/[\tr\{(\bm{V}^o)^2\}]^{3/2}\to0\), then
\(\Pp_0(\mathrm{reject})\to\alpha\).
\end{theorem}
The test retains \(B_{\rm full}\) and
\(F_{\rm full}=1+U_{\rm full}/B_{\rm full}\); only the trace estimates
are projected onto the nonnegative half-line. This projection preserves
\(\widehat R_{\rm full}\le\widehat J_{\rm full}/(n_{\min}-1)\), so
\(\widehat d_{2,\rm full}\ge(n_{\min}-1)\widehat d_{1,\rm full}\)
whenever calibration is defined.

We denote the procedure with $M=20$ in \eqref{eq:pooledtraces} by
$F_{20}$. The four-shift procedure based on the original row order is
reported as F in the numerical comparisons. For $n_1,n_2\ge8$ and
$0<\alpha<1$, $F_{20}$ is computed as follows.
\begin{enumerate}[label=\arabic*.,leftmargin=*]
\item Center each group and compute $D_{n_i}$ by \eqref{eq:Dgram}, followed
by $Q_{\rm full}$, $B_{\rm full}$ and $U_{\rm full}$ in
\eqref{eq:Bfull}--\eqref{eq:Ffull}. All observations enter these quantities.
\item Fix $M=20$ pairs of permutations independently of the observed values.
For each pair, apply shifts $q=0,1,2,3$ and form consecutive blocks of four.
Compute the three raw estimates in
\eqref{eq:fullshat}--\eqref{eq:fullchat} for each of the 80 partitions.
\item Average the raw estimates as in \eqref{eq:pooledtraces}, then take
componentwise positive parts once. Use the original sample sizes in
\eqref{eq:fullblockJR} and set
$\widehat d_1=B_{\rm full}^2/\widehat J_{\rm full}$ and
$\widehat d_2=B_{\rm full}^2/\widehat R_{\rm full}$.
\item Apply the failure rule in Theorem~\ref{thm:fulltransfer}. On failure,
return an explicit failure flag and no inferential $p$-value. Otherwise,
report $F_{\rm full}$ and
$p_F=1-F_{\widehat d_1,\widehat d_2}(F_{\rm full})$, and reject when
$p_F<\alpha$.
\end{enumerate}
Here $F_{\nu_1,\nu_2}(\cdot)$ denotes the ordinary $F$ distribution
function. The rejection probability satisfies \eqref{eq:fullsize}. Failed calibrations count as nonrejections in the
simulations and have no reported inferential $p$-value.

The partition plan is chosen independently of the data. Finite partition
averaging reduces sensitivity to row order but does not eliminate it.
The computation requires $O(M(n_1+n_2)^2p)$ operations; implementation
details are given in Appendix~\ref{app:simulationcomputation}.

\subsection{Asymptotic null distribution}
\label{sec:unifiedlimit}
The limiting shape of the reference distribution depends on the relative
sizes of the eigenvalues of \(\bm{V}^o\). Let \(\lambda_{nr}\) denote these
eigenvalues in decreasing order, appending zeros beyond the feature
dimension, and put
\[
 c_{nr}=\lambda_{nr}/\sqrt{J^o},\qquad r=1,2,\ldots.
\]
The normalization gives \(\sum_r c_{nr}^2=1\). If each coefficient has
a limit, the Gaussian comparison yields the following null distribution.
The comparison and size bound in Theorem~\ref{thm:fulltransfer} remain
valid without this additional convergence assumption.

\begin{corollary}[Unified full-sample null limit]
\label{cor:unifiedfull}
Under the assumptions of Theorem~\ref{thm:fulltransfer}, suppose also that
\(c_{nr}\to c_r\) for every fixed \(r\). Then \(c_r\ge0\),
\(\sum_r c_r^2\le1\), and
\begin{equation}\label{eq:unifiedfull}
 \frac{b^o}{a^o}(F_{\rm full}-1)\ \Rightarrow\ \zeta,
 \qquad
 \zeta=\sqrt{1-\sum_{r\ge1}c_r^2}\,Z_0
       +\frac1{\sqrt2}\sum_{r\ge1}c_r(Z_r^2-1),
\end{equation}
where \(Z_0,Z_1,\ldots\) are independent standard normal variables.
The series converges in \(L^2\), and \(\zeta\) has mean zero, variance
one, and a continuous distribution function \(G_\zeta\).
The same limit holds for \(U_0^{o,*}/a^o\).
On calibration-failure events the statistic can be assigned any fixed value;
these events have vanishing probability.
\end{corollary}
The limit is standard normal if all \(c_r\) vanish and nonnormal
otherwise; when \(\sum_r c_r^2=1\), it consists entirely of centered
chi-square components. Corollary~\ref{cor:unifiedfull} extends the spectral
representation of \citet[Corollary~1]{WangXu2022} to the full-sample ratio,
in agreement with the covariance normal-reference results of
\citet{WangZhuZhang2024,WangZhuZhang2026}. Ordinary \(F\) calibration
still incurs the matching error in \eqref{eq:fullsize}.

\subsubsection{Sufficient conditions under factor models}
\label{sec:verifiedmodels}

The feature-moment and trace conditions can be verified directly from
an independent-coordinate factor model. Consider
\begin{equation}\label{eq:verifiedfactormodel}
 \bm{\epsilon}_i=\bm{\Sigma}_i^{1/2}\bm{z}_i,
 \qquad \bm{z}_i=(z_{i1},\ldots,z_{ip})^\top,
\end{equation}
where the coordinates of $\bm{z}_i$ are independent, have mean zero and variance
one, and need not have the same distribution. Assume that constants
$M<\infty$ and $c>0$, independent of $p$ and the sample sizes, satisfy
\begin{equation}\label{eq:latentmoments}
 \max_{i,j}\E z_{ij}^8\le M,
 \qquad \min_{i,j}\{\E z_{ij}^4-1\}\ge c.
\end{equation}
The lower bound excludes latent coordinates with constant squared values;
the eigenvalues of $\bm{\Sigma}_i$ may nevertheless approach zero.

By Proposition~\ref{prop:verifiedmodels} in
Appendix~\ref{app:momentconditions}, these assumptions imply both the
directional fourth-moment bound and \eqref{eq:blockdominance}.
In particular, for constants \(0<c_0\le C_0<\infty\) uniform across
groups and dimensions,
\(c_0\bm{\Gamma}_i\preceq\bm{\Omega}_i^o\preceq C_0\bm{\Gamma}_i\).
This permits singular covariance matrices and imposes no restriction on
\(p/n_{\min}\). Gaussian, standardized \(t_{10}\), and
standardized centered chi-square innovations satisfy
\eqref{eq:latentmoments}; their fourth and eighth moments are reported
in Table~\ref{tab:modelmoments} in the appendix. The groups may have
different innovation distributions under \(H_0\).

\subsubsection{Accuracy of the \texorpdfstring{$F$}{F} approximation}
\label{sec:originalspectralcheck}

We next express a sufficient condition for vanishing matching error in
terms of the original covariance matrix. Under $H_0$, write
$\bm{\Sigma}_1=\bm{\Sigma}_2=\bm{\Sigma}$ and assume that the constants in
\eqref{eq:Omegasandwich} are uniform across groups. Then it suffices that
\begin{equation}\label{eq:originaldiffuse}
 \frac{\|\bm{\Sigma}\|_{\mathrm{op}}^2}{\tr(\bm{\Sigma}^2)}\longrightarrow0.
\end{equation}
Indeed, if $w_n=n_1^{-1}+n_2^{-1}$, then
$c_0w_n\bm{\Gamma}\preceq \bm{V}^o\preceq C_0w_n\bm{\Gamma}$, and hence
\begin{align}
 \frac{\tr\{(\bm{V}^o)^3\}}{\tr\{(\bm{V}^o)^2\}^{3/2}}
 &\le \frac{\|\bm{V}^o\|_{\mathrm{op}}}{\sqrt{\tr\{(\bm{V}^o)^2\}}}
 \le \frac{C_0}{c_0}
       \frac{\|\bm{\Gamma}\|_{\mathrm{op}}}{\sqrt{\tr(\bm{\Gamma}^2)}}\notag\\
 &\le \sqrt2\frac{C_0}{c_0}
       \frac{\|\bm{\Sigma}\|_{\mathrm{op}}^2}{\tr(\bm{\Sigma}^2)}.
 \label{eq:originaldiffusebound}
\end{align}
Here $\|\bm{\Gamma}\|_{\mathrm{op}}=2\|\bm{\Sigma}\|_{\mathrm{op}}^2$.
The same conclusion holds for the paired features in
Appendix~\ref{sec:independent}, with their applicable comparison constants.
Identity covariance and autoregressive covariance
$(\rho^{|j-k|})_{j,k}$ with fixed $|\rho|<1$ satisfy
\eqref{eq:originaldiffuse}: the operator norms are bounded and
$\tr(\bm{\Sigma}^2)\ge p\to\infty$.

The matching error behaves differently when one covariance eigenvalue
dominates. To see this, consider Gaussian observations with
$\bm{\Sigma}=(1-\eta)\bm{I}_p+\eta p\bm{v}\bm{v}^\top$, where $0<\eta<1$ is fixed and
$\|\bm{v}\|=1$, as in Table~\ref{tab:simulation_covariance}.
Its eigenvalues are $\lambda_1=\eta p+1-\eta$ and
$\lambda_2=\cdots=\lambda_p=1-\eta$. The largest eigenvalue of $\bm{\Gamma}$
is $2\lambda_1^2$, while
\[
 p^{-2}\tr(\bm{\Gamma})\longrightarrow1+\eta^2,
 \qquad
 p^{-4}\tr(\bm{\Gamma}^2)\longrightarrow4\eta^4.
\]
The other feature eigenvalues contribute to the limiting mean but
negligibly to the variance at this scale. Consequently,
\begin{equation}\label{eq:spikedQlimit}
 \frac{Q^{o,*}}{b^o}\ \Rightarrow\
 \frac{1-\eta^2+2\eta^2\chi_1^2}{1+\eta^2},
 \qquad
 d_1^o\longrightarrow\frac{(1+\eta^2)^2}{4\eta^4}.
\end{equation}
As $n_{\min}\to\infty$, $B^{o,*}/b^o\to_{\rm p}1$ and $d_2^o\to\infty$.
The limit on the right of \eqref{eq:spikedQlimit} has a strictly positive
lower endpoint, whereas the limiting law of the matched distribution
$\chi_{d_1^o}^2/d_1^o$ has support starting at zero. Thus their Kolmogorov
distance cannot vanish.
More explicitly, put \(\ell_\eta=(1-\eta^2)/(1+\eta^2)>0\) and
\(d_\eta=(1+\eta^2)^2/(4\eta^4)\). At the lower endpoint of the
limiting reference law,
\[
 \liminf_{p,n_{\min}\to\infty}\varepsilon_{\rm WS}^o
 \ge F_{\chi^2_{d_\eta}}(d_\eta\ell_\eta)>0.
\]
Combining \eqref{eq:spikedQlimit} with Gaussian comparison and relative
parameter consistency gives the limiting rejection probability of the
feasible test:
\begin{equation}\label{eq:spikedsize}
 \alpha_\infty(\eta)
 =1-F_{\chi_1^2}\left[
  1+\frac{1+\eta^2}{2\eta^2}
      \left\{\frac{\chi^2_{d_\eta,1-\alpha}}{d_\eta}-1\right\}
 \right].
\end{equation}
For \(\alpha=0.05\), this equals approximately \(0.06727\),
\(0.06143\), and \(0.05488\) for \(\eta=0.1,0.3,0.5\), respectively.
These limits take both \(p\) and \(n_{\min}\) to infinity with balanced
group sample sizes. They quantify the persistent calibration error in
the spiked simulation settings.

\begin{corollary}[Asymptotic level for verified original-data models]
\label{cor:verifiedlevel}
Under \(H_0\), suppose \eqref{eq:verifiedfactormodel}--\eqref{eq:latentmoments}
hold, \(\bm{\Sigma}\ne0\), and the sample sizes are balanced and tend to infinity.
Then all moment, bridge, denominator, and calibration-estimation conditions
of Theorem~\ref{thm:fulltransfer} hold. If also \eqref{eq:originaldiffuse}
holds, the fully specified test has rejection probability tending to
\(\alpha\). In particular, this conclusion applies to the Gaussian,
standardized \(t_{10}\), and standardized centered chi-square innovation
models with identity covariance or fixed-parameter autoregressive covariance
as \(p\to\infty\), without an additional restriction on \(p/n_{\min}\).
\end{corollary}

\subsection{Local power and consistency}
\label{sec:fullpower}
We study power by comparing the squared covariance distance with the standard
deviation of the centered statistic. Along a sequence of alternatives,
let \(\bm{\Delta}=\bm{\theta}_1-\bm{\theta}_2\) and
\(\delta_n=\norm{\bm{\Delta}}^2=\norm{\bm{\Sigma}_1-\bm{\Sigma}_2}_F^2\), and evaluate
all oracle quantities at the current populations. Under the conditions
of Lemma~\ref{lem:alternativebridge}, the unknown-mean expansion and
ratio consistency continue to hold. For the central Gaussian reference,
define
\[
 G_n^o(x)=\Pp(U_0^{o,*}/a^o\le x),\qquad
 t_n^o=\frac{b^o}{a^o}\{q_{1-\alpha}(d_1^o,d_2^o)-1\}.
\]
Here \(\bm{u}_{ij}^{o,*}\sim N_d(0,\bm{\Omega}_i^{\,o})\) independently, and
\begin{equation}\label{eq:centraloraclereference}
 Q_0^{o,*}=\norm{\bm{\bar u}_1^{o,*}-\bm{\bar u}_2^{o,*}}^2,
 \qquad U_0^{o,*}=Q_0^{o,*}-B^{o,*}.
\end{equation}
The reference therefore has zero feature means even when \(\bm{\Delta}\ne0\),
with covariances determined by the alternative.

\begin{theorem}[Full-sample local power]
\label{thm:fullpower}
Suppose the sampling and moment conditions in
Section~\ref{sec:oraclereference} and \eqref{eq:blockdominance} hold
uniformly along a sequence of alternatives. If
\begin{equation}\label{eq:fulllocal}
 h_n=\delta_n/a^o=O(1),\qquad
 \bm{\Delta}^\top \bm{V}^o\bm{\Delta}=o\{(a^o)^2\},
\end{equation}
then
\begin{equation}\label{eq:fullpowersequence}
 \Pp(\mathrm{reject})=1-G_n^o(t_n^o-h_n)+o(1).
\end{equation}
If also \(c_{nr}\to c_r\), \(h_n\to h\), and
\(d_1^o\to d\in[1,\infty]\), then
\begin{equation}\label{eq:fullpowerlimit}
 \Pp(\mathrm{reject})\longrightarrow
  1-G_\zeta\{t_\alpha(d)-h\},
 \quad
 t_\alpha(d)=
 \begin{cases}
  \sqrt{d/2}\{\chi^2_{d,1-\alpha}/d-1\},&d<\infty,\\
  z_{1-\alpha},&d=\infty.
 \end{cases}
\end{equation}
In particular, if
\(\tr\{(\bm{V}^o)^3\}/(J^o)^{3/2}\to0\) and \(h_n\to h\), then
\begin{equation}\label{eq:fullnormalpower}
 \Pp(\mathrm{reject})\longrightarrow1-\Phi(z_{1-\alpha}-h).
\end{equation}
\end{theorem}
The first condition in \eqref{eq:fulllocal} places the signal on the
standard-deviation scale. The second makes the linear signal term
negligible, as in the covariance normal-reference power analyses of
\citet{WangZhuZhang2024,WangZhuZhang2026}. Together they yield
\eqref{eq:fullpowersequence} without requiring spectral convergence.
With a limiting spectrum, \eqref{eq:fullpowerlimit} gives rejection
probability \(1-G_\zeta\{t_\alpha(d)\}\) at \(h=0\); this can differ from
\(\alpha\) if matching error persists. A normal limit requires control
of the squared eigenvalue weights, beyond divergence of the
first-moment effective rank \(d_1^o\).

The linear-term condition has a simple sufficient criterion. Since
\(a^o\asymp\sqrt{J^o}\),
\begin{equation}\label{eq:localspectralbound}
 \frac{\bm{\Delta}^\top \bm{V}^o\bm{\Delta}}{(a^o)^2}
 \le \frac{\delta_n}{a^o}\,
       \frac{\|\bm{V}^o\|_{\rm op}}{a^o}=o(1)
\end{equation}
whenever \(\delta_n=O(a^o)\) and
\(\|\bm{V}^o\|_{\rm op}/\sqrt{J^o}\to0\).
For Gaussian populations with uniformly bounded covariance eigenvalues
bounded away from zero and balanced sample sizes of order \(n\),
\(a^o\asymp p/n\) and \(\|\bm{V}^o\|_{\rm op}=O(n^{-1})\).
Consequently, alternatives with
\(\|\bm{\Sigma}_1-\bm{\Sigma}_2\|_F^2=O(p/n)\) satisfy the condition as
\(p\to\infty\). This includes local perturbations of identity covariance
that preserve these eigenvalue bounds.

For spike-aligned signals, however, the linear term can remain appreciable.
The Gaussian example in Appendix~\ref{app:powerexample} has bounded
\(\delta_n/a^o\), while
\(\bm{\Delta}^\top \bm{V}^o\bm{\Delta}/(a^o)^2\) converges to a positive constant.
Thus the two conditions in \eqref{eq:fulllocal} must be checked separately.

\begin{proposition}[Full-sample consistency]
\label{prop:fullconsistency}
Suppose the sampling and moment conditions in
Section~\ref{sec:oraclereference} and \eqref{eq:blockdominance} hold
uniformly along a sequence of alternatives. If
\(\delta_n>0\) eventually and
\begin{equation}\label{eq:fullconsistency}
 \frac{(a^o)^2+4\bm{\Delta}^\top \bm{V}^o\bm{\Delta}}{\delta_n^2}\longrightarrow0,
\end{equation}
then the rejection probability of the feasible full-sample test tends to one.
No spectral diffuseness assumption is required for this assertion.
\end{proposition}
Under balanced sample sizes, \(a^o\) is of order
\(\sqrt{S_o}/n_{\min}\), where
\(S_o=\tr\{(\bm{\Omega}_1^{\,o})^2\}+\tr\{(\bm{\Omega}_2^{\,o})^2\}\).
Thus \(\norm{\bm{\Sigma}_1-\bm{\Sigma}_2}_F^2\) of order
\(\sqrt{S_o}/n_{\min}\) is the local signal scale in
Theorem~\ref{thm:fullpower}, subject to its linear-term condition.
Establishing a minimax detection boundary would also require a lower
bound over a class of tests.

\section{Simulation studies}\label{sec:simulation}
\subsection{Design and competing procedures}
In this section, we compare the size and power of $F_{20}$ with those of
three competing procedures: its chi-square counterpart
$\mathrm{Chi}_{20}$, the Li--Chen test (LC), and a three-cumulant
normal-reference test (NR3). The nominal level is $\alpha=0.05$.
$F_{20}$ and $\mathrm{Chi}_{20}$ use the same statistic and pooled trace
estimates, differing only in their reference distributions. LC uses a
normal approximation to $U_{\mathrm{full}}/
(2D_{n_1}/n_2+2D_{n_2}/n_1)$. The implementation of $F_{20}$ is specified
in Section~\ref{sec:releasealgorithm}.

For NR3, we use the statistic and calibration in Sections 2.1 and 2.3
of arXiv:2212.12338v1, the preprint underlying \citet{WangZhuZhang2024}.
As suggested there, group sample means are subtracted before constructing
the outer-product features. This statistic differs from
$U_{\mathrm{full}}$. The implementation was checked against the preprint
equations and an independent evaluation of the HDNRA three-cumulant
mean-test code applied to the same outer-product features. The complete
published covariance-test implementation was not available for comparison.
Invalid calibrations are treated as nonrejections and their frequencies
are reported separately. Details of the numerical checks are given in
Appendix~\ref{app:simulationcomputation}.

The two independent samples are generated from
\[
 \bm{X}_{1j}=\bm{\Sigma}_1^{1/2}\bm{z}_{1j},\qquad
 \bm{X}_{2j}=\tfrac12\bm{1}_p+\bm{\Sigma}_2^{1/2}\bm{z}_{2j},
\]
where $\bm{\Sigma}_i^{1/2}$ is the symmetric positive square root.
The coordinates of $\bm{z}_{ij}$ are independent. For the main comparison,
both groups use standard normal innovations (G),
$\sqrt{8/10}\,t_{10}$ innovations (T), or
$(\chi_1^2-1)/\sqrt2$ innovations (K). These choices represent normal,
symmetric nonnormal, and skewed observations, respectively, while
retaining finite eighth moments. Both group means are estimated.
Additional comparisons with different innovation laws in the two groups
are described in Section~\ref{sec:v13extensions}.

For size, we set $\bm{\Sigma}_1=\bm{\Sigma}_2=\bm{\Sigma}_0$ and consider
\[
 \mathrm I:\ \bm{I}_p,\qquad
 \mathrm A:\ (0.6^{|r-s|})_{r,s=1}^p,\qquad
 \mathrm S:\ 0.7\bm{I}_p+0.3\bm{1}_p\bm{1}_p^\top.
\]
The first two structures have diffuse spectra, whereas the leading
eigenvalue of S remains nonnegligible relative to its Frobenius norm as
$p$ increases. We take $p\in\{50,200\}$ and
$(n_1,n_2)\in\{(20,30),(50,80)\}$. Each of the resulting 36 settings uses
$10{,}000$ independent datasets. The four methods are applied to the
same observations in each run.

\subsection{Empirical size and approximation parameters}
Tables~\ref{tab:newsize20} and \ref{tab:newsize50} display the empirical
sizes for the two sample-size pairs. Overall size accuracy is summarized
by the average relative error
\[
 \mathrm{ARE}=\frac{100}{36}\sum_{k=1}^{36}
       \frac{|\widehat\pi_k-0.05|}{0.05},
\]
where $\widehat\pi_k$ denotes the empirical size in setting $k$.
The ARE values of $F_{20}$, $\mathrm{Chi}_{20}$, LC, and NR3 are
29.21, 32.49, 62.31, and 63.86\%, respectively.
Their empirical size ranges are 4.99--9.13\%, 5.11--9.31\%,
4.69--15.29\%, and 0.00--5.11\%, respectively. At a rejection probability of 5\%,
the Monte Carlo standard error is 0.218 percentage points.

It is seen that $F_{20}$ behaves differently under the three covariance
structures. For I and A, its sizes at the larger sample-size pair are
4.99--5.64\% across the three innovation models and two
dimensions. Under S, it remains liberal. For G with $p=50$, increasing
$(n_1,n_2)$ from $(20,30)$ to $(50,80)$ changes its size from
8.48\% to 6.60\%. Under G with I and $p=200$, its size
at $(50,80)$ is 5.35\%. Thus, increasing sample size improves
estimation but does not by itself remove the difficulty caused by a
persistent spike.

LC is particularly liberal for skewed observations. Under K with I,
$p=50$, and $(n_1,n_2)=(50,80)$, its size is 15.29\%, compared with
5.23\% for $F_{20}$. NR3 is conservative in several diffuse-spectrum
settings. Its largest calibration failure rate is 14.49\%,
compared with 0.60\% for the two pooled procedures; LC has
no failures. The larger size errors of LC and NR3 therefore arise from
different departures from the nominal level. For Gaussian observations
under I, the known-mean diagnostic in Appendix~\ref{app:simulationcomputation}
shows upward finite-sample bias in NR3's second-cumulant estimate,
which helps explain its conservative sizes.

The estimated degrees of freedom help describe the reference
distributions. Under G with $p=50$ and $(n_1,n_2)=(50,80)$, the median
numerator degrees of freedom of $F_{20}$ are 1328.2, 301.5, and
36.3 for I, A, and S, respectively. The corresponding NR3 medians
are 152.9, 36.7, and 1.75. Within each approximation,
the decrease under S indicates greater skewness. These parameters refer
to different statistics and are not directly interchangeable. The
persistent size error of $F_{20}$ under S agrees with the
moment-matching analysis in Section~\ref{sec:originalspectralcheck}.
The sizes of $F_{20}$ and $\mathrm{Chi}_{20}$ remain close, consistent
with their first-order equivalence.

\begin{table}[!ht]
\centering\small
\setlength{\tabcolsep}{4pt}
\caption{Empirical sizes (percent) at $(n_1,n_2)=(20,30)$; $N=10{,}000$ per setting. The last two columns report failure percentages. G, T, and K denote normal, standardized $t_{10}$, and standardized $\chi^2_1$ innovations.}
\label{tab:newsize20}
\begin{tabular*}{\textwidth}{@{\extracolsep{\fill}}ccrrrrrrr@{}}
\toprule
Dist. & Cov. & $p$ & $F_{20}$ & $\mathrm{Chi}_{20}$ & LC & NR3 & $F_{20}$ fail & NR3 fail \\
\midrule
G & I & 50 & 5.76 & 5.98 & 4.94 & 0.72 & 0.00 & 0.20 \\
G & I & 200 & 6.45 & 6.68 & 5.45 & 2.39 & 0.00 & 14.49 \\
G & A & 50 & 6.29 & 6.49 & 5.93 & 1.55 & 0.00 & 0.00 \\
G & A & 200 & 6.18 & 6.37 & 5.26 & 0.78 & 0.00 & 1.98 \\
G & S & 50 & 8.48 & 8.66 & 9.12 & 4.86 & 0.26 & 0.00 \\
G & S & 200 & 8.98 & 9.31 & 10.03 & 5.09 & 0.56 & 0.00 \\
T & I & 50 & 5.68 & 5.94 & 5.28 & 0.16 & 0.00 & 0.09 \\
T & I & 200 & 5.76 & 6.08 & 4.69 & 0.26 & 0.00 & 7.04 \\
T & A & 50 & 6.33 & 6.66 & 6.45 & 1.03 & 0.00 & 0.00 \\
T & A & 200 & 6.18 & 6.51 & 5.51 & 0.29 & 0.00 & 1.45 \\
T & S & 50 & 8.34 & 8.51 & 9.24 & 4.59 & 0.26 & 0.00 \\
T & S & 200 & 9.03 & 9.22 & 9.67 & 5.11 & 0.60 & 0.00 \\
K & I & 50 & 5.93 & 6.32 & 14.80 & 0.00 & 0.00 & 0.00 \\
K & I & 200 & 6.26 & 6.47 & 8.45 & 0.00 & 0.00 & 0.22 \\
K & A & 50 & 6.20 & 6.44 & 10.52 & 0.04 & 0.00 & 0.00 \\
K & A & 200 & 5.97 & 6.19 & 6.72 & 0.00 & 0.00 & 0.03 \\
K & S & 50 & 8.14 & 8.32 & 10.42 & 2.48 & 0.27 & 0.00 \\
K & S & 200 & 9.13 & 9.29 & 9.77 & 3.95 & 0.60 & 0.00 \\
\bottomrule
\end{tabular*}
\end{table}
\begin{table}[!ht]
\centering\small
\setlength{\tabcolsep}{4pt}
\caption{Empirical sizes (percent) at $(n_1,n_2)=(50,80)$; $N=10{,}000$ per setting. The last two columns report failure percentages. G, T, and K denote normal, standardized $t_{10}$, and standardized $\chi^2_1$ innovations.}
\label{tab:newsize50}
\begin{tabular*}{\textwidth}{@{\extracolsep{\fill}}ccrrrrrrr@{}}
\toprule
Dist. & Cov. & $p$ & $F_{20}$ & $\mathrm{Chi}_{20}$ & LC & NR3 & $F_{20}$ fail & NR3 fail \\
\midrule
G & I & 50 & 4.99 & 5.11 & 5.17 & 0.72 & 0.00 & 0.00 \\
G & I & 200 & 5.35 & 5.48 & 5.32 & 0.22 & 0.00 & 3.55 \\
G & A & 50 & 5.53 & 5.66 & 6.19 & 2.09 & 0.00 & 0.00 \\
G & A & 200 & 5.39 & 5.45 & 5.35 & 0.47 & 0.00 & 0.00 \\
G & S & 50 & 6.60 & 6.72 & 9.29 & 4.44 & 0.00 & 0.00 \\
G & S & 200 & 7.39 & 7.49 & 10.34 & 4.81 & 0.00 & 0.00 \\
T & I & 50 & 5.17 & 5.24 & 5.74 & 0.36 & 0.00 & 0.00 \\
T & I & 200 & 5.39 & 5.46 & 5.42 & 0.00 & 0.00 & 0.66 \\
T & A & 50 & 5.61 & 5.69 & 6.74 & 1.61 & 0.00 & 0.00 \\
T & A & 200 & 5.57 & 5.68 & 5.73 & 0.19 & 0.00 & 0.00 \\
T & S & 50 & 7.23 & 7.31 & 9.78 & 4.78 & 0.00 & 0.00 \\
T & S & 200 & 7.39 & 7.46 & 10.23 & 4.80 & 0.00 & 0.00 \\
K & I & 50 & 5.23 & 5.33 & 15.29 & 0.00 & 0.00 & 0.00 \\
K & I & 200 & 5.20 & 5.28 & 8.62 & 0.00 & 0.00 & 0.00 \\
K & A & 50 & 5.64 & 5.74 & 11.93 & 0.12 & 0.00 & 0.00 \\
K & A & 200 & 5.47 & 5.55 & 7.12 & 0.00 & 0.00 & 0.00 \\
K & S & 50 & 6.72 & 6.76 & 10.60 & 2.88 & 0.00 & 0.00 \\
K & S & 200 & 7.59 & 7.64 & 10.31 & 4.67 & 0.00 & 0.00 \\
\bottomrule
\end{tabular*}
\end{table}

\subsection{Raw and size-adjusted power}
We first compare power under three forms of covariance difference,
taking $p=50$, $(n_1,n_2)=(50,80)$, G or K innovations, and covariance
I or S. With $\bm{\Sigma}_1=\bm{\Sigma}_0$, the alternatives are
\[
\begin{array}{ll}
\text{dense:}&\bm{\Sigma}_2=(1+\delta)\bm{\Sigma}_0,
       \quad\delta\in\{0.15,0.30,0.50\},\\
\text{sparse:}&\bm{\Sigma}_2=\bm{\Sigma}_0+\delta\diag(\bm{1}_5,\bm{0}_{p-5}),
       \quad\delta\in\{0.5,1,2\},\\
\text{rank one:}&\bm{\Sigma}_2=\bm{\Sigma}_0+
       \delta\lambda_{\max}(\bm{\Sigma}_0)\bm{u}\bm{u}^\top,
       \quad\delta\in\{0.5,1,2\},
\end{array}
\]
where $\bm{u}$ is a leading unit eigenvector: a fixed coordinate
direction for I and $\bm{1}_p/\sqrt p$ for S. The parameter $\delta$
controls the covariance difference within each family. Signal strengths
are not matched across families or covariance structures, so comparisons
concern methods under the same setting. Each of the 36 settings uses
$10{,}000$ independent datasets.

The size differences must be taken into account when interpreting raw
power. We therefore also use $10{,}000$ independent null datasets for
each innovation--covariance combination to estimate the fifth percentile
of each method's approximate $p$-values. This percentile replaces 0.05
when computing size-adjusted power; calibration failures retain the
value one. On the independent null samples used in the size study, the
adjusted rejection rates for the four calibrated models are
4.36--5.08\%. The adjustment uses the known simulation
model and permits comparisons at similar null rejection rates.
The maximum conditional Monte Carlo standard error for power is
0.50 percentage points, excluding uncertainty in the estimated thresholds.

Tables~\ref{tab:newpowerI} and \ref{tab:newpowerC} report raw and
size-adjusted powers. Under K with I and the dense alternative at
$\delta=0.30$, the raw powers of LC and $F_{20}$ are
33.61\% and 19.86\%, whereas their adjusted
powers are 14.46\% and 19.34\%. The higher
raw rejection rate of LC must therefore be read together with its more
liberal size. For the same alternative under G, NR3 has raw power
13.23\%, compared with 28.60\% for $F_{20}$,
but the adjusted powers are 31.93\% and
27.94\%, respectively.

Under K with S and the rank-one alternative at $\delta=1$, the adjusted
powers of NR3 and $F_{20}$ are 77.63\% and 73.17\%.
These comparisons do not establish a uniform power advantage for
$F_{20}$. The form of the covariance difference also matters. Under S,
changing only five marginal variances produces little power at moderate
signals: for G and $\delta=1$, the four adjusted powers are
6.02--6.64\%. Such a change is small relative to the variation
associated with the leading covariance direction. $F_{20}$ and
$\mathrm{Chi}_{20}$ have similar adjusted powers across these settings.

\begin{table}[!ht]
\centering\small
\setlength{\tabcolsep}{4pt}
\caption{Raw/size-adjusted power (percent), covariance I, $p=50$, $(n_1,n_2)=(50,80)$, and $N=10{,}000$. Adjusted thresholds use a separate set of $10{,}000$ null datasets for each innovation model and covariance structure.}
\label{tab:newpowerI}
\begin{tabular*}{\textwidth}{@{\extracolsep{\fill}}clrrrrr@{}}
\toprule
Dist. & Alternative & $\delta$ & $F_{20}$ & $\mathrm{Chi}_{20}$ & LC & NR3 \\
\midrule
G & dense & 0.15 & 9.4/9.1 & 9.5/9.1 & 7.7/7.4 & 2.6/10.3 \\
G & dense & 0.3 & 28.6/27.9 & 29.0/27.9 & 22.7/21.8 & 13.2/31.9 \\
G & dense & 0.5 & 79.0/78.5 & 79.3/78.5 & 71.8/71.0 & 61.3/82.1 \\
G & sparse & 0.5 & 11.0/10.7 & 11.2/10.7 & 10.3/9.9 & 3.1/11.2 \\
G & sparse & 1 & 39.6/38.6 & 39.8/38.6 & 36.0/35.1 & 18.2/40.3 \\
G & sparse & 2 & 99.0/98.9 & 99.0/98.9 & 98.6/98.5 & 96.0/99.1 \\
G & rank & 0.5 & 6.1/5.9 & 6.2/5.9 & 6.1/5.8 & 1.4/5.8 \\
G & rank & 1 & 9.4/9.2 & 9.6/9.2 & 9.2/8.7 & 2.3/9.3 \\
G & rank & 2 & 32.6/31.9 & 32.8/31.8 & 31.1/30.4 & 14.4/32.1 \\
K & dense & 0.15 & 8.8/8.5 & 8.9/8.5 & 19.6/6.7 & 0.0/8.8 \\
K & dense & 0.3 & 19.9/19.3 & 20.0/19.3 & 33.6/14.5 & 0.1/20.6 \\
K & dense & 0.5 & 45.6/44.8 & 45.7/44.7 & 61.5/35.9 & 1.3/46.1 \\
K & sparse & 0.5 & 9.0/8.6 & 9.1/8.6 & 21.0/8.0 & 0.0/8.8 \\
K & sparse & 1 & 23.3/22.7 & 23.6/22.7 & 40.1/20.9 & 0.5/23.7 \\
K & sparse & 2 & 69.4/69.0 & 69.7/68.9 & 82.9/68.3 & 11.1/70.9 \\
K & rank & 0.5 & 6.1/5.8 & 6.3/5.8 & 16.9/5.8 & 0.0/6.0 \\
K & rank & 1 & 8.7/8.3 & 8.8/8.3 & 20.9/8.1 & 0.0/8.6 \\
K & rank & 2 & 19.9/19.4 & 20.1/19.4 & 35.5/19.1 & 0.5/20.1 \\
\bottomrule
\end{tabular*}
\end{table}
\begin{table}[!ht]
\centering\small
\setlength{\tabcolsep}{4pt}
\caption{Raw/size-adjusted power (percent), covariance S, $p=50$, $(n_1,n_2)=(50,80)$, and $N=10{,}000$. Adjusted thresholds use a separate set of $10{,}000$ null datasets for each innovation model and covariance structure.}
\label{tab:newpowerC}
\begin{tabular*}{\textwidth}{@{\extracolsep{\fill}}clrrrrr@{}}
\toprule
Dist. & Alternative & $\delta$ & $F_{20}$ & $\mathrm{Chi}_{20}$ & LC & NR3 \\
\midrule
G & dense & 0.15 & 12.4/9.4 & 12.5/9.4 & 14.5/7.8 & 9.7/10.1 \\
G & dense & 0.3 & 24.2/19.6 & 24.3/19.6 & 27.0/16.7 & 20.2/21.0 \\
G & dense & 0.5 & 44.7/39.0 & 44.8/39.0 & 48.0/34.2 & 39.9/40.7 \\
G & sparse & 0.5 & 7.8/5.4 & 7.9/5.4 & 10.4/5.2 & 5.0/5.4 \\
G & sparse & 1 & 9.2/6.5 & 9.2/6.5 & 11.8/6.0 & 6.3/6.6 \\
G & sparse & 2 & 18.4/12.9 & 18.5/12.9 & 23.7/10.3 & 12.2/12.9 \\
G & rank & 0.5 & 41.5/36.1 & 41.6/36.0 & 45.3/33.0 & 37.5/38.4 \\
G & rank & 1 & 81.4/76.7 & 81.5/76.6 & 84.9/75.4 & 78.8/79.5 \\
G & rank & 2 & 99.3/98.8 & 99.3/98.8 & 99.6/99.0 & 99.2/99.2 \\
K & dense & 0.15 & 12.0/8.9 & 12.1/8.9 & 15.5/7.6 & 6.7/9.7 \\
K & dense & 0.3 & 22.6/17.8 & 22.7/17.8 & 26.7/15.5 & 14.5/18.8 \\
K & dense & 0.5 & 42.4/35.9 & 42.5/35.8 & 47.6/32.7 & 30.8/37.6 \\
K & sparse & 0.5 & 7.6/4.8 & 7.7/4.8 & 11.5/5.1 & 3.1/4.8 \\
K & sparse & 1 & 8.9/6.3 & 9.0/6.3 & 13.3/6.1 & 3.9/6.0 \\
K & sparse & 2 & 18.6/12.6 & 18.7/12.6 & 27.8/12.2 & 6.9/10.9 \\
K & rank & 0.5 & 39.5/33.4 & 39.7/33.4 & 45.0/30.8 & 31.1/37.3 \\
K & rank & 1 & 78.6/73.2 & 78.8/73.2 & 83.2/72.3 & 72.3/77.6 \\
K & rank & 2 & 98.6/97.7 & 98.6/97.7 & 99.2/98.0 & 97.8/98.6 \\
\bottomrule
\end{tabular*}
\end{table}

\subsection{Additional alternatives and mixed innovation distributions}
\label{sec:v13extensions}
To examine the effect of dimension and sample size, we repeat selected
power comparisons under I at $(p,n_1,n_2)=(200,50,80)$ and $(50,80,120)$,
using G and K innovations. The selected alternatives are dense with
$\delta=0.30$, sparse with $\delta=1$, and rank one with $\delta=2$.
For example, under G and the dense alternative, the adjusted power
of $F_{20}$ is 27.94\% at $(50,50,80)$,
29.38\% at $(200,50,80)$, and
48.64\% at $(50,80,120)$.

We also consider changes in five off-diagonal correlations while keeping
all marginal variances equal to one. These use five disjoint $2\times2$
blocks with off-diagonal entries $\delta$, where $\delta=0.30$ or 0.60.
The resulting covariance matrix is positive definite, with nonunit
eigenvalues $1-\delta$ and $1+\delta$. Under G and $\delta=0.60$, the
adjusted power of $F_{20}$ is 28.66\% at $p=50$ and
8.12\% at $p=200$, with $(n_1,n_2)=(50,80)$ in both cases.

Finally, under covariance equality we generate Gaussian innovations
in group 1 and standardized $\chi^2_1$ innovations in group 2. At
$p=200$, the sizes of $F_{20}$ for $(n_1,n_2)=(20,30)$ and $(50,80)$
are 5.87\% and 5.30\% under I, and
9.20\% and 6.97\% under S. This comparison allows
different fourth-order structures in the two groups while keeping
their covariance matrices equal. All additional size and power
settings use $10{,}000$ datasets. Full designs, comparisons with the
other methods, and independent null-calibration checks appear in
Appendix~\ref{app:v13extensions}.

\subsection{Local power and concentrated spectra}
\label{sec:theorysimulation}
We next examine the local-power approximation and the spiked-null calculation
along the path $(p,n_1,n_2)=(50,25,40),(100,50,80),(200,100,160)$.
For Gaussian observations with $\bm{\Sigma}_1=\bm{I}_p$ and
$\bm{\Sigma}_2=(1+t)\bm{I}_p$, we choose $t$ so that
$\delta_n/a^o=h$, for $h=0,1,2$. The parameters and the linear-term
ratio are given in Appendix~\ref{app:theoryexperiments}.
Each original-data setting uses 10,000 independent datasets.

\begin{figure}[!htbp]
\centering
\includegraphics[width=\textwidth]{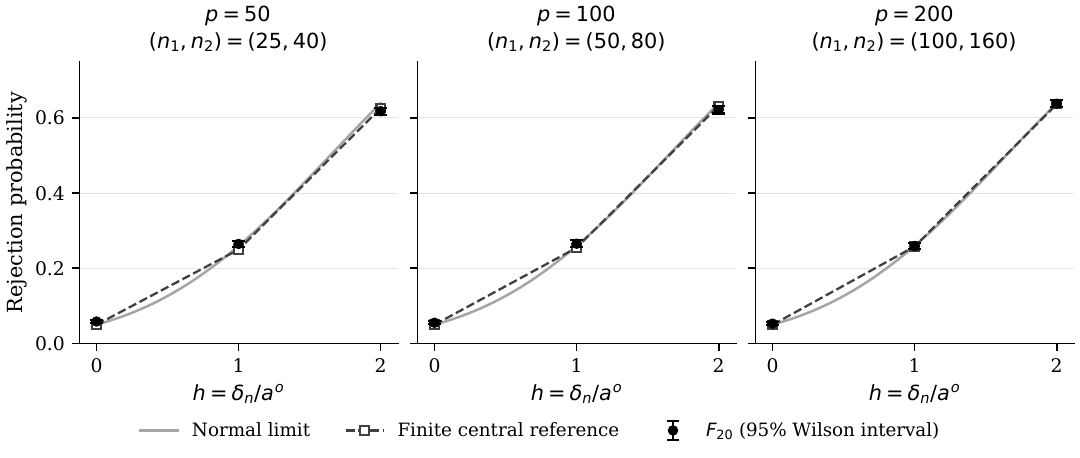}
\caption{Gaussian local alternatives along the displayed dimension and
sample-size path. The intervals for $F_{20}$ are 95\% Wilson intervals
from 10,000 datasets per setting. The finite central-reference prediction
uses 1,000,000 independent reference draws per setting; it approximates
power after omitting the linear signal term.}
\label{fig:localpower}
\end{figure}

Figure~\ref{fig:localpower} compares the empirical rejection probabilities
with the finite central-reference prediction and the normal-limit power
$1-\Phi(z_{0.95}-h)$. At $p=200$, the empirical powers for $h=1$ and
$h=2$ are 25.91\% and 63.78\%, compared with normal-limit predictions
25.95\% and 63.88\%. At $p=50$, the power for $h=2$ is 61.79\%,
and the linear-term ratio is larger. These finite-sample differences
are relevant when using the local-power approximation.

\begin{figure}[!htbp]
\centering
\includegraphics[width=.83\textwidth]{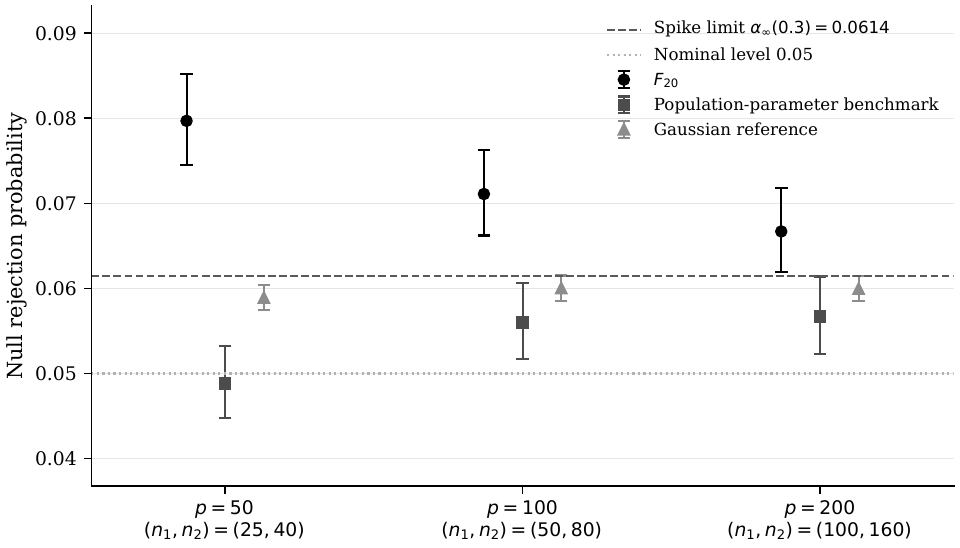}
\caption{Gaussian spiked null with $\eta=0.3$. The original-data procedures
use 10,000 datasets per setting; the Gaussian reference uses 100,000
independent draws. Error bars are 95\% Wilson intervals. The horizontal
lines mark the nominal level and the nonnominal limiting rejection probability.}
\label{fig:spikecalibration}
\end{figure}

For the spiked null, we use
$\bm{\Sigma}_1=\bm{\Sigma}_2=0.7\bm{I}_p+0.3\bm{1}_p\bm{1}_p^{\top}$.
Figure~\ref{fig:spikecalibration} compares $F_{20}$ with the same
original-data statistic calibrated using population degrees of freedom
and with the exact Gaussian reference at that population critical value.
The empirical size of $F_{20}$ decreases from 7.97\% to 6.67\% along
the path. The Gaussian-reference rejection probabilities are about
5.9--6.0\%, while the limiting value in \eqref{eq:spikedsize} is 6.14\%.
The comparison separates error in estimating the calibration parameters,
finite-sample effects in the original-data statistic, and the persistent
matching error. The finite-sample curves remain distinct at the largest
setting considered.

\FloatBarrier
\section{Concluding remarks}\label{sec:discussion}
In this paper, an \(F\)-type calibration has been studied for testing
equality of high-dimensional covariance matrices with unknown means.
The statistic is based on the complete Li--Chen discrepancy. An exact
mean-correction expansion and independent-block trace estimators justify
its Gaussian reference under the stated conditions. The reference allows
both normal and nonnormal limits, and yields local-power and consistency
results for the feasible test.

The simulations show that averaging the trace estimates improves
stability. They also show that the ordinary \(F\) approximation can
remain liberal under strong correlation, in agreement with the spiked
example. The \(F\) and chi-square calibrations have similar empirical
sizes and powers. After adjustment for null rejection rates, neither has
a uniform power advantage over the specified three-cumulant comparator.
In the diffuse Gaussian examples, the local-power approximation is
accurate at the largest setting considered. The NR3 diagnostics identify
a finite-sample bias in cumulant estimation, distinct from error in the
Gaussian-reference approximation.

Further work could replace moment matching by direct numerical evaluation
of the weighted chi-square reference \citep{Imhof1961,Davies1980}.
This would require suitable estimates of the reference spectrum.
Multiplier-bootstrap comparisons and a real-data application would also
help assess finite-sample performance. On the theoretical side,
alternatives with a nonnegligible linear term require a separate power
approximation.

\section*{Computational reproducibility}
The accompanying code provides the test implementations, simulation
settings, random seeds, and numerical summaries. Instructions for
reproducing the tables are given in \texttt{README.txt}.

\section*{Acknowledgements}
The author used OpenAI's ChatGPT for language polishing and assistance
with implementing and checking supplementary simulation code, checking supplementary
mathematical calculations, and checking proofs. The author reviewed and corrected all AI-assisted output and takes full responsibility for the content of this paper.

{\small
\setlength{\bibsep}{2pt}
\bibliographystyle{plainnat}
\bibliography{references}
}

\clearpage
\appendix
\numberwithin{equation}{section}
\numberwithin{theorem}{section}
\numberwithin{proposition}{section}
\numberwithin{lemma}{section}
\numberwithin{corollary}{section}
\numberwithin{assumption}{section}
\numberwithin{remark}{section}
\numberwithin{table}{section}
\section{Independent-feature reference theory}\label{app:features}
The results below apply to independent feature observations. For the oracle
features, take \(\bm{Z}_{ij}=\bm{Z}_{ij}^o\), \(m_i=n_i\) and
\(\bm{\Omega}_i=\bm{\Omega}_i^o\). The proofs are given in
Appendices~\ref{app:identities}--\ref{app:momentconditions}.

\subsection{Independent features and the normal reference}\label{sec:reference}
\subsubsection{Test statistic}\label{sec:independent}
Consider independent feature samples \(\bm{Z}_{ij}\in\R^d\),
\(j=1,\ldots,m_i\), with mean \(\bm{\theta}_i\) and covariance \(\bm{\Omega}_i\).
Besides the oracle features in \eqref{eq:oraclefeatures}, this formulation
includes features obtained from disjoint pair differences:
\begin{equation}\label{eq:paired}
 \bm{A}_{ij}=\frac{\bm{X}_{i,2j}-\bm{X}_{i,2j-1}}{\sqrt2},\qquad
 \bm{Z}_{ij}=\svec(\bm{A}_{ij}\bm{A}_{ij}^\top),\quad
 m_i=\lfloor n_i/2\rfloor.
\end{equation}
These features are independent with mean \(\bm{\theta}_i\) and covariance
\begin{equation}\label{eq:paircov}
 \bm{\Omega}_i=\frac12(\bm{\Omega}_i^{\,o}+\bm{\Gamma}_i).
\end{equation}
The corresponding reference parameters are computed from \(m_i\) and
\(\bm{\Omega}_i\). For the full-sample test, these are replaced by
\(n_i\) and \(\bm{\Omega}_i^{\,o}\).

Let
\begin{align}
 \bm{\bar Z}_i&=\frac1{m_i}\sum_j \bm{Z}_{ij},&
 \bm{\widehat\Omega}_i&=\frac1{m_i-1}\sum_j
        (\bm{Z}_{ij}-\bm{\bar Z}_i)(\bm{Z}_{ij}-\bm{\bar Z}_i)^\top,\notag\\
 \bm{V}&=\frac{\bm{\Omega}_1}{m_1}+\frac{\bm{\Omega}_2}{m_2},&
 b&=\tr(\bm{V}),&
 B&=\frac{\tr(\bm{\widehat\Omega}_1)}{m_1}
      +\frac{\tr(\bm{\widehat\Omega}_2)}{m_2}.\label{eq:VB}
\end{align}
Define
\begin{equation}\label{eq:F}
 Q=\norm{\bm{\bar Z}_1-\bm{\bar Z}_2}^2,\qquad
 U=Q-B,\qquad F=\frac QB=1+\frac UB.
\end{equation}
The centered statistic \(U\) is the two-sample inner-product U-statistic
of \citet{ChenQin2010}, applied here to independent features.
Under \(H_0\), \(\E Q=\E B=b\).

\subsubsection{Moment conditions}
Write \(\bm{u}_{ij}=\bm{Z}_{ij}-\bm{\theta}_i\), \(t_i=\tr(\bm{\Omega}_i)\), \(s_i=\tr(\bm{\Omega}_i^2)\), and
\begin{align}
 J&=\tr(\bm{V}^2)
   =\frac{s_1}{m_1^2}
    +\frac{2\tr(\bm{\Omega}_1\bm{\Omega}_2)}{m_1m_2}
    +\frac{s_2}{m_2^2},\label{eq:J}\\
 R&=\sum_{i=1}^2\frac{s_i}{m_i^2(m_i-1)},&
 a^2&=2(J+R),\label{eq:a}\\
 \kappa_i&=\E\norm{\bm{u}_{i1}}^4-t_i^2-2s_i,&
 C_\kappa&=\sum_{i=1}^2\frac{\kappa_i}{m_i^3}.\label{eq:kappa}
\end{align}
Dependence on sample sizes and dimension is suppressed in the notation.

\begin{assumption}\label{ass:sampling}
The feature samples are independent and identically distributed within groups, with \(b>0\), \(a>0\), \(m_{\min}=\min(m_1,m_2)\to\infty\), and \(m_1/(m_1+m_2)\) bounded away from zero and one.
\end{assumption}

\begin{assumption}\label{ass:fourth}
There is a constant \(K<\infty\), independent of the dimension and sample sizes, such that
\[
 \E(\bm{v}^\top \bm{u}_{i1})^4\le K(\bm{v}^\top\bm{\Omega}_i \bm{v})^2
 \quad\text{for all }\bm{v}\in\R^d,\quad i=1,2.
\]
\end{assumption}
Assumption~\ref{ass:fourth} also gives
\(\E\norm{\bm{M} \bm{u}_{i1}}^4\le K\{\tr(\bm{M}\bm{\Omega}_i\bm{M}^\top)\}^2\) for every
deterministic matrix \(\bm{M}\), by expansion in an orthonormal basis and
Cauchy--Schwarz. For quadratic features, it controls fourth moments of
centered quadratic forms and thus requires uniform control of the relevant
original-data eighth moments.

\begin{remark}[Covariance normal-reference conditions]
\label{rem:frameworkconditions}
For the oracle features, \(\bm{v}^\top \bm{u}=\bm{\epsilon}^\top \bm{A}\bm{\epsilon}-\tr(\bm{A}\bm{\Sigma})\)
for some symmetric \(\bm{A}\). Assumption~\ref{ass:fourth} is therefore the
bounded directional kurtosis condition used in covariance normal-reference
analysis \citep{WangZhuZhang2024,WangZhuZhang2026}; the matrix formulation
above is equivalent by taking \(\bm{M}=\bm{v}^\top\). Isometric \(\svec\) removes
redundant coordinates while preserving inner products and nonzero covariance
eigenvalues. Section~\ref{sec:verifiedmodels} verifies this condition for
independent-coordinate factor models. The additional trace condition
\eqref{eq:blockdominance} is used for unknown-mean transfer and full-sample
calibration. No dimension-to-sample-size ratio is imposed, while group sample
sizes increase as specified in Assumption~\ref{ass:sampling}.
\end{remark}

\begin{proposition}[Cancellation in the centered statistic]\label{prop:moments}
Under \(H_0\), provided the relevant fourth feature moments are finite,
\begin{align}
 \Var(Q)&=2J+C_\kappa,&
 \Var(B)&=2R+C_\kappa,&
 \Cov(Q,B)&=C_\kappa,\label{eq:jointmom}\\
 \Var(U)&=a^2=2(J+R),&
 \Cov(U,B)&=-2R.\label{eq:Umom}
\end{align}
Under \(H_0\), \(U=U_0\), where the statistic formed from the centered feature errors is
\begin{equation}\label{eq:Udeg}
 U_0=
 \sum_{i=1}^2\frac1{m_i(m_i-1)}
     \sum_{j\ne k}\bm{u}_{ij}^\top \bm{u}_{ik}
 -\frac2{m_1m_2}\sum_{j,k}\bm{u}_{1j}^\top \bm{u}_{2k}.
\end{equation}
\end{proposition}
The fourth-cumulant terms cancel in \(U=Q-B\), whose variance therefore
depends only on the feature covariance matrices.

\begin{remark}
For a scalar \(X\sim N(0,\sigma^2)\), \(Z=X^2\) has variance \(2\sigma^4\) and fourth centered moment \(60\sigma^8\). Hence \(\kappa=48\sigma^8\), even though the original observation is Gaussian. If \(\widehat\Omega\) is the sample variance of \(m\) such features,
\[
 \Var(\widehat\Omega)
 =\frac{56\sigma^8}{m}
  +\frac{8\sigma^8}{m(m-1)}.
\]
In contrast, the Gaussian-feature reference variance is
\(8\sigma^8/(m-1)\). Thus normality of the original observations does not
make the reference distribution exact.
\end{remark}

\subsubsection{The Gaussian reference}
Construct independent samples
\[
 \bm{Z}_{ij}^*\sim N_d(\bm{\theta}_i,\bm{\Omega}_i)
\]
and calculate \(Q^*,B^*,U^*,F^*\) using the same definitions. A possibly singular Gaussian distribution is allowed.

\begin{proposition}[Independent-mixture representation]\label{prop:mixtures}
Under \(H_0\), let \(\lambda_r\) be the eigenvalues of \(\bm{V}\), and \(\lambda_{ir}\) those of \(\bm{\Omega}_i\), including zeros. Then
\begin{align}
 Q^*&\ \overset d=\ \sum_{r=1}^d\lambda_r \xi_r,
       &\xi_r&\stackrel{\mathrm{iid}}{\sim}\chi_1^2,\label{eq:Qstar}\\
 B^*&\ \overset d=\ \sum_{i=1}^2\sum_{r=1}^d
       \frac{\lambda_{ir}}{m_i(m_i-1)}\eta_{ir},
       &\eta_{ir}&\sim\chi_{m_i-1}^2.\label{eq:Bstar}
\end{align}
All variables on the right-hand sides may be taken mutually independent. In particular, \(Q^*\) and \(B^*\) are independent, with means \(b\) and variances \(2J\) and \(2R\), respectively.
\end{proposition}
Moment matching of the two independent chi-square mixtures gives the
\(F\) approximation in Appendix~\ref{sec:calibration}.

\subsubsection{Gaussian approximation}

For probability laws on \(\R\), let \(d_{\rm K}\) denote Kolmogorov distance.
Write \(G_n(x)=\Pp(U_0^*/a\le x)\). The two properties needed below are
\begin{align}
 d_{\rm K}\{\mathcal L(U_0/a),\mathcal L(U_0^*/a)\}&\longrightarrow0,
 \label{eq:Gcondition}\\
 \lim_{h\downarrow0}\limsup_n\sup_x
 \{G_n(x+h)-G_n(x-h)\}&=0.\label{eq:AC}
\end{align}
Here \(U_0^*\) is formed from \eqref{eq:Udeg} using independent
\(N_d(\bm{0},\bm{\Omega}_i)\) vectors. Theorem~\ref{thm:feature-universality}
establishes the Gaussian comparison, and
Lemma~\ref{lem:gaussian-anticoncentration} gives the required
anti-concentration bound. The proof uses Gaussian replacement and the
small-ball inequality of \citet{CarberyWright2001}. Related comparison
results were obtained by \citet{WangXu2022} and applied to covariance
features by \citet{WangZhuZhang2024}.

\begin{theorem}[Gaussian universality of independent features]
\label{thm:feature-universality}
Suppose the two feature samples are independent, are identically
distributed within groups, satisfy Assumption~\ref{ass:fourth}, and have
\(a^2=2(J+R)>0\). For \(m_1,m_2\ge2\), there is a constant
\(C_K<\infty\), depending only on \(K\), such that, for every
\(f\in C_b^3(\R)\),
\begin{equation}\label{eq:universality-smooth}
 \left|\E f(U_0/a)-\E f(U_0^*/a)\right|
 \le C_K m_{\min}^{-1/2}\norm{f'''}_\infty.
\end{equation}
Moreover,
\begin{equation}\label{eq:universality-kolm}
 d_{\rm K}\{\mathcal L(U_0/a),\mathcal L(U_0^*/a)\}
 \le C_K m_{\min}^{-1/14}.
\end{equation}
The constants are independent of the feature dimension and covariance
spectra.  No restriction on the ratio \(m_1/m_2\) is needed.
\end{theorem}

\begin{lemma}[Uniform anti-concentration of the Gaussian companion]
\label{lem:gaussian-anticoncentration}
There is an absolute constant \(C<\infty\) such that
\begin{equation}\label{eq:gaussian-anticoncentration}
 \sup_{x\in\R}\Pp\{|U_0^*/a-x|\le h\}\le C h^{1/2},
 \qquad h>0.
\end{equation}
This holds for all \(m_1,m_2\ge2\), all feature dimensions, and all
positive-semidefinite \(\bm{\Omega}_1,\bm{\Omega}_2\) for which \(a>0\).
\end{lemma}

Together, these results imply \eqref{eq:Gcondition} and \eqref{eq:AC}
under Assumptions~\ref{ass:sampling}--\ref{ass:fourth}. The rate in
\eqref{eq:universality-kolm} is sufficient for the calibration results.

\begin{remark}[Gaussian original observations]
\label{rem:gaussian-feature-fourth}
For \(\bm{G}\sim N_p(\bm{0},\bm{\Sigma})\), the centered feature
\(\bm{u}=\svec(\bm{G}\bm{G}^\top-\bm{\Sigma})\) satisfies Assumption~\ref{ass:fourth}
with \(K=15\), uniformly in \(p\) and \(\bm{\Sigma}\).
Indeed, for each feature direction \(\bm{v}\), there is a symmetric matrix
\(\bm{A}\) such that \(\bm{v}^\top \bm{u}=\bm{G}^\top \bm{A} \bm{G}-\tr(\bm{A}\bm{\Sigma})\).
If \(\lambda_r\) are the eigenvalues of
\(\bm{\Sigma}^{1/2}\bm{A}\bm{\Sigma}^{1/2}\), this variable has the distribution of
\(\sum_r\lambda_r(g_r^2-1)\), with independent standard normal
\(g_r\).  Its variance is \(2\sum_r\lambda_r^2\), and its fourth
moment is
\[
 12\left(\sum_r\lambda_r^2\right)^2+48\sum_r\lambda_r^4
 \le 60\left(\sum_r\lambda_r^2\right)^2.
\]
The same moment bound holds for features formed from Gaussian disjoint
pair differences.
\end{remark}

\subsubsection{Approximation of the ratio statistic}
\begin{lemma}[Relative consistency of the denominator]\label{lem:B}
Under Assumptions~\ref{ass:sampling}--\ref{ass:fourth},
\[
 B/b\longrightarrow_{\rm p}1,\qquad B^*/b\longrightarrow_{\rm p}1,
 \qquad R/J\le(m_{\min}-1)^{-1}.
\]
In particular, \(a^2=2J\{1+o(1)\}\).
\end{lemma}

\begin{theorem}[Ratio transfer]\label{thm:transfer}
Under \(H_0\) and Assumptions~\ref{ass:sampling}--\ref{ass:fourth},
\begin{equation}\label{eq:ratioK}
 d_{\rm K}\left\{\mathcal L\left(\frac b a(F-1)\right),
                  \mathcal L\left(\frac b a(F^*-1)\right)\right\}
 \longrightarrow0.
\end{equation}
Any deterministic convention on \(B\le0\) leaves this conclusion unchanged.
\end{theorem}
The transfer follows from
\[
 \frac b a(F-1)=\frac{U_0/a}{B/b},
\]
together with denominator consistency and Gaussian anti-concentration.

\subsection{\texorpdfstring{\(F\)}{F} calibration}\label{sec:calibration}
\subsubsection{Moment matching}
Matching the first two moments in \eqref{eq:Qstar}--\eqref{eq:Bstar} gives
\begin{equation}\label{eq:df}
 Q^*\approx\frac b{d_1}\chi_{d_1}^2,\qquad
 B^*\approx\frac b{d_2}\chi_{d_2}^2,\qquad
 d_1=\frac{b^2}{J},\quad d_2=\frac{b^2}{R}.
\end{equation}
This Welch--Satterthwaite approximation
\citep{Satterthwaite1946,Welch1947} yields an \(F_{d_1,d_2}\) ratio.
Its distributional error is
\begin{equation}\label{eq:WSerror}
 \varepsilon_{\rm WS}
 =d_{\rm K}\{\mathcal L(F^*),F_{d_1,d_2}\}.
\end{equation}

\begin{proposition}[Separation of the matching errors]\label{prop:WS}
Let
\[
 \varepsilon_Q=d_{\rm K}\{\mathcal L(Q^*/b),\mathcal L(\chi_{d_1}^2/d_1)\},
 \quad
 \varepsilon_B=d_{\rm K}\{\mathcal L(B^*/b),\mathcal L(\chi_{d_2}^2/d_2)\}.
\]
Then \(\varepsilon_{\rm WS}\le\varepsilon_Q+\varepsilon_B\). Moreover,
\begin{equation}\label{eq:dfbound}
 d_1\ge1,\qquad d_2\ge(m_{\min}-1)d_1.
\end{equation}
If
\begin{equation}\label{eq:rho3}
 \rho_3=\frac{\tr(\bm{V}^3)}{\{\tr(\bm{V}^2)\}^{3/2}}\longrightarrow0,
\end{equation}
then \(\varepsilon_Q\to0\). As \(m_{\min}\to\infty\), \(\varepsilon_B\to0\). Consequently, \eqref{eq:rho3} is a sufficient spectral condition for \(\varepsilon_{\rm WS}\to0\).
\end{proposition}

Condition \eqref{eq:rho3} describes a diffuse spectrum. The numerator
matching is also exact when all nonzero eigenvalues of \(\bm{V}\) are equal.
For a general spiked spectrum, \(\varepsilon_Q\) may persist. A simple
reference example has two nonzero eigenvalues proportional to \(2\) and
\(1\), giving
\[
 Q^*/b\ \overset d=\ (2\chi_1^2+\widetilde\chi_1^2)/3,
 \qquad d_1=9/5.
\]
The third cumulants are \(8/3\) and \(200/81\), respectively, so the
reference law and its moment match differ. Section~\ref{sec:originalspectralcheck}
examines such approximation error for covariance-feature models.
Since \(d_2\to\infty\), the denominator correction vanishes to first
order; see Proposition~\ref{prop:Fchi}.

\subsubsection{Estimating the feature traces}\label{sec:traces}
For independent features, distribution-free trace estimators follow from
pairing observations. Let \(L_i=\lfloor m_i/2\rfloor\), and set
\begin{equation}\label{eq:featurecontrasts}
 \bm{T}_{ir}=(\bm{Z}_{i,2r}-\bm{Z}_{i,2r-1})/\sqrt2,\qquad r=1,\ldots,L_i.
\end{equation}
These contrasts have mean zero and covariance \(\bm{\Omega}_i\). The statistic
\(F\) continues to use all \(m_i\) observations.
Define
\begin{align}
 \widehat s_i&=\frac1{L_i(L_i-1)}
      \sum_{r\ne s}(\bm{T}_{ir}^\top \bm{T}_{is})^2,\label{eq:shat}\\
 \widehat c_{12}&=\frac1{L_1L_2}
      \sum_{r=1}^{L_1}\sum_{s=1}^{L_2}(\bm{T}_{1r}^\top \bm{T}_{2s})^2,\label{eq:chat}\\
 \widehat J&=\frac{\widehat s_1}{m_1^2}
       +\frac{2\widehat c_{12}}{m_1m_2}
       +\frac{\widehat s_2}{m_2^2},&
 \widehat R&=\sum_{i=1}^2
       \frac{\widehat s_i}{m_i^2(m_i-1)}.\label{eq:JRhat}
\end{align}
Each summand uses distinct contrasts, and the resulting estimators are
nonnegative. The estimated degrees of freedom are
\begin{equation}\label{eq:dfhat}
 \widehat d_1=B^2/\widehat J,\qquad
 \widehat d_2=B^2/\widehat R.
\end{equation}
At least \(m_i\ge4\) is needed for the within-group estimators.

\begin{proposition}[Trace estimation]\label{prop:traces}
The estimators \(\widehat s_i\) and \(\widehat c_{12}\) are unbiased for \(s_i\) and \(\tr(\bm{\Omega}_1\bm{\Omega}_2)\), respectively. Under Assumptions~\ref{ass:sampling}--\ref{ass:fourth}, their nonzero targets are estimated with relative mean-square error \(O(L_{\min}^{-1})\). Zero targets yield identically zero estimators. Consequently,
\[
 \widehat J/J\to_{\rm p}1,\qquad
 \widehat R/R\to_{\rm p}1,\qquad
 \widehat d_j/d_j\to_{\rm p}1,\quad j=1,2.
\]
\end{proposition}
These estimators give a feasible calibration for the paired features.
Section~\ref{sec:fullblocktraces} develops the corresponding estimators for
the full-sample statistic.

\subsubsection{Feasible critical values and size}
Let \(q_{1-\alpha}(d_1,d_2)\) denote the lower \(1-\alpha\) quantile of \(F_{d_1,d_2}\). The proposed independent-feature procedure rejects if
\begin{equation}\label{eq:test}
 F>q_{1-\alpha}(\widehat d_1,\widehat d_2).
\end{equation}
The rule assigns nonrejection and records a calibration failure when
\(B,\widehat J,\widehat R\) do not yield positive, finite degrees of
freedom. The probability of this event tends to zero under the stated
consistency conditions.

\begin{lemma}[Random critical values]\label{lem:quantile}
For fixed \(0<\alpha<1\), under \eqref{eq:dfbound} with \(m_{\min}\to\infty\), relative consistency of both degrees of freedom implies
\[
 \frac b a\{q_{1-\alpha}(\widehat d_1,\widehat d_2)
       -q_{1-\alpha}(d_1,d_2)\}\to_{\rm p}0.
\]
\end{lemma}

\begin{proposition}[First-order equivalence of reference critical values]
\label{prop:Fchi}
Set \(q_{1-\alpha}(d,\infty)=\chi^2_{d,1-\alpha}/d\), where
\(\chi^2_{d,1-\alpha}\) is the lower \(1-\alpha\) chi-square quantile.
For fixed \(0<\alpha<1\), under \eqref{eq:dfbound} and
\(m_{\min}\to\infty\),
\[
 \frac b a\{q_{1-\alpha}(d_1,d_2)-q_{1-\alpha}(d_1,\infty)\}
 \longrightarrow0.
\]
The same conclusion holds in probability with both degrees of freedom
replaced by relatively consistent estimates, retaining the deterministic
normalization \(b/a\).
\end{proposition}
The \(F\) and chi-square critical values therefore agree to first order.
Their finite-sample difference is examined in Section~\ref{sec:simulation}.

\begin{theorem}[Size approximation]\label{thm:size}
Under the assumptions of Theorem~\ref{thm:transfer}, Proposition~\ref{prop:traces}, and Lemma~\ref{lem:quantile}, the rejection probability in \eqref{eq:test} satisfies
\begin{equation}\label{eq:size}
 \left|\Pp_0\{F>q_{1-\alpha}(\widehat d_1,\widehat d_2)\}-\alpha\right|
 \le \varepsilon_{\rm WS}+o(1).
\end{equation}
In particular, the test has asymptotic size \(\alpha\) when the additional conditions in Proposition~\ref{prop:WS} make \(\varepsilon_{\rm WS}\to0\).
\end{theorem}
The proof uses deterministic brackets to handle the estimated critical value.

\subsubsection{Feature-Gram calculations}
For features generated from vectors \(\bm{a},\bm{b}\in\R^p\),
\[
 \svec(\bm{a}\bm{a}^\top)^\top\svec(\bm{b}\bm{b}^\top)=(\bm{a}^\top \bm{b})^2.
\]
Thus all terms involving feature means, centered traces, and contrasts in
\eqref{eq:VB}--\eqref{eq:JRhat} can be evaluated from a feature Gram matrix.
For \(M=m_1+m_2\), forming the original-vector Gram matrix requires
\(O(M^2p)\) operations and \(O(M^2)\) memory. The subsequent trace sums
are computed from this Gram matrix without constructing a
\(d\times d\) feature covariance matrix. Blocking reduces the memory
requirement.

\subsection{Independent-feature alternatives}\label{sec:power}
Let \(\bm{\Delta}=\bm{\theta}_1-\bm{\theta}_2\) and \(\bm{w}=\bm{\bar u}_1-\bm{\bar u}_2\). Then
\begin{equation}\label{eq:alternative}
 U=U_0+2\bm{\Delta}^\top \bm{w}+\norm{\bm{\Delta}}^2,\qquad
 \Var(U)=a^2+4\bm{\Delta}^\top \bm{V}\bm{\Delta}.
\end{equation}
Degeneracy gives \(\Cov(U_0,\bm{w})=\bm{0}\), so the variance identity holds without
a symmetry assumption.

\begin{theorem}[Quadratic-dominated local alternatives]\label{thm:power}
Suppose Assumptions~\ref{ass:sampling}--\ref{ass:fourth} hold along a sequence of alternatives for the centered feature distributions, and use the trace estimators in Appendix~\ref{sec:traces}. If
\[
 h_n=\norm{\bm{\Delta}}^2/a=O(1),\qquad
 \bm{\Delta}^\top \bm{V}\bm{\Delta}=o(a^2),
\]
then, with \(G_n\) as in \eqref{eq:AC} and
\[
 t_n=\frac b a\{q_{1-\alpha}(d_1,d_2)-1\},
\]
the rejection probability is
\begin{equation}\label{eq:power}
 \Pp\{F>q_{1-\alpha}(\widehat d_1,\widehat d_2)\}
       =1-G_n(t_n-h_n)+o(1).
\end{equation}
If additionally \eqref{eq:rho3} holds and \(h_n\to h\), this converges to
\(1-\Phi(z_{1-\alpha}-h)\).
\end{theorem}
The condition \(\bm{\Delta}^\top \bm{V}\bm{\Delta}=o(a^2)\) makes the linear term
negligible. For signals aligned with a dominant eigenvector, this term can
remain of leading order, giving a weighted noncentral chi-square reference.

\begin{proposition}[A sufficient consistency condition]\label{prop:consistency}
Suppose \(B/b\to_{\rm p}1\), the plug-in critical value satisfies
\(\widehat q-1=\Op(a/b)\), and
\[
 \delta_n=\norm{\bm{\Delta}}^2,\qquad
 \frac{a^2+4\bm{\Delta}^\top \bm{V}\bm{\Delta}}{\delta_n^2}\longrightarrow0.
\]
Then the rejection probability tends to one.
\end{proposition}

\section{Finite-sample identities}\label{app:identities}
Throughout the proofs, constants denoted by \(C\) may change from line to
line. Unless stated otherwise, they are independent of the dimension and
sample sizes.

\subsection{Proof of Proposition~\ref{prop:full}}
Consider one group and write \(\bm{\epsilon}_j=\bm{X}_j-\bm{\mu}\),
\(\bm{Y}_j=\svec(\bm{\epsilon}_j\bm{\epsilon}_j^\top)\), and
\(\bm{\theta}=\svec(\bm{\Sigma})\). Applying \(\svec\) to the sample covariance
decomposition gives
\[
 \svec(\bm{S})-\bm{\theta}
 =\frac1n\sum_j(\bm{Y}_j-\bm{\theta})
 -\frac1{n(n-1)}\sum_{j\ne k}
       \svec\{(\bm{\epsilon}_j\bm{\epsilon}_k^\top+
                    \bm{\epsilon}_k\bm{\epsilon}_j^\top)/2\}.
\]
The two terms are uncorrelated. Indeed, every mixed product contains
an independent centered observation that occurs only once. The covariance
of the first term is \(\bm{\Omega}^o/n\). For symmetric matrices \(\bm{A},\bm{B}\),
independence also gives
\[
 \E\{(\bm{\epsilon}_1^\top \bm{A}\bm{\epsilon}_2)
       (\bm{\epsilon}_1^\top \bm{B}\bm{\epsilon}_2)\}
       =\tr(\bm{A}\bm{\Sigma} \bm{B}\bm{\Sigma}).
\]
Only identical or reversed ordered pairs contribute to the covariance
of the second term. Their sum is \(\bm{\Gamma}/\{n(n-1)\}\), which proves
\eqref{eq:covS}.

For four distinct indices, \(\bm{h}_{ab}\) and \(\bm{h}_{ce}\) are independent
with common expectation \(\bm{\theta}\). Thus \(\E D_n=\norm{\bm{\theta}}^2\).
Unbiasedness of \(\bm{S}\) further yields
\[
 \E\norm{\bm{S}}_F^2-\norm{\bm{\Sigma}}_F^2
       =\tr\Cov\{\svec(\bm{S})\}.
\]
Combining these identities across the two independent groups proves
the expectation formulas for \(B_{\rm full}\) and \(U_{\rm full}\).
Both the pair differences and the sample covariance matrices are
translation invariant.

For \eqref{eq:paircov}, let \(\bm{\widetilde\epsilon}\) be an independent
copy of \(\bm{\epsilon}\). The contrast
\((\bm{\epsilon}-\bm{\widetilde\epsilon})/\sqrt2\) has covariance \(\bm{\Sigma}\)
and fourth cumulant tensor equal to half that of \(\bm{\epsilon}\).
Its quadratic-feature covariance is therefore
\(\bm{\Gamma}+(\bm{\Omega}^o-\bm{\Gamma})/2\).

\subsection{Proof of Proposition~\ref{prop:LiChen}}
Let \(\sum'\) denote summation over ordered, mutually distinct
indices within a group. Expansion of \eqref{eq:D}, followed by
relabelling the indices, gives
\begin{align*}
 D_n
 &=\frac1{4(n)_4}\sum_{a,b,c,e}'
       \{(\bm{X}_a-\bm{X}_b)^\top(\bm{X}_c-\bm{X}_e)\}^2\\
 &=\frac1{(n)_2}\sum_{a,b}'(\bm{X}_a^\top \bm{X}_b)^2
   -\frac2{(n)_3}\sum_{a,b,c}'
          (\bm{X}_a^\top \bm{X}_b)(\bm{X}_b^\top \bm{X}_c)\\
 &\quad+\frac1{(n)_4}\sum_{a,b,c,e}'
          (\bm{X}_a^\top \bm{X}_b)(\bm{X}_c^\top \bm{X}_e)
   =A_n^{\rm LC}.
\end{align*}
For the cross-group estimator, define
\[
 \bm{H}_i=\frac1{n_i}\sum_a \bm{X}_{ia}\bm{X}_{ia}^\top,
 \qquad
 \bm{L}_i=\frac1{(n_i)_2}\sum_{a\ne b}\bm{X}_{ia}\bm{X}_{ib}^\top.
\]
The exact identity \(\bm{S}_i=\bm{H}_i-\bm{L}_i\) and the four terms of the
Li--Chen cross-group estimator yield
\[
 C_{n_1,n_2}^{\rm LC}
 =\tr(\bm{H}_1\bm{H}_2)-\tr(\bm{H}_1\bm{L}_2)-\tr(\bm{L}_1\bm{H}_2)+\tr(\bm{L}_1\bm{L}_2)
 =\tr(\bm{S}_1\bm{S}_2).
\]
Combining these identities with \eqref{eq:Ffull} proves the assertion.

\subsection{Proof of Proposition~\ref{prop:moments}}
For centered independent vectors \(\bm{u}_1,\ldots,\bm{u}_m\) with common
covariance \(\bm{\Omega}\), set \(t=\tr(\bm{\Omega})\),
\(s=\tr(\bm{\Omega}^2)\), and
\(\kappa=\E\norm{\bm{u}_1}^4-t^2-2s\). Expanding
\(\norm{\bm{\bar u}}^2=m^{-2}\sum_{j,k}\bm{u}_j^\top \bm{u}_k\) and retaining
only terms with matching indices gives
\[
 \Var(\norm{\bm{\bar u}}^2)=\frac{2s}{m^2}+\frac{\kappa}{m^3}.
\]
The same calculation for the unbiased sample covariance
\(\bm{\widehat\Omega}\) gives
\begin{align*}
 \Var\{\tr(\bm{\widehat\Omega})\}
 &=\frac{\E\norm{\bm{u}_1}^4-t^2}{m}
       +\frac{2s}{m(m-1)}
  =\frac{2s}{m-1}+\frac{\kappa}{m},\\
 \Cov\{\norm{\bm{\bar u}}^2,\tr(\bm{\widehat\Omega})\}
 &=\frac{\kappa}{m^2}.
\end{align*}
Independence of the two groups implies that the cross-group term in
\(Q\) has variance \(4\tr(\bm{\Omega}_1\bm{\Omega}_2)/(m_1m_2)\) and is
uncorrelated with the within-group squared means and with \(B\).
The joint moment identities \eqref{eq:jointmom}--\eqref{eq:Umom}
follow by addition. Subtracting the sample-trace terms from the squared
means gives the degenerate representation \eqref{eq:Udeg}.

\subsection{Proof of Proposition~\ref{prop:mixtures}}
Under the reference model,
\(\bm{\bar Z}_1^*-\bm{\bar Z}_2^*\sim N_d(\bm{0},\bm{V})\), so spectral decomposition
of \(\bm{V}\) gives \eqref{eq:Qstar}. The Gaussian sample mean and sample
covariance are independent within each group. For a singular covariance,
the same assertion follows by working on its support. Diagonalizing
\(\bm{\Omega}_i\) gives \eqref{eq:Bstar}; independence across groups
then proves independence of the numerator and denominator.

\subsection{Computational identities}\label{app:blockcomputation}
The full-sample statistic also admits a Gram calculation. Decomposing
pairs according to their overlap reduces \eqref{eq:D} as follows. If \(P=\binom n2\), \(A_n=P^{-1}\sum_{a<b}\norm{\bm{h}_{ab}}^2\), and \(E_n\) is the average of \(\bm{h}_{ab}^\top \bm{h}_{ce}\) over ordered distinct unordered pairs sharing one index, then
\begin{equation}\label{eq:Doverlap}
 D_n=\frac{n(n-1)\norm{\bm{S}}_F^2-2A_n-4(n-2)E_n}
            {(n-2)(n-3)}.
\end{equation}

The estimators can be evaluated from Gram matrices without forming a
\(d\times d\) matrix. For blocks represented by \((\bm{a},\bm{b})\) and
\((\bm{c},\bm{d})\), the required inner products are
\begin{align}
 \bm{v}(\bm{a},\bm{b})^\top \bm{v}(\bm{c},\bm{d})
  &=(\bm{a}^\top \bm{c})^2-(\bm{a}^\top \bm{d})^2-(\bm{b}^\top \bm{c})^2+(\bm{b}^\top \bm{d})^2,\notag\\
 \bm{v}(\bm{a},\bm{b})^\top \bm{w}(\bm{c},\bm{d})
  &=\sqrt2\{(\bm{a}^\top \bm{c})(\bm{a}^\top \bm{d})-(\bm{b}^\top \bm{c})(\bm{b}^\top \bm{d})\},\notag\\
 \bm{w}(\bm{a},\bm{b})^\top \bm{w}(\bm{c},\bm{d})
  &=(\bm{a}^\top \bm{c})(\bm{b}^\top \bm{d})+(\bm{a}^\top \bm{d})(\bm{b}^\top \bm{c}).
 \label{eq:fourblockgram}
\end{align}
Thus
\[
 \tr(\bm{K}_r\bm{K}_s)
  =(\bm{v}_r^\top \bm{v}_s)^2-(\bm{v}_r^\top \bm{w}_s)^2
    -(\bm{w}_r^\top \bm{v}_s)^2+(\bm{w}_r^\top \bm{w}_s)^2.
\]

For \eqref{eq:Doverlap}, there are \(P=\binom n2\) unordered
observation pairs. In the squared norm of their average feature,
the ordered pairs of edges are identical, share one index, or are
disjoint. Their respective counts are \(P\), \(2P(n-2)\), and
\(P\binom{n-2}2\). Hence
\[
 \norm{\bm{S}}_F^2
 =\frac{2A_n+4(n-2)E_n+(n-2)(n-3)D_n}{n(n-1)}.
\]
Solving gives \eqref{eq:Doverlap}. To obtain \eqref{eq:Dgram},
use translation invariance to replace the observations in the
Li--Chen expansion by the residuals. Since $\bm{G}\mathbf1=\bm{0}$,
\begin{align*}
 \sum_{a\ne b}G_{ab}^2&=u_G-q_G,\\
 \sum_{a,b,c}^{\ne}G_{ab}G_{bc}&=2q_G-u_G,\\
 \sum_{a,b,c,d}^{\ne}G_{ab}G_{cd}&=t_G^2+2u_G-6q_G.
\end{align*}
For the second identity, fix $b$ and subtract the terms $a=c$ from
$(\sum_{a\ne b}G_{ab})^2=G_{bb}^2$. For the third, start from
$(\sum_{a\ne b}G_{ab})^2=t_G^2$ and remove the two identical-edge
arrangements and four shared-index arrangements. Substitution into
\[
 D_n=\frac{u_G-q_G}{(n)_2}
      -\frac{2(2q_G-u_G)}{(n)_3}
      +\frac{t_G^2+2u_G-6q_G}{(n)_4}
\]
and collection of coefficients proves \eqref{eq:Dgram}.

\subsection{Residual covariance estimation}
\paragraph{Comparison with residual covariance estimation.}
For comparison with the block estimator, consider the sample covariance
of the residual outer products:
\[
 \bm{R}_j=\svec\{(\bm{X}_j-\bm{\bar X})(\bm{X}_j-\bm{\bar X})^\top\},\qquad
 \bm{C}_r=\frac1{n-1}\sum_j(\bm{R}_j-\bm{\bar R})(\bm{R}_j-\bm{\bar R})^\top.
\]
For finite fourth original-data moments and \(n>2\),
\begin{equation}\label{eq:residualbias}
 \E \bm{C}_r=\frac{(n-2)^2}{n^2}\bm{\Omega}^o
          +\frac{2(n-2)}{n^2}\bm{\Gamma}.
\end{equation}
Thus \(\bm{C}_r\) is generally biased for \(\bm{\Omega}^o\). Estimating
\(\tr\{(\bm{\Omega}^o)^2\}\) additionally requires control of the dependence
between the matrix factors.

Write \(c_{ja}=\ind\{j=a\}-1/n\) and
\(\bm{r}_j=\sum_a c_{ja}\bm{\epsilon}_a\). The fourth-moment expansion of
\(\bm{R}_j=\svec(\bm{r}_j\bm{r}_j^\top)\) gives
\[
 \Cov(\bm{R}_j,\bm{R}_k)=
 \left(\sum_a c_{ja}^2c_{ka}^2\right)(\bm{\Omega}^o-\bm{\Gamma})
 +\left(\sum_a c_{ja}c_{ka}\right)^2\bm{\Gamma}.
\]
The residual features are exchangeable with a common mean, so
\(\E \bm{C}_r=\Cov(\bm{R}_1,\bm{R}_1)-\Cov(\bm{R}_1,\bm{R}_2)\). Substituting the coefficients yields
\[
 \E \bm{C}_r=\frac{(n-2)^2}{n^2}\bm{\Omega}^o
             +\frac{2(n-2)}{n^2}\bm{\Gamma},
\]
which proves \eqref{eq:residualbias}. For the jackknife calculation,
let \(\bm{s}_{(-j)}\) be the vectorized unbiased covariance matrix after
observation \(j\) is deleted, and set
\(\bm{\bar s}_{(-\cdot)}=n^{-1}\sum_j\bm{s}_{(-j)}\). The delete-one covariance
jackknife is
\[
 \bm{J}_{\rm jack}=\frac{n-1}{n}\sum_j
 (\bm{s}_{(-j)}-\bm{\bar s}_{(-\cdot)})(\bm{s}_{(-j)}-\bm{\bar s}_{(-\cdot)})^\top.
\]
It satisfies
\[
 \bm{J}_{\rm jack}=\frac{n}{(n-2)^2}C_r,\qquad
 \E \bm{J}_{\rm jack}=\frac{\bm{\Omega}^o}{n}
                      +\frac{2\bm{\Gamma}}{n(n-2)}.
\]
The jackknife therefore differs from \eqref{eq:covS} in its second-order term.

\section{Gaussian comparison and reference calibration}\label{app:transfer}

\subsection{Proof of Lemma~\ref{lem:gaussian-anticoncentration}}
Represent the Gaussian feature vectors as deterministic linear images
of a vector \(\bm{g}\) of independent standard normal variables. This
representation also covers singular covariances. Then
\(P(\bm{g})=U_0^*/a\) is a polynomial of degree at most two, with
\(\E P=0\) and \(\E P^2=1\).

Theorem~8 of \citet{CarberyWright2001}, applied to a polynomial \(q\)
of degree at most two under standard Gaussian measure, gives
\[
 \Pp\{|q(\bm{g})|\le h\}
 \le C h^{1/2}\{\E |q(\bm{g})|^2\}^{-1/4}.
\]
This is the specialization of their theorem to degree two and
integrability parameter four. Apply it to \(P(\bm{g})-x\). Since
\(\E\{P(\bm{g})-x\}^2=1+x^2\), we obtain
\[
 \Pp\{|P(\bm{g})-x|\le h\}
 \le C h^{1/2}(1+x^2)^{-1/4}\le C h^{1/2}.
\]
This proves the bound uniformly in the coefficients.

\subsection{Proof of Theorem~\ref{thm:feature-universality}}
Set \(K_0=\max(K,3)\). The directional fourth-moment bound with
constant \(K_0\) holds for both the original centered features and
the matching Gaussian vectors. It also holds for any sum
\(\bm{L}=\sum_r c_r\bm{Y}_r\) of independent vectors drawn from these two
classes, because
\[
 \E\left(\sum_r c_r \bm{v}^\top \bm{Y}_r\right)^4
 =\sum_r c_r^4\E(\bm{v}^\top \bm{Y}_r)^4
  +6\sum_{r<s}c_r^2c_s^2
       \E(\bm{v}^\top \bm{Y}_r)^2\E(\bm{v}^\top \bm{Y}_s)^2
 \le K_0\{\Var(\bm{v}^\top \bm{L})\}^2.
\]
A second consequence is that a centered vector \(\bm{Y}\), with covariance
\(\bm{\Psi}\) and the same directional bound, satisfies
\begin{equation}\label{eq:psd-fourth-bound}
 \E(\bm{Y}^\top \bm{M} \bm{Y})^2\le K_0\{\tr(\bm{M}\bm{\Psi})\}^2
\end{equation}
for every positive semidefinite matrix \(\bm{M}\). To obtain this bound,
diagonalize \(\bm{M}\), expand the square, and apply Cauchy--Schwarz to
each mixed fourth moment.

Enumerate the \(m_1+m_2\) centered observations by \(\bm{Y}_\ell\), and
let \(\bm{\Psi}_\ell\) be the covariance of \(\bm{Y}_\ell\).  Write
\[
 T(\bm{Y}):=U_0/a=\sum_{\ell<r}\gamma_{\ell r}\bm{Y}_\ell^\top \bm{Y}_r,
\]
where \(\gamma_{\ell r}=2/\{a m_i(m_i-1)\}\) if both observations
are in group \(i\), and
\(\gamma_{\ell r}=-2/(a m_1m_2)\) if they are in different groups.
For each observation define its influence by
\[
 \nu_\ell=\sum_{r\ne\ell}\gamma_{\ell r}^2
                      \tr(\bm{\Psi}_\ell\bm{\Psi}_r).
\]
Distinct unordered-pair summands in \(T\) are uncorrelated, because
an observation occurring only once has conditional mean zero.
Consequently,
\begin{equation}\label{eq:feature-influences}
 \sum_\ell\nu_\ell=2\Var(T)=2,
 \qquad \max_\ell\nu_\ell\le\frac2{m_{\min}}.
\end{equation}
For example, when \(\ell\) belongs to group one and
\(c_{12}=\tr(\bm{\Omega}_1\bm{\Omega}_2)\), its influence is
\[
 \nu_\ell=\frac4{a^2}\left\{
  \frac{s_1}{m_1^2(m_1-1)}+\frac{c_{12}}{m_1^2m_2}\right\}
 \le\frac2{m_1}.
\]
The inequality follows from
\[
 a^2=\frac{2s_1}{m_1(m_1-1)}
       +\frac{2s_2}{m_2(m_2-1)}
       +\frac{4c_{12}}{m_1m_2}.
\]
Interchanging the group indices gives the other case.

Replace the observations one at a time by independent matching
Gaussian vectors.  At the step replacing observation \(\ell\),
condition on every other vector in the current mixture of original
and Gaussian observations.  The statistic has the form
\[
 A_\ell+\bm{Y}_\ell^\top \bm{L}_\ell,
 \qquad \bm{L}_\ell=\sum_{r\ne\ell}\gamma_{\ell r}\bm{\widetilde Y}_r.
\]
Conditionally on the remaining observations, \(A_\ell\) and
\(\bm{L}_\ell\) are fixed, and the first two moments of
\(\bm{Y}_\ell^\top \bm{L}_\ell\) agree before and after replacement.
Let \(W_\ell=\bm{Y}_\ell^\top \bm{L}_\ell\).  Its second moment is
\(\nu_\ell\).  Conditioning first on \(\bm{L}_\ell\), and then using
\eqref{eq:psd-fourth-bound} for \(\bm{L}_\ell\), gives
\[
 \E W_\ell^4
 \le K_0\E(\bm{L}_\ell^\top\bm{\Psi}_\ell \bm{L}_\ell)^2
 \le K_0^2\nu_\ell^2.
\]
Hence
\(\E|W_\ell|^3\le
\{\E W_\ell^2\E W_\ell^4\}^{1/2}
\le K_0\nu_\ell^{3/2}\).
The same moment bounds hold for the Gaussian replacement.

Taylor expansion of \(f(A_\ell+W_\ell)\) to second order about
\(A_\ell\) cancels the conditional terms through order two.
The difference at this step is therefore at most
\(C_K\norm{f'''}_\infty\nu_\ell^{3/2}\).
Summing the replacement errors and using
\eqref{eq:feature-influences} proves
\[
 |\E f(U_0/a)-\E f(U_0^*/a)|
 \le C_K\norm{f'''}_\infty
       (\max_\ell\nu_\ell)^{1/2}\sum_\ell\nu_\ell
 \le C_Km_{\min}^{-1/2}\norm{f'''}_\infty.
\]

To pass to distribution functions, fix a smooth nonincreasing
function \(\psi\) that equals one on \(( -\infty,0]\), equals zero
on \([1,\infty)\), and takes values in \([0,1]\).  The functions
\(f_{x,h}(t)=\psi\{(t-x)/h\}\) satisfy
\[
 \ind\{t\le x\}\le f_{x,h}(t)\le\ind\{t\le x+h\},
 \qquad \norm{f_{x,h}'''}_\infty\le C h^{-3}.
\]
Apply the smooth comparison to \(f_{x,h}\) for an upper bound
and to \(f_{x-h,h}\) for a lower bound.  Lemma~
\ref{lem:gaussian-anticoncentration} then yields, uniformly in \(x\),
\[
 \left|\Pp(U_0/a\le x)-\Pp(U_0^*/a\le x)\right|
 \le C h^{1/2}+C_Km_{\min}^{-1/2}h^{-3}.
\]
Balancing the two terms by taking \(h=m_{\min}^{-1/7}\) proves
\eqref{eq:universality-kolm}.

\subsection{Proof of Lemma~\ref{lem:B}}
By the trace variance identity in Appendix~\ref{app:identities} and
\(\E\norm{\bm{u}_{i1}}^4\le Kt_i^2\), we have
\[
 \Var(B)\le C\sum_i\frac{t_i^2}{m_i^3}
            \le \frac{C}{m_{\min}}b^2.
\]
Chebyshev's inequality proves relative consistency; the Gaussian
reference satisfies the same bound. Since \(s_i\ge0\) and
\(\tr(\bm{\Omega}_1\bm{\Omega}_2)\ge0\),
\[
 R\le \frac1{m_{\min}-1}\sum_i\frac{s_i}{m_i^2}
       \le \frac{J}{m_{\min}-1}.
\]

\subsection{Proof of Theorem~\ref{thm:transfer}}
The centered statistics \(U_0/a\) and \(U_0^*/a\) have mean zero
and variance one, so both sequences are tight. Lemma~\ref{lem:B}
and the ratio identity imply
\[
 \frac b a(F-1)-\frac{U_0}{a}\to_{\rm p}0,\qquad
 \frac b a(F^*-1)-\frac{U_0^*}{a}\to_{\rm p}0.
\]
For fixed \(h>0\), a perturbation bounded by \(h\) changes a
distribution function by at most its concentration function over an
interval of radius \(h\). An arbitrary perturbation adds the
probability that its absolute value exceeds \(h\). Apply this bound
to both displayed approximations and then use \eqref{eq:Gcondition}.
Letting the sample sizes tend to infinity and then \(h\downarrow0\),
\eqref{eq:AC} yields the asserted Kolmogorov convergence. Moreover, relative consistency gives \(\Pp(B\le0)\to0\).

\subsection{Proof of Proposition~\ref{prop:WS}}
For independent positive variables \(A,D\), condition first on
\(D\) in \(\Pp(A/D\le x)\), \(x>0\). Replacing the distribution
of \(A\) changes this probability by at most its Kolmogorov distance.
Conditioning instead on \(A\) gives the corresponding bound for a
replacement of \(D\), since the conditional function is monotone
with range in \([0,1]\). Replacing the two components successively
proves \(\varepsilon_{\rm WS}\le\varepsilon_Q+\varepsilon_B\).

Positive semidefiniteness of \(\bm{V}\) gives \(J\le b^2\), and
Lemma~\ref{lem:B} supplies the bound for \(d_2\). A Berry--Esseen
bound for the weighted independent variables \(\chi_1^2-1\)
controls the standardized numerator's normal approximation error by
\(C\rho_3\). The standardized matching gamma variable tends to
normal when \(d_1\to\infty\). The trace inequality
\[
 \{\tr(\bm{V}^2)\}^2\le\tr(\bm{V})\tr(\bm{V}^3)
\]
implies \(1/\sqrt{d_1}\le\rho_3\). Thus \eqref{eq:rho3}
ensures both normal approximations and hence \(\varepsilon_Q\to0\).

For the denominator write each \(\chi_{m_i-1}^2\) as the sum of \(m_i-1\) independent squared standard normals. If \(w_{ir}=\lambda_{ir}/\{m_i(m_i-1)\}\), the normalized third absolute moment bound is proportional to
\[
 \frac{\sum_{i,r}(m_i-1)w_{ir}^3}
      {\{\sum_{i,r}(m_i-1)w_{ir}^2\}^{3/2}}
 \le \frac1{\sqrt{m_{\min}-1}}.
\]
Its standardized distribution therefore tends to normal. Since \(d_2\to\infty\), so does that of the matching gamma variable; they share the same mean and variance. This proves \(\varepsilon_B\to0\).

\subsection{Proof of Proposition~\ref{prop:traces}}
Independence and iterated expectation give
\[
 \E(\bm{T}_{ir}^\top \bm{T}_{is})^2=\tr(\bm{\Omega}_i^2),\qquad
 \E(\bm{T}_{1r}^\top \bm{T}_{2s})^2=\tr(\bm{\Omega}_1\bm{\Omega}_2).
\]
The contrasts obey a uniform directional fourth-moment bound with a constant depending only on \(K\). For independent contrasts \(\bm{T},\bm{T}'\) with covariances \(\bm{\Omega},\bm{\Omega}'\), conditioning twice and expanding a positive semidefinite quadratic form imply
\[
 \E(\bm{T}^\top \bm{T}')^4\le C\{\tr(\bm{\Omega}\bm{\Omega}')\}^2.
\]
Disjoint kernel summands are independent. The fraction of overlapping
terms is \(O(L_i^{-1})\) for a within-group order-two average and
\(O(L_1^{-1}+L_2^{-1})\) for a two-sample average. Bounding the
remaining covariances by Cauchy--Schwarz gives the stated relative
mean-square rates. If a target is zero, unbiasedness and nonnegativity
make the corresponding estimator zero almost surely. Positive weighted
sums therefore give relative consistency of \(\widehat J\) and
\(\widehat R\). Combining this conclusion with Lemma~\ref{lem:B}
proves relative consistency of \(\widehat d_1,\widehat d_2\).

\subsection{Proof of Lemma~\ref{lem:quantile}}
Consider an arbitrary subsequence. Passing to a further subsequence,
either \(d_1\to d\in[1,\infty)\) or \(d_1\to\infty\).
In the first case, \(d_2\to\infty\), \(\widehat d_1\to_{\rm p}d\), and
\(\widehat d_2\to_{\rm p}\infty\). The limiting distribution
\(\chi_d^2/d\) has a continuous, strictly increasing distribution
function on its support. Quantile continuity therefore gives convergence
of both critical values to \(\chi^2_{d,1-\alpha}/d\). Since
\(b/a=\sqrt{d_1/2}\{1+o(1)\}\), the assertion follows.

For the second case, first consider deterministic sequences
\(\nu\to\infty\) and \(\nu/\mu\to0\). Let
\(A_\nu\sim\chi_\nu^2/\nu\) and
\(D_\mu\sim\chi_\mu^2/\mu\) independently. Then
\[
 \sqrt{\nu/2}\left(\frac{A_\nu}{D_\mu}-1\right)
 =\frac{\sqrt{\nu/2}(A_\nu-1)
       -\sqrt{\nu/2}(D_\mu-1)}{D_\mu}
 \Rightarrow N(0,1),
\]
because the second term in the numerator has variance \(\nu/\mu\to0\)
and \(D_\mu\to_{\rm p}1\). Thus
\[
 q_{1-\alpha}(\nu,\mu)
 =1+\sqrt{2/\nu}\{z_{1-\alpha}+r(\nu,\mu)\},
 \qquad r(\nu,\mu)\to0
\]
along every such deterministic sequence. This sequential statement also
implies \(r(\widehat d_1,\widehat d_2)=o_{\rm p}(1)\) when
\(\widehat d_1\to_{\rm p}\infty\) and
\(\widehat d_1/\widehat d_2\to_{\rm p}0\): every subsequence has
a further subsequence on which these convergences hold almost surely,
and the deterministic statement applies pointwise there. Relative
consistency gives these two conditions and
\(\sqrt{d_1/\widehat d_1}\to_{\rm p}1\). Subtracting the two
quantile expansions and using
\(b/a=\sqrt{d_1/2}\{1+o(1)\}\) completes the proof.

\subsection{Proof of Proposition~\ref{prop:Fchi}}
It suffices to consider arbitrary subsequences. If \(d_1\) is bounded,
pass to a further subsequence with \(d_1\to d\ge1\).
Then \(d_2\to\infty\), and both critical values converge to
\(\chi^2_{d,1-\alpha}/d\); the factor \(b/a\) is bounded.
If \(d_1\to\infty\), then \(d_1/d_2\to0\), so both quantiles equal
\(1+z_{1-\alpha}\sqrt{2/d_1}+o(d_1^{-1/2})\).
Since \(b/a=\sqrt{d_1/2}\{1+o(1)\}\), their scaled difference tends
to zero. For the estimated parameters, relative consistency implies
\(\widehat d_1/\widehat d_2\to_{\rm p}0\) and preserves the same
bounded- and divergent-degree arguments in probability. This proves the
plug-in claim for estimates computed from the same observations.

\subsection{Proof of Theorem~\ref{thm:size}}
Let \(c_n=(b/a)\{q_{1-\alpha}(d_1,d_2)-1\}\), and let \(\widehat c_n\) denote the corresponding estimated threshold. For fixed \(h>0\),
\begin{align*}
 &\{b(F-1)/a>c_n+h\}\setminus\{|\widehat c_n-c_n|>h\}
 \subseteq \{b(F-1)/a>\widehat c_n\},\\
 &\{b(F-1)/a>\widehat c_n\}\subseteq
 \{b(F-1)/a>c_n-h\}\cup\{|\widehat c_n-c_n|>h\}.
\end{align*}
The standardized reference ratio satisfies asymptotic uniform
anti-concentration because it differs from \(U_0^*/a\) by
\(o_{\rm p}(1)\): for fixed \(h>0\),
its concentration over an interval of radius \(h\) is bounded by that
of \(U_0^*/a\) over radius \(2h\), plus the probability that the
perturbation exceeds \(h\). Calibration-failure events contribute
\(o(1)\) and may be omitted from these probability comparisons.
Lemma~\ref{lem:quantile}, Theorem~\ref{thm:transfer}, and uniform
anti-concentration allow the bracket width to vanish by first taking
\(n\to\infty\) and then \(h\downarrow0\). At the deterministic
threshold, the reference rejection probability differs from \(\alpha\)
by at most \(\varepsilon_{\rm WS}\), proving \eqref{eq:size}.

\section{Unknown-mean transfer and alternatives}\label{app:unknownmeans}

\subsection{The effect of estimating the group means}
Throughout this subsection, \(H_0\) holds, \(\bm{\Sigma}_1=\bm{\Sigma}_2=\bm{\Sigma}\), and
\(\bm{\Gamma}_1=\bm{\Gamma}_2=\bm{\Gamma}\). Write
\[
 \bm{M}_i=\frac1{n_i}\sum_j\bm{\epsilon}_{ij}\bm{\epsilon}_{ij}^\top,
 \qquad \bm{e}_i=\svec(\bm{M}_i-\bm{\Sigma}),\qquad
 \bm{v}_i=\svec(\bm{S}_i-\bm{M}_i).
\]
Here \(\bm{M}_i\) uses the population mean, whereas \(\bm{S}_i\) uses the sample
mean. Put \(\bm{Y}_{ij}=\bm{\epsilon}_{ij}\bm{\epsilon}_{ij}^\top-\bm{\Sigma}\) and define,
with sums over ordered distinct indices,
\begin{align}
 A_{3i}&=\frac1{(n_i)_3}\sum_{a,b,c}^{\ne}
              \bm{\epsilon}_{ia}^\top \bm{Y}_{ib}\bm{\epsilon}_{ic},\label{eq:A3bridge}\\
 A_{4i}&=\frac1{(n_i)_4}\sum_{a,b,c,d}^{\ne}
       (\bm{\epsilon}_{ia}^\top\bm{\epsilon}_{ib})
       (\bm{\epsilon}_{ic}^\top\bm{\epsilon}_{id}).\label{eq:A4bridge}
\end{align}
Let \(\bm{\mathcal T}_i\) be the \(p\times p^2\) matrix with
\((a,(b,c))\) entry \(\E(\epsilon_{i1,a}\epsilon_{i1,b}\epsilon_{i1,c})\), and set
\[
 \eta_i=\norm{\bm{\Sigma}^{1/2}\bm{\mathcal T}_i}_F^2,
 \qquad c_\Sigma=\{\tr(\bm{\Sigma}^2)\}^2+2\tr(\bm{\Sigma}^4).
\]

\begin{proposition}[Exact centered remainder]\label{prop:fullbridge}
If \(n_i\ge4\) and \(\E\norm{\bm{\epsilon}_{i1}}^4<\infty\), then
\begin{align}
 U_{\rm full}-U^o
 ={}&-2A_{31}-2A_{32}+A_{41}+A_{42}\notag\\
    &-2\bm{e}_1^\top \bm{v}_2-2\bm{v}_1^\top \bm{e}_2-2\bm{v}_1^\top \bm{v}_2.
       \label{eq:exactbridge}
\end{align}
The seven terms on the right are mutually uncorrelated and are uncorrelated
with \(U^o\). In particular, writing
\(r_n^2=\E(U_{\rm full}-U^o)^2\),
\begin{align}
 r_n^2={}&4\sum_{i=1}^2
       \frac{\tr(\bm{\Gamma}\bm{\Omega}_i^{\,o})+4\eta_i}{(n_i)_3}
       +8c_\Sigma\sum_{i=1}^2\frac1{(n_i)_4}\notag\\
 &+\frac{4\tr(\bm{\Omega}_1^{\,o}\bm{\Gamma})}{n_1n_2(n_2-1)}
  +\frac{4\tr(\bm{\Gamma}\bm{\Omega}_2^{\,o})}{n_1(n_1-1)n_2}
  +\frac{4\tr(\bm{\Gamma}^2)}{(n_1)_2(n_2)_2}.
       \label{eq:bridgevariance}
\end{align}
Moreover, \(0\le\eta_i\le\tfrac12\tr(\bm{\Gamma}\bm{\Omega}_i^{\,o})\).
Consequently, \(\Var(U_{\rm full})=(a^o)^2+r_n^2\).
\end{proposition}

\begin{corollary}[Sufficient full-sample transfer conditions]
\label{cor:fullbridge}
Suppose Assumptions~\ref{ass:sampling} and \ref{ass:fourth} hold for the
oracle features. If
\begin{equation}\label{eq:bridgeconditions}
 \frac{r_n^2}{(a^o)^2}\longrightarrow0,
 \qquad
 \frac{\rho_n}{b^o}\longrightarrow0,
 \qquad
 \rho_n=\sum_{i=1}^2\frac{\tr(\bm{\Gamma})}{n_i(n_i-1)},
\end{equation}
then
\[
 \frac{U_{\rm full}-U^o}{a^o}\to_{\rm p}0,
 \qquad \frac{B_{\rm full}}{b^o}\to_{\rm p}1,
 \qquad \Pp_0(B_{\rm full}\le0)\longrightarrow0.
\]
For balanced sample sizes tending to infinity, a convenient sufficient
condition for the first limit in \eqref{eq:bridgeconditions} is
\begin{equation}\label{eq:GammaOmegacondition}
 \tr(\bm{\Gamma}^2)=O\!\left[
     \tr\{(\bm{\Omega}_1^{\,o})^2\}+\tr\{(\bm{\Omega}_2^{\,o})^2\}\right].
\end{equation}
This gives \(r_n^2/(a^o)^2=O(n_{\min}^{-1})\), without a restriction on
\(p/n_{\min}\).
\end{corollary}

A sufficient condition for the mean-scale approximation can also be
expressed in terms of \(\bm{\Sigma}\). If
\(\{\tr(\bm{\Sigma})\}^2/\tr(\bm{\Sigma}^2)\ge1+\delta\) for a fixed \(\delta>0\),
Jensen's inequality gives
\[
 \tr(\bm{\Omega}_i^{\,o})\ge\{\tr(\bm{\Sigma})\}^2-\tr(\bm{\Sigma}^2),
 \qquad
 \frac{\tr(\bm{\Gamma})}{\tr(\bm{\Omega}_i^{\,o})}\le\frac{2+\delta}{\delta}.
\]
Thus \(\rho_n/b^o=O(n_{\min}^{-1})\). The trace condition used in
Section~\ref{sec:fullblocktraces} for estimating \(J^o\) and \(R^o\)
also implies \eqref{eq:bridgeconditions}.

\subsection{Proofs of Proposition~\ref{prop:fullbridge} and
Corollary~\ref{cor:fullbridge}}
For one group, define
\[
 O_i=\frac1{(n_i)_2}\sum_{a\ne b}
                 (\bm{\epsilon}_{ia}^\top\bm{\epsilon}_{ib})^2.
\]
Expanding the four differences in \eqref{eq:D}, and grouping terms with
two, three and four distinct indices, gives
\[
 D_{n_i}=O_i+2\bm{\theta}^\top \bm{v}_i-2A_{3i}+A_{4i},
 \qquad \bm{\theta}=\svec(\bm{\Sigma}).
\]
Here we used
\(\bm{S}_i-\bm{M}_i=-(n_i)_2^{-1}\sum_{a\ne b}
\bm{\epsilon}_{ia}\bm{\epsilon}_{ib}^\top\).
Since \(U^o=O_1+O_2-2\svec(\bm{M}_1)^\top\svec(\bm{M}_2)\), substitution into
\eqref{eq:Ffull} cancels all terms containing \(\bm{\theta}\) and proves
\eqref{eq:exactbridge}.

To calculate the remainder variance, symmetrize the kernel of
\(A_{3i}\):
\[
 \frac13\{\bm{\epsilon}_1^\top \bm{Y}_2\bm{\epsilon}_3
          +\bm{\epsilon}_1^\top \bm{Y}_3\bm{\epsilon}_2
          +\bm{\epsilon}_2^\top \bm{Y}_1\bm{\epsilon}_3\}.
\]
Integrating any argument gives zero, so this kernel is canonical of
order three. A squared summand has expectation
\(A_i=\E\tr(\bm{Y}_{i1}\bm{\Sigma} \bm{Y}_{i1}\bm{\Sigma})
=\tfrac12\tr(\bm{\Gamma}\bm{\Omega}_i^{\,o})\).
The product of two distinct summands has expectation
\[
 \sum_{a,c}\Sigma_{ac}\sum_{b,d}
  \E(\epsilon_a\epsilon_b\epsilon_d)
  \E(\epsilon_c\epsilon_b\epsilon_d)=\eta_i.
\]
The displayed expression is \(\norm{\bm{\Sigma}^{1/2}\bm{\mathcal T}_i}_F^2\),
and is therefore nonnegative. Cauchy--Schwarz gives \(\eta_i\le A_i\).
The variance formula for a canonical U-statistic now yields
\[
 \Var(A_{3i})
 =\frac{(A_i+2\eta_i)/3}{\binom{n_i}{3}}
 =\frac{\tr(\bm{\Gamma}\bm{\Omega}_i^{\,o})+4\eta_i}{(n_i)_3}.
\]
The symmetric kernel of \(A_{4i}\) averages the three pairings of four
arguments. It is canonical of order four. A squared pairing has expectation
\(\{\tr(\bm{\Sigma}^2)\}^2\), whereas the product of two different pairings has
expectation \(\tr(\bm{\Sigma}^4)\). Hence
\(\Var(A_{4i})=8c_\Sigma/(n_i)_4\).

Proposition~\ref{prop:full} gives
\(\Cov(\bm{e}_i)=\bm{\Omega}_i^{\,o}/n_i\) and
\(\Cov(\bm{v}_i)=\bm{\Gamma}/(n_i)_2\). Independence between groups gives
\begin{align*}
 \Var(\bm{e}_1^\top \bm{v}_2)
 &=\frac{\tr(\bm{\Omega}_1^{\,o}\bm{\Gamma})}{n_1(n_2)_2},&
 \Var(\bm{v}_1^\top \bm{e}_2)
 &=\frac{\tr(\bm{\Gamma}\bm{\Omega}_2^{\,o})}{(n_1)_2n_2},\\
 \Var(\bm{v}_1^\top \bm{v}_2)
 &=\frac{\tr(\bm{\Gamma}^2)}{(n_1)_2(n_2)_2}.
\end{align*}
The corresponding canonical kernels have within-group orders
\[
 (3,0),\ (0,3),\ (4,0),\ (0,4),\ (1,2),\ (2,1),\ (2,2),
\]
respectively.
Distinct patterns are orthogonal: a product of their kernels has an
unmatched observation, and integrating that argument gives zero.
The same reasoning proves orthogonality to the null oracle statistic,
whose patterns are \((2,0),(0,2),(1,1)\). Summing the variances proves
\eqref{eq:bridgevariance} and Proposition~\ref{prop:fullbridge}.

For the corollary, Chebyshev's inequality first gives
\((U_{\rm full}-U^o)/a^o\to_{\rm p}0\). Assumption~\ref{ass:fourth}
implies \(\E\norm{\bm{Z}_{i1}^{\,o}-\bm{\theta}}^4
\le K\{\tr(\bm{\Omega}_i^{\,o})\}^2\); the sample-trace variance identity
therefore gives \(\Var(B^o)/(b^o)^2=O(n_{\min}^{-1})\), so
\(B^o/b^o\to_{\rm p}1\). With \(\bm{w}=\bm{e}_1-\bm{e}_2\) and \(\bm{v}=\bm{v}_1-\bm{v}_2\),
\[
 Q_{\rm full}-Q^o=2\bm{w}^\top \bm{v}+\norm{\bm{v}}^2,
 \qquad
 \E|Q_{\rm full}-Q^o|\le2\sqrt{b^o\rho_n}+\rho_n.
\]
Consequently \((Q_{\rm full}-Q^o)/b^o\to_{\rm p}0\).
Since \(a^o\le2b^o\), the exact identity
\[
 B_{\rm full}-B^o
  =(Q_{\rm full}-Q^o)-(U_{\rm full}-U^o)
\]
proves denominator consistency and the claimed positivity probability.
Finally, \(\eta_i\le\tr(\bm{\Gamma}\bm{\Omega}_i^{\,o})/2\),
\(\tr(\bm{\Gamma}\bm{\Omega}_i^{\,o})\le
\sqrt{\tr(\bm{\Gamma}^2)\tr\{(\bm{\Omega}_i^{\,o})^2\}}\), and
\(c_\Sigma\le3\tr(\bm{\Gamma}^2)/4\). Under balanced sample sizes and
\eqref{eq:GammaOmegacondition}, these bounds give
\(r_n^2/(a^o)^2=O(n_{\min}^{-1})\).

\subsection{Proof of Theorem~\ref{thm:fulltransfer}}
We first verify that
\eqref{eq:blockdominance} also controls the mean correction in
\eqref{eq:bfull}. Put \(T_o=\sum_i\tr(\bm{\Omega}_i^{\,o})\) and
\(G_o=\sum_i\tr(\bm{\Gamma}_i^2)\). Jensen's inequality gives
\[
 \tr(\bm{\Gamma}_i)\le\tr(\bm{\Omega}_i^{\,o})+2\tr(\bm{\Sigma}_i^2)
 \le\tr(\bm{\Omega}_i^{\,o})+\sqrt2\{\tr(\bm{\Gamma}_i^2)\}^{1/2}.
\]
If \(G_o\le C S_o\), summation and \(\sqrt{S_o}\le T_o\) imply
\(\sum_i\tr(\bm{\Gamma}_i)\le(1+2\sqrt C)T_o\). Under balanced sample sizes,
\begin{equation}\label{eq:rhoautomatic}
 \rho_n/b^o=O(n_{\min}^{-1}).
\end{equation}
Thus the denominator condition of Corollary~\ref{cor:fullbridge}
follows from the trace condition used for calibration.

The trace condition, balanced sample sizes,
and \eqref{eq:rhoautomatic} allow Corollary~\ref{cor:fullbridge} to be applied.
The variance bound in Proposition~\ref{prop:fullbridge} yields the rate in
\eqref{eq:FB1}. Hence
\[
 \frac{b^o}{a^o}(F_{\rm full}-1)
 =\frac{U^o/a^o+\op(1)}{B_{\rm full}/b^o}
 =U^o/a^o+\op(1).
\]
The Gaussian comparison and anti-concentration follow from
Theorem~\ref{thm:feature-universality} and
Lemma~\ref{lem:gaussian-anticoncentration}. Proposition~\ref{prop:fullblocktraces}
and denominator consistency establish relative consistency of both estimated
degrees of freedom. Lemma~\ref{lem:quantile} and the deterministic threshold
bracketing argument used for Theorem~\ref{thm:size} prove \eqref{eq:fullsize}.
All calibration-failure events have vanishing probability by relative
consistency. Finally, under \eqref{eq:rho3} for \(\bm{V}^o\),
Proposition~\ref{prop:WS} gives vanishing matching error and hence
asymptotic size \(\alpha\). Corollary~\ref{cor:verifiedlevel} follows from
Proposition~\ref{prop:verifiedmodels} and \eqref{eq:originaldiffusebound}.

\subsection{Proof of Corollary~\ref{cor:unifiedfull}}
For the centered Gaussian reference, \(R^o/J^o\to0\) and
\eqref{eq:centraloraclereference} give
\[
 \frac{U_0^{o,*}}{a^o}
 =\frac{Q_0^{o,*}-b^o}{\sqrt{2J^o}}+\op(1)
 =\sum_{r\ge1}c_{nr}\xi_r+\op(1),
 \qquad \xi_r=(Z_r^2-1)/\sqrt2.
\]
Indeed, the variance of \((B^{o,*}-b^o)/a^o\) is
\(R^o/(J^o+R^o)\to0\), and \(a^o/\sqrt{2J^o}\to1\).
Fatou's lemma gives \(\sum_r c_r^2\le1\), while the ordering implies
\(c_{n,K+1}\le(K+1)^{-1/2}\).
For fixed \(t\) and sufficiently small \(|c|\),
\[
 \log\E e^{itc\xi_1}=-t^2c^2/2+O_t(|c|^3).
\]
Consequently, the logarithm of the characteristic function of the tail
\(\sum_{r>K}c_{nr}\xi_r\) differs from
\(-t^2\sum_{r>K}c_{nr}^2/2\) by at most
\(C_t c_{n,K+1}\sum_{r>K}c_{nr}^2\).
For fixed \(K\), let \(n\to\infty\) in the leading coefficients
and in the tail variance, then let \(K\to\infty\). The resulting
characteristic function is that of \(\zeta\) in
\eqref{eq:unifiedfull}: the limiting leading coordinates form the
chi-square series, and the remaining variance produces the independent
normal component.
The chi-square series converges in \(L^2\), since \(\sum_r c_r^2<\infty\).
If the Gaussian coefficient is positive, its convolution has a continuous
distribution. Otherwise, at least one \(c_r\) is positive, and convolution with
that continuous chi-square component again gives a continuous
distribution. For a finite coefficient sequence \((b_r)\), independence
gives
\[
 \E\left(\sum_r b_r\xi_r\right)^4
 =3\left(\sum_r b_r^2\right)^2+12\sum_r b_r^4
 \le15\left(\sum_r b_r^2\right)^2.
\]
Hence the chi-square series also converges in \(L^4\), and the third
moment of \(\zeta\) is \(2\sqrt2\sum_r c_r^3\). This is strictly
positive whenever some \(c_r>0\), proving that the limit is then
nonnormal.

The standardized full-sample statistic differs from \(U_0^o/a^o\) by
\(\op(1)\) under \(H_0\), by Theorem~\ref{thm:fulltransfer} and its
proof. Gaussian comparison therefore transfers the reference limit to
that statistic. Convergence to the continuous distribution also gives
\(\sup_x|G_n^o(x)-G_\zeta(x)|\to0\).

\subsection{Independent-feature alternatives}
To obtain \eqref{eq:alternative}, write the difference of feature means
as \(\bm{\Delta}+\bm{w}\) and expand its squared norm. Every summand in
\eqref{eq:Udeg} is degenerate in each observation index; its product
with any component of \(\bm{w}\) therefore has expectation zero. This
proves \(\Cov(U_0,\bm{w})=\bm{0}\) and the stated variance formula.

Under the conditions of Theorem~\ref{thm:power}, Chebyshev's inequality gives
\(2\bm{\Delta}^\top \bm{w}/a\to_{\rm p}0\). Since \(h_n=O(1)\),
multiplication by \(b/B=1+\op(1)\) gives
\[
 \frac b a(F-1)=U_0/a+h_n+\op(1).
\]
The threshold bracketing argument of Theorem~\ref{thm:size} now
proves the power formula. Under spectral diffuseness, the centered
Gaussian quadratic form tends to normal, \(R/J\to0\), and
\(t_n\to z_{1-\alpha}\), giving the normal limit.

Under the conditions of Proposition~\ref{prop:consistency},
\eqref{eq:alternative} gives \(U/\delta_n\to_{\rm p}1\) and
\(a/\delta_n\to0\). The assumptions \(B/b\to_{\rm p}1\) and
\(\widehat q-1=\Op(a/b)\) give
\(B(\widehat q-1)=\Op(a)\). Since $\Pp(B>0)\to1$, the rejection inequality is
$U>B(\widehat q-1)$ with probability tending to one. Its right-hand side
is negligible relative to $\delta_n$, so the rejection probability tends
to one.

\subsection{Approximation under alternatives}
\label{sec:alternativebridge}
We now allow \(\bm{\Sigma}_1\ne\bm{\Sigma}_2\). Write
\[
 \bm{\theta}_i=\svec(\bm{\Sigma}_i),\quad \bm{\Delta}=\bm{\theta}_1-\bm{\theta}_2,
 \quad \delta_n=\norm{\bm{\Delta}}^2,
 \quad \bm{e}_i=\svec(\bm{M}_i-\bm{\Sigma}_i),\quad \bm{v}_i=\svec(\bm{S}_i-\bm{M}_i).
\]
Define \(A_{3i},A_{4i}\) as in \eqref{eq:A3bridge}--\eqref{eq:A4bridge},
now with \(\bm{Y}_{ij}=\bm{\epsilon}_{ij}\bm{\epsilon}_{ij}^\top-\bm{\Sigma}_i\).
All reference parameters below are evaluated under the sequence of
alternatives.

\begin{lemma}[Alternative bridge and stable calibration]
\label{lem:alternativebridge}
The exact decomposition is
\begin{align}
 U_{\rm full}-U^o
  &=\mathcal R_{0n}+2\bm{\Delta}^\top(\bm{v}_1-\bm{v}_2),\label{eq:altbridge}\\
 \mathcal R_{0n}
  &=-2A_{31}-2A_{32}+A_{41}+A_{42}
    -2\bm{e}_1^\top \bm{v}_2-2\bm{v}_1^\top \bm{e}_2-2\bm{v}_1^\top \bm{v}_2.\notag
\end{align}
The two terms on the right of \eqref{eq:altbridge} have mean zero and are
uncorrelated, with
\begin{equation}\label{eq:altsignalvar}
 \Var\{2\bm{\Delta}^\top(\bm{v}_1-\bm{v}_2)\}
  =4\sum_{i=1}^2\frac{\bm{\Delta}^\top\bm{\Gamma}_i\bm{\Delta}}{(n_i)_2}.
\end{equation}
Suppose the oracle sampling and moment assumptions and
\eqref{eq:blockdominance} hold uniformly along the sequence, without
imposing \(H_0\). There is a constant \(C\), independent of the dimension
and signal, such that, for all sufficiently large sample sizes,
\begin{equation}\label{eq:altbridgebound}
 \frac{\E(U_{\rm full}-U^o)^2}{(a^o)^2}
 \le \frac{C}{n_{\min}}\left(1+\frac{\delta_n}{a^o}\right).
\end{equation}
Moreover,
\begin{equation}\label{eq:altcalibration}
 B_{\rm full}/b^o\to_{\rm p}1,
 \qquad \widehat d_{j,\rm full}/d_j^o\to_{\rm p}1,
 \quad j=1,2,
\end{equation}
and the probability of calibration failure tends to zero. These calibration conclusions hold for arbitrary signal size \(\delta_n\).
\end{lemma}
For \(\delta_n=O(a^o)\), the standardized approximation error is
\(\Op(n_{\min}^{-1/2})\), as under the null, without an additional
condition on \(\bm{\Delta}^\top\bm{\Gamma}_i\bm{\Delta}\). The signal-dependent term may share order-two components with \(U^o\);
their covariance need not vanish under alternatives.

\subsection{Proof of Lemma~\ref{lem:alternativebridge}}
For each group, the expansion used for Proposition~\ref{prop:fullbridge}
holds with its own \(\bm{\theta}_i\):
\[
 D_{n_i}=O_i+2\bm{\theta}_i^\top \bm{v}_i-2A_{3i}+A_{4i}.
\]
Substituting \(\svec(\bm{M}_i)=\bm{\theta}_i+\bm{e}_i\) into the two cross-products
leaves \(2(\bm{\theta}_1-\bm{\theta}_2)^\top(\bm{v}_1-\bm{v}_2)\) and the seven centered
terms in \eqref{eq:altbridge}. This proves the identity.
The seven canonical order patterns of \(\mathcal R_{0n}\) remain
\((3,0),(0,3),(4,0),(0,4),(1,2),(2,1),(2,2)\).
The additional term has patterns \((2,0),(0,2)\), so it is orthogonal
to \(\mathcal R_{0n}\). Independence of the two groups and
\(\Cov(\bm{v}_i)=\bm{\Gamma}_i/(n_i)_2\) prove \eqref{eq:altsignalvar}.

Write \(N=n_{\min}\). Balanced sample sizes and positive semidefiniteness
give \((a^o)^2\asymp S_o/N^2\). Under \eqref{eq:blockdominance},
\(\norm{\bm{\Omega}_i^{\,o}}_F\le\sqrt{S_o}\) and
\(\norm{\bm{\Gamma}_i}_F\le C\sqrt{S_o}\).
The variance computation for Proposition~\ref{prop:fullbridge} applies
groupwise with \(\bm{\Sigma}_i,\bm{\Gamma}_i\). Its third-order terms are bounded by
\(C\tr(\bm{\Gamma}_i\bm{\Omega}_i^{\,o})/N^3\), and its fourth-order terms by
\(C\tr(\bm{\Gamma}_i^2)/N^4\). The cross-group terms are bounded by
\[
 \frac{C}{N^3}\{
  \tr(\bm{\Omega}_1^{\,o}\bm{\Gamma}_2)+\tr(\bm{\Gamma}_1\bm{\Omega}_2^{\,o})\}
 +\frac{C}{N^4}\tr(\bm{\Gamma}_1\bm{\Gamma}_2).
\]
Trace Cauchy--Schwarz thus gives
\(\E\mathcal R_{0n}^2\le C S_o/N^3\le C(a^o)^2/N\).
Also,
\[
 \sum_i\frac{\bm{\Delta}^\top\bm{\Gamma}_i\bm{\Delta}}{(n_i)_2}
 \le\frac{C\delta_n\sqrt{S_o}}{N^2}
 \le\frac{C\delta_n a^o}{N}.
\]
Combining these bounds with orthogonality proves \eqref{eq:altbridgebound}.

Denominator consistency follows from a within-group cancellation of
the signal terms. Let
\(B_i^o=\norm{\svec(\bm{M}_i)}^2-O_i\). Then
\(B^o=B_1^o+B_2^o\), and
\begin{equation}\label{eq:altdenominatoridentity}
 \widehat b_i-B_i^o
  =2\bm{e}_i^\top \bm{v}_i+\norm{\bm{v}_i}^2+2A_{3i}-A_{4i}.
\end{equation}
Set \(\rho_n=\sum_i\tr(\bm{\Gamma}_i)/(n_i)_2\), now allowing distinct
\(\bm{\Gamma}_i\). The trace argument preceding \eqref{eq:rhoautomatic}
does not use \(H_0\), so \(\rho_n/b^o=O(N^{-1})\).
Cauchy--Schwarz, \(\sum_i\E\norm{\bm{e}_i}^2=b^o\), and the preceding
variance bounds imply
\[
 \E|B_{\rm full}-B^o|
 \le 2\sqrt{b^o\rho_n}+\rho_n+C a^o/\sqrt N.
\]
Since \(a^o\le2b^o\), this bound is \(O(b^o/\sqrt N)\).
The sample-trace variance identity applied to the centered features
gives \(B^o/b^o\to_{\rm p}1\) under alternatives. Thus
\(B_{\rm full}/b^o\to_{\rm p}1\) regardless of the signal size.
Proposition~\ref{prop:fullblocktraces} likewise does not require \(H_0\).
Its relative consistency, together with that of \(B_{\rm full}\), proves
\eqref{eq:altcalibration} and the vanishing failure probability.

\subsection{Proofs of Theorem~\ref{thm:fullpower} and
Proposition~\ref{prop:fullconsistency}}
The oracle identity under alternatives is
\[
 U^o=U_0^o+2\bm{\Delta}^\top(\bm{e}_1-\bm{e}_2)+\delta_n,
 \qquad
 \Var(U^o)=(a^o)^2+4\bm{\Delta}^\top \bm{V}^o\bm{\Delta}.
\]
Under \eqref{eq:fulllocal}, the linear term divided by \(a^o\)
converges to zero in probability. Lemma~\ref{lem:alternativebridge} gives
\((U_{\rm full}-U^o)/a^o=\op(1)\) and denominator consistency. Hence
\[
 \frac{b^o}{a^o}(F_{\rm full}-1)
 =U_0^o/a^o+h_n+\op(1).
\]
Gaussian comparison and uniform anti-concentration apply to the centered
features along the alternative. Lemma~\ref{lem:quantile}, combined with
\eqref{eq:altcalibration}, removes the random threshold by deterministic
bracketing. These statements prove \eqref{eq:fullpowersequence}, with the
failure events contributing only \(o(1)\).
The argument in the proof of Corollary~\ref{cor:unifiedfull} concerns
the central reference \(Q_0^{o,*}\), so it applies to the feature
covariances of the current alternative without imposing \(H_0\).
If \(d_1^o\to d<\infty\), then
\(b^o/a^o\to\sqrt{d/2}\) and \(d_2^o\to\infty\).
If \(d_1^o\to\infty\), the standardized reference quantile tends to
\(z_{1-\alpha}\), as in Lemma~\ref{lem:quantile}.
Thus \(t_n^o\to t_\alpha(d)\); continuity of \(G_\zeta\) proves
\eqref{eq:fullpowerlimit}. Under the stated diffuseness condition,
\(\max_r c_{nr}^3\le\sum_r c_{nr}^3\to0\), and
\(d_1^o\to\infty\) by Proposition~\ref{prop:WS}.
Then \(\zeta\) is standard normal and \eqref{eq:fullnormalpower} follows.

For consistency, \eqref{eq:fullconsistency} gives
\(a^o/\delta_n\to0\) and \(U^o/\delta_n\to_{\rm p}1\).
By \eqref{eq:altbridgebound},
\[
 \E\left(\frac{U_{\rm full}-U^o}{\delta_n}\right)^2
 \le\frac C N\left\{
       \left(\frac{a^o}{\delta_n}\right)^2
       +\frac{a^o}{\delta_n}\right\}\longrightarrow0.
\]
Therefore \(U_{\rm full}/\delta_n\to_{\rm p}1\).
The bounded- and divergent-degree arguments of Lemma~\ref{lem:quantile}
also show \(t_n^o=O(1)\). Consequently,
\(\widehat q-1=\Op(a^o/b^o)\) for the actual estimated critical value,
and \(B_{\rm full}(\widehat q-1)=\Op(a^o)=\op(\delta_n)\).
The rejection inequality
\(U_{\rm full}>B_{\rm full}(\widehat q-1)\) therefore holds with
probability tending to one. Together with the vanishing calibration-failure
probability, this proves consistency.

\subsection{A spike-aligned local alternative}\label{app:powerexample}
A signal aligned with a strong spike can behave differently. Let
\(n_1=n_2=n\), assume Gaussian sampling, and take
\[
 \bm{\Sigma}_1=(1-\eta)\bm{I}_p+\eta p\bm{v}\bm{v}^\top,\qquad
 \bm{\Sigma}_2=\bm{\Sigma}_1+\kappa p n^{-1/2}\bm{v}\bm{v}^\top,
 \quad \|\bm{v}\|=1,
\]
where \(0<\eta<1\) and \(\kappa>0\) are fixed. As \(p,n\to\infty\),
the leading eigenvalue of \(\bm{V}^o\), in direction \(\svec(\bm{v}\bm{v}^\top)\),
is asymptotic to \(4\eta^2p^2/n\), while the sum of squares of the
remaining eigenvalues is negligible relative to its square. Hence
\[
 a^o\sim4\sqrt2\eta^2p^2/n,\qquad
 \frac{\delta_n}{a^o}\longrightarrow
       \frac{\kappa^2}{4\sqrt2\eta^2},\qquad
 \frac{\bm{\Delta}^\top \bm{V}^o\bm{\Delta}}{(a^o)^2}
       \longrightarrow\frac{\kappa^2}{8\eta^2}>0.
\]
Thus a bounded standardized signal can coexist with a nonnegligible
linear term. Formula \eqref{eq:fullpowersequence} does not cover this
spike-aligned alternative.

\section{Trace estimation and sufficient moment conditions}\label{app:momentconditions}
\subsection{Proof of Proposition~\ref{prop:fullblocktraces}}
The independent, centered pair contrasts \(\bm{a},\bm{b}\) have covariance
\(\bm{\Sigma}\). Equation~\eqref{eq:paircov} gives
\(\Cov\{\svec(\bm{a}\bm{a}^\top)\}=(\bm{\Omega}^o+\bm{\Gamma})/2\). Taking the
difference of two independent copies gives the first identity in
\eqref{eq:fourblockmean}. For symmetric matrices \(\bm{A},\bm{B}\),
\[
 \E\{(\sqrt2\bm{a}^\top \bm{A} \bm{b})(\sqrt2\bm{a}^\top \bm{B} \bm{b})\}
       =2\tr(\bm{A}\bm{\Sigma} \bm{B}\bm{\Sigma}).
\]
This proves the other identities in \eqref{eq:fourblockmean}.
Products of independent blocks therefore give unbiased estimates in
\eqref{eq:fullshat}--\eqref{eq:fullchat}.

For consistency, first derive a linear fourth-moment bound from the
quadratic-feature assumption. Let \(W=(\bm{t}^\top\bm{\epsilon})^2\), \(m=\E W\), and
\(\sigma_W^2=\Var(W)\). When \(\sigma_W>0\), H\"older's inequality and
the feature assumption give
\[
 \E|W-m|\ge
  \frac{\{\E(W-m)^2\}^{3/2}}{\{\E(W-m)^4\}^{1/2}}
  \ge\frac{\sigma_W}{\sqrt K}.
\]
Since \(W\ge0\), \(\E|W-m|=2\E(m-W)_+\le2m\). Consequently,
\begin{equation}\label{eq:originallinear4}
 \E(\bm{t}^\top\bm{\epsilon})^4\le(1+4K)(\bm{t}^\top\bm{\Sigma} \bm{t})^2.
\end{equation}
The same inequality is immediate when \(\sigma_W=0\); a dimension-free
constant also applies to the independent pair contrasts \(\bm{a},\bm{b}\).

Put \(\bm{\Lambda}_i=\bm{\Omega}_i^{\,o}+\bm{\Gamma}_i\). Expansion of
\[
 \bm{a}^\top \bm{A} \bm{a}-\tr(\bm{A}\bm{\Sigma})
 =\tfrac12\sum_{j=1}^2
       \{\bm{\epsilon}_j^\top \bm{A}\bm{\epsilon}_j-\tr(\bm{A}\bm{\Sigma})\}
       -\bm{\epsilon}_1^\top \bm{A}\bm{\epsilon}_2
\]
and conditioning in the bilinear term, using \eqref{eq:originallinear4},
show that a constant \(C_K\), independent of dimension, satisfies
\begin{equation}\label{eq:blockfourth}
 \E(\bm{t}^\top \bm{v}_{ir})^4+\E(\bm{t}^\top \bm{w}_{ir})^4
       \le C_K(\bm{t}^\top\bm{\Lambda}_i \bm{t})^2.
\end{equation}
Specifically, the bilinear fourth moment is at most
\(C\{\tr(\bm{A}\bm{\Sigma} \bm{A}\bm{\Sigma})\}^2\), whereas the centered quadratic
terms are bounded by \(K\{\svec(\bm{A})^\top\bm{\Omega}_i^{\,o}\svec(\bm{A})\}^2\).
Their sum is controlled by \(\bm{\Lambda}_i\), as required.

For a positive semidefinite matrix \(\bm{H}\), spectral decomposition and
Cauchy--Schwarz in \eqref{eq:blockfourth} imply
\[
 \E(\bm{v}_{ir}^\top \bm{H} \bm{v}_{ir})^2+
 \E(\bm{w}_{ir}^\top \bm{H} \bm{w}_{ir})^2
 \le C_K\{\tr(\bm{H}\bm{\Lambda}_i)\}^2.
\]
Conditioning on one of two independent blocks likewise bounds every
fourth power of a block inner product by
\(C_K\{\tr(\bm{\Lambda}_i\bm{\Lambda}_j)\}^2\). Expanding
\(\tr(\bm{K}_{ir}\bm{K}_{js})\) into its four terms therefore gives
\begin{align*}
 \E\{\tr(\bm{K}_{ir}\bm{K}_{js})\}^2
 &\le C_K\{\tr(\bm{\Lambda}_i\bm{\Lambda}_j)\}^2,\\
 \E\{\tr(\bm{K}_{ir}\bm{\Omega}_j^{\,o})\}^2
 &\le C_K\{\tr(\bm{\Lambda}_i\bm{\Omega}_j^{\,o})\}^2.
\end{align*}
The one-sample and two-sample U-statistic variance decompositions yield
\begin{align}
 \E(\widetilde s_i-s_i^o)^2
 &\le C_K\left[
  \frac{\{\tr(\bm{\Lambda}_i\bm{\Omega}_i^{\,o})\}^2}{L_i}
  +\frac{\{\tr(\bm{\Lambda}_i^2)\}^2}{L_i^2}\right],\label{eq:blockvar1}\\
 \E(\widetilde c_{12}-c_{12}^o)^2
 &\le C_K\left[
  \frac{\{\tr(\bm{\Lambda}_1\bm{\Omega}_2^{\,o})\}^2}{L_1}
  +\frac{\{\tr(\bm{\Lambda}_2\bm{\Omega}_1^{\,o})\}^2}{L_2}
  +\frac{\{\tr(\bm{\Lambda}_1\bm{\Lambda}_2)\}^2}{L_1L_2}\right].
  \label{eq:blockvar2}
\end{align}
Here \(L_i(L_i-1)\asymp L_i^2\) for \(L_i\ge2\).

Condition \eqref{eq:blockdominance} implies
\(\tr(\bm{\Lambda}_i^2)\le2s_i^o+2\tr(\bm{\Gamma}_i^2)=O(S_o)\).
By the trace Cauchy--Schwarz inequality, each right-hand side of
\eqref{eq:blockvar1}--\eqref{eq:blockvar2} is
\(O(S_o^2/L_{\min})\). Projection onto \([0,\infty)\) cannot increase
squared error from a nonnegative target. Thus the same absolute-error
bounds apply to \(\widehat s_i\) and \(\widehat c_{12}\), including
when an individual target is zero. Balanced sample sizes give
\(J^o\ge S_o/n_{\max}^2\) and
\(R^o\ge S_o/\{n_{\max}^2(n_{\max}-1)\}\).
Combining these bounds with \eqref{eq:fullblockJR} proves the asserted
relative mean-square rates. Chebyshev's inequality gives the positivity
probability statements. Jensen's inequality gives the same mean-square bound after averaging
over prespecified partitions.

\subsection{Verification of the factor-model conditions}
\begin{proposition}[Quadratic-feature moment and trace bounds]
\label{prop:verifiedmodels}
Under \eqref{eq:verifiedfactormodel}--\eqref{eq:latentmoments}, there is a
constant $K=K(M,c)<\infty$ such that, for every symmetric matrix $\bm{A}$,
\begin{equation}\label{eq:verifiedquadfourth}
 \E\left\{\bm{\epsilon}_i^\top \bm{A}\bm{\epsilon}_i-\tr(\bm{A}\bm{\Sigma}_i)\right\}^4
 \le K\left[\Var(\bm{\epsilon}_i^\top \bm{A}\bm{\epsilon}_i)\right]^2.
\end{equation}
Consequently, Assumption~\ref{ass:fourth} holds for the oracle features
$\svec(\bm{\epsilon}_i\bm{\epsilon}_i^\top)$, uniformly in the dimension. Moreover,
\begin{equation}\label{eq:Omegasandwich}
 c_0\bm{\Gamma}_i\preceq\bm{\Omega}_i^{\,o}\preceq C_0\bm{\Gamma}_i,
 \qquad
 c_0=\min\{1,c/2\},\quad
 C_0=\max\{1,(\sqrt M-1)/2\}.
\end{equation}
In particular, for every nonzero relevant covariance matrix,
\begin{equation}\label{eq:verifiedscalebounds}
 \frac{\tr(\bm{\Gamma}_i)}{\tr(\bm{\Omega}_i^{\,o})}\le c_0^{-1},
 \qquad
 \frac{\tr(\bm{\Gamma}_i^2)}{\tr\{(\bm{\Omega}_i^{\,o})^2\}}
 \le c_0^{-2}.
\end{equation}
The paired quadratic features also satisfy Assumption~\ref{ass:fourth}.
For those features alone, the lower bound in \eqref{eq:latentmoments} is
unnecessary: independent standardized latent coordinates with uniformly
bounded eighth moments suffice.
\end{proposition}

Under covariance equality, \eqref{eq:verifiedscalebounds} yields
$\bm{\Gamma}_1=\bm{\Gamma}_2=\bm{\Gamma}$,
$\rho_n/b^o=O(n_{\min}^{-1})$, and
$\tr(\bm{\Gamma}^2)=O[\tr\{(\bm{\Omega}_1^{\,o})^2\}
+\tr\{(\bm{\Omega}_2^{\,o})^2\}]$.
It therefore verifies the transfer and calibration conditions without
restricting $p/n_{\min}$ or the condition number of $\bm{\Sigma}_i$.

Table~\ref{tab:modelmoments} gives the standardized innovation moments
for the distributions used in the simulations.

\begin{table}[!ht]
\centering
\small
\caption{Standardized innovation moments used to verify the assumptions.}
\label{tab:modelmoments}
\begin{tabular}{@{}lrrr@{}}
\toprule
Moment & \(N(0,1)\) & \(\sqrt{8/10}\,t_{10}\) & \((\chi_1^2-1)/\sqrt2\)\\
\midrule
\(\E z^4\) & 3 & 4 & 15\\
\(\E z^8\) & 105 & 1120 & 74417\\
\bottomrule
\end{tabular}
\end{table}

All three innovation distributions satisfy \eqref{eq:latentmoments},
although the moment bounds are much larger for the skewed distribution.
The two groups may also have different innovation distributions under
the covariance null hypothesis.

\subsection{Proof of Proposition~\ref{prop:verifiedmodels}}
Suppress the group index. In this proof only, write
$\bm{B}=\bm{\Sigma}^{1/2}\bm{A}\bm{\Sigma}^{1/2}$ for the matrix of the quadratic form.
Let $\mu_{4j}=\E z_j^4$ and decompose
\[
 Q=\bm{z}^\top \bm{B}\bm{z}-\tr(\bm{B})=D+O,
 \qquad
 D=\sum_j b_{jj}(z_j^2-1),\quad
 O=\sum_{j\ne k}b_{jk}z_jz_k.
\]
Independence and centering give the exact variance formula
\begin{equation}\label{eq:factorquadraticvariance}
 \Var(Q)=2\sum_{j\ne k}b_{jk}^2
          +\sum_j(\mu_{4j}-1)b_{jj}^2
 \ge \min\{2,c\}\tr(\bm{B}^2).
\end{equation}
Each diagonal--off-diagonal product contains an independent centered
coordinate to the first power, so these covariance terms vanish even
for skewed latent coordinates.

For the diagonal part, independence gives
\[
 \E D^4
 =\sum_j b_{jj}^4\E(z_j^2-1)^4
  +6\sum_{j<k}b_{jj}^2b_{kk}^2(\mu_{4j}-1)(\mu_{4k}-1).
\]
Since $\E(z_j^2-1)^4\le8(M+1)$ and
$\mu_{4j}\le\sqrt M$,
\begin{equation}\label{eq:diagfourthbound}
 \E D^4\le8(M+1)\left(\sum_j b_{jj}^2\right)^2.
\end{equation}

For the off-diagonal part, introduce independent Bernoulli$(1/2)$
variables $\delta_j$, independent of $\bm{z}$, and put
\[
 T_\delta=\sum_{j,k}\delta_j(1-\delta_k)b_{jk}z_jz_k.
\]
Then $O=4\E_\delta(T_\delta\mid \bm{z})$, and Jensen's inequality gives
$\E O^4\le256\E_\delta\E T_\delta^4$. For fixed $\delta$, let
$I=\{j:\delta_j=1\}$ and $J=I^c$; then
$T_\delta=\bm{z}_I^\top \bm{B}_{IJ}\bm{z}_J$ has independent arguments.
Every linear form in standardized independent coordinates obeys
\[
 \E(\bm{a}^\top \bm{z})^4
 =3\|\bm{a}\|^4+\sum_j(\mu_{4j}-3)a_j^4
 \le K_4\|\bm{a}\|^4,
 \qquad K_4=\max\{3,\sqrt M\}.
\]
Conditional expectation, followed by Cauchy--Schwarz applied to pairs
of linear forms in $\bm{z}_I$, therefore yields
\[
 \E T_\delta^4
 \le K_4\E\|\bm{B}_{IJ}^\top \bm{z}_I\|^4
 \le K_4^2\|\bm{B}_{IJ}\|_F^4.
\]
Consequently,
\begin{equation}\label{eq:offdiagfourthbound}
 \E O^4\le256K_4^2
          \left(\sum_{j\ne k}b_{jk}^2\right)^2.
\end{equation}
Combining \eqref{eq:diagfourthbound}--\eqref{eq:offdiagfourthbound} with
$(D+O)^4\le8(D^4+O^4)$ proves
$\E Q^4\le C_M\{\tr(\bm{B}^2)\}^2$, where one admissible constant is
\[
 C_M=8\max\{8(M+1),256K_4^2\}.
\]
Equation~\eqref{eq:factorquadraticvariance} gives
\eqref{eq:verifiedquadfourth} with
$K=C_M/\min\{2,c\}^2$.

The Gaussian quadratic-form variance is $2\tr(\bm{B}^2)$. Comparing this
quantity with \eqref{eq:factorquadraticvariance}, from above and below,
proves \eqref{eq:Omegasandwich}. Eigenvalue monotonicity for positive
semidefinite matrices then gives \eqref{eq:verifiedscalebounds}.

For a paired feature, the latent coordinates become
$w_j=(z_j-z'_j)/\sqrt2$, where $\bm{z}'$ is an independent copy.
They remain independent, centered, and standardized, and satisfy
\[
 \E w_j^4=(\mu_{4j}+3)/2,
 \qquad \Var(w_j^2)=(\mu_{4j}+1)/2\ge1,
 \qquad \E w_j^8\le16M.
\]
Applying the preceding argument with $c=1$ and $M$ replaced by $16M$
proves the assertion for paired features.

\clearpage
\section{Additional calibration and partition results}\label{app:earlierstudies}

We examine trace estimation, partition averaging and reference accuracy
at nominal level $\alpha=0.05$. F denotes the four-partition procedure and
$F_{20}$ the pooled procedure in Section~\ref{sec:releasealgorithm}.

Observations are generated from
\[
 \bm{X}_{ij}=\bm{\Sigma}_0^{1/2}\bm{z}_{ij},\qquad
 j=1,\ldots,n,\quad i=1,2,
\]
where $\bm{\Sigma}_0^{1/2}$ is the symmetric positive square root. Both groups
have sample size $n$, and their means are estimated in every test. The
coordinates of $\bm{z}_{ij}$ are independent and standardized. We use two
innovation models:
\begin{description}
\item[Model 1.] Standard normal coordinates in both groups.
\item[Model 2.] Standard normal coordinates in group~1 and
$\sqrt{8/10}\,t_{10}$ coordinates in group~2.
\end{description}
Both models satisfy the covariance null. In Model~2, the two groups have
different quadratic-feature covariances, and the population trace parameters
include their respective fourth moments. The covariance matrices in
Table~\ref{tab:simulation_covariance} have trace $p$; I and A have diffuse
spectra, while S has a nonvanishing spike.

\begin{table}[!ht]
\centering\small
\caption{Covariance structures used in the simulations.}
\label{tab:simulation_covariance}
\begin{tabular*}{\textwidth}{@{\extracolsep{\fill}}cll@{}}
\toprule
Code & Covariance matrix $\bm{\Sigma}_0$ & Structure \\
\midrule
I & $\bm{I}_p$ & Independent coordinates \\
A & $(0.6^{|r-s|})_{r,s=1}^p$ & Autoregressive correlation \\
S & $0.7\bm{I}_p+0.3\mathbf1_p\mathbf1_p^\top$ & A nonvanishing spectral spike \\
\bottomrule
\end{tabular*}
\end{table}

The comparison includes the normal-calibrated test of
\citet{LiChen2012}, denoted by LC, and a chi-square calibration, denoted
by Chi. LC applies the standard normal critical value to
$U_{\mathrm{full}}$ divided by $2D_{n_1}/n_2+2D_{n_2}/n_1$.
Chi uses the same ratio and trace estimates as F, with critical value
$\chi^2_{\widehat d_1,1-\alpha}/\widehat d_1$. The F--Chi comparison isolates the denominator correction.
A population-parameter benchmark, $F_{\mathrm o}$, uses the observed
ratio $Q_{\mathrm{full}}/B_{\mathrm{full}}$ and the population degrees
of freedom. It provides a benchmark for assessing trace-estimation error.

Failed calibrations count as nonrejections. F and Chi require valid trace
estimates and an observed ratio; $F_{\mathrm o}$ requires only the latter.
For $N$ independent datasets, the Monte Carlo standard error of an
estimated rejection probability $\widehat\pi$ is
$\{\widehat\pi(1-\widehat\pi)/N\}^{1/2}$. Comparisons on the same datasets
use paired standard errors. The ancillary summaries also provide 95\%
Wilson confidence intervals.

\subsection{Empirical size and trace estimation}

For Model~1, we use $p\in\{50,200\}$, $n\in\{20,50\}$, and all three
covariance structures, with $10{,}000$ replications per setting.
For Model~2, we use $p=50$, $n\in\{20,50\}$, and I or S, with
$5{,}000$ replications per setting. These designs comprise $140{,}000$
independent datasets. The Monte Carlo standard errors at a rejection
probability of 5\% are 0.218 and 0.308 percentage points, respectively.
Following \citet{WangZhuZhang2024}, we summarize size accuracy over the
$K$ settings specified in Table~\ref{tab:calibration_baseline} by the
average relative error
\[
 \operatorname{ARE}=\frac{100}{K}\sum_{k=1}^{K}
 \frac{|\widehat\pi_k-\alpha|}{\alpha}.
\]

\begin{table}[!ht]
\centering\small
\setlength{\tabcolsep}{3.5pt}
\caption{Empirical sizes (\%) at the 5\% nominal level. Each ARE row uses
the four settings per covariance in Model~1 or the two in Model~2.
A dash denotes a combination not included in the study.}
\label{tab:calibration_baseline}
\begin{tabular*}{\textwidth}{@{\extracolsep{\fill}}rr*{9}{r}@{}}
\toprule
& & \multicolumn{3}{c}{I} & \multicolumn{3}{c}{A} & \multicolumn{3}{c}{S} \\
\cmidrule(lr){3-5}\cmidrule(lr){6-8}\cmidrule(l){9-11}
$p$ & $n$ & LC & Chi & F & LC & Chi & F & LC & Chi & F \\
\midrule
\multicolumn{11}{l}{Model~1} \\
50 & 20 & 4.63 & 5.98 & 5.65 & 6.14 & 6.80 & 6.47 & 9.08 & 8.65 & 8.44 \\
50 & 50 & 5.47 & 6.71 & 6.61 & 6.08 & 6.61 & 6.57 & 9.45 & 8.33 & 8.25 \\
200 & 20 & 4.89 & 6.10 & 5.92 & 5.13 & 6.19 & 5.95 & 10.08 & 9.05 & 8.84 \\
200 & 50 & 4.99 & 6.09 & 5.96 & 5.31 & 6.30 & 6.17 & 9.82 & 8.20 & 8.09 \\
\multicolumn{2}{l}{ARE (\%)} & 4.80 & 24.40 & 20.70 & 13.30 & 29.50 & 25.80 & 92.15 & 71.15 & 68.10 \\
\midrule
\multicolumn{11}{l}{Model~2} \\
50 & 20 & 5.30 & 6.94 & 6.72 & -- & -- & -- & 8.50 & 8.36 & 8.14 \\
50 & 50 & 5.56 & 6.18 & 6.10 & -- & -- & -- & 9.52 & 8.52 & 8.48 \\
\multicolumn{2}{l}{ARE (\%)} & 8.60 & 31.20 & 28.20 & -- & -- & -- & 80.20 & 68.80 & 66.20 \\
\bottomrule\end{tabular*}\end{table}

For Model~1, F has empirical sizes of 5.65--8.84\%. The overall AREs
of LC, Chi and F are 36.75\%, 41.68\% and 38.20\%, respectively, compared
with 44.40\%, 50.00\% and 47.20\% for Model~2. LC performs better for I
but has sizes of 9.08--10.08\% for S under Model~1, consistent with the
effect of spectral concentration on normal calibration. F is less liberal
than Chi by 0.04--0.33 percentage points. The largest difference occurs
for Model~1 with I and $(p,n)=(50,20)$ and has paired standard error
0.057 percentage points. The denominator correction is therefore small
in these settings.

To examine the approximation parameters, Table~\ref{tab:calibration_diagnostics}
reports the median degrees of freedom among valid calibrations and the
relative root mean squared error,
\[
 \operatorname{RRMSE}(\widehat J)
 =\left\{\frac1N\sum_{b=1}^{N}
       \left(\frac{\widehat J^{(b)}}{J}-1\right)^2\right\}^{1/2},
\]
computed over all replications. S has smaller effective degrees of freedom
and more variable trace estimates. Under Model~1 with $p=50$, increasing
$n$ from 20 to 50 reduces RRMSE from 0.888 to 0.396 for I and from 3.143
to 1.604 for S. The baseline failure rate reaches 14.73\%, whereas LC has
no failures.

\begin{table}[!ht]
\centering\small
\setlength{\tabcolsep}{4pt}
\caption{Calibration diagnostics for the four-partition baseline. The two degrees of
freedom are medians conditional on a valid F calibration; the second is
reported in thousands. $F_{\mathrm o}$ size and failure rates are percentages.
The relative root mean squared error of $\widehat J$ includes every replication.}
\label{tab:calibration_diagnostics}
\begin{tabular*}{\textwidth}{@{\extracolsep{\fill}}rrcrrrrr@{}}
\toprule
$p$ & $n$ & $\bm{\Sigma}_0$ & $\widehat d_1$ & $\widehat d_2/10^3$
& $F_{\mathrm o}$ & Failure & RRMSE$(\widehat J)$ \\
\midrule
\multicolumn{8}{l}{Model~1} \\
50 & 20 & I & 1228.5 & 46.698 & 4.65 & 8.71 & 0.888 \\
50 & 20 & A & 293.9 & 11.059 & 5.34 & 9.52 & 1.090 \\
50 & 20 & S & 54.9 & 2.197 & 4.81 & 14.26 & 3.143 \\
50 & 50 & I & 1330.1 & 130.896 & 5.30 & 0.84 & 0.396 \\
50 & 50 & A & 309.1 & 30.671 & 5.04 & 1.48 & 0.502 \\
50 & 50 & S & 44.3 & 4.464 & 5.52 & 5.99 & 1.604 \\
200 & 20 & I & 18945.7 & 701.735 & 5.03 & 9.04 & 0.860 \\
200 & 20 & A & 4284.7 & 160.153 & 5.16 & 9.57 & 0.898 \\
200 & 20 & S & 68.6 & 2.718 & 5.30 & 14.73 & 2.918 \\
200 & 50 & I & 20777.2 & 2053.208 & 5.02 & 0.69 & 0.370 \\
200 & 50 & A & 4660.3 & 459.020 & 5.11 & 0.83 & 0.397 \\
200 & 50 & S & 53.9 & 5.505 & 5.81 & 8.05 & 1.766 \\
\midrule
\multicolumn{8}{l}{Model~2} \\
50 & 20 & I & 1251.3 & 47.400 & 5.16 & 10.08 & 0.891 \\
50 & 20 & S & 53.6 & 2.157 & 4.58 & 13.76 & 3.318 \\
50 & 50 & I & 1343.9 & 132.546 & 5.22 & 0.80 & 0.400 \\
50 & 50 & S & 44.5 & 4.573 & 5.06 & 6.74 & 1.561 \\
\bottomrule\end{tabular*}\end{table}

The population-parameter benchmark has sizes of 4.65--5.81\% under
Model~1, suggesting a substantial contribution from trace-estimation error.
It still includes the mean correction
$\E(B_{\mathrm{full}})-b^o=2\tr(\bm{\Gamma})/\{n(n-1)\}$ and the
moment-matching error examined below.

\subsection{Partition averaging and stability}

A separate study evaluates $F_{20}$ on $3{,}000$ datasets in each of six
settings at $p=50$. These $18{,}000$ datasets are independent of the size
and ordering studies. The value $M=20$ was fixed before this evaluation;
the 80 raw trace estimates are averaged before taking positive parts.
The corresponding chi-square procedure is denoted by $\mathrm{Chi}_{20}$.
Both pooled and baseline procedures use the same $Q_{\mathrm{full}}$,
$B_{\mathrm{full}}$ and $U_{\mathrm{full}}$.

\begin{table}[!ht]
\centering\small
\setlength{\tabcolsep}{4pt}
\caption{Empirical sizes and failure rates for pooled calibration at $p=50$,
$M=20$, and $N=3{,}000$ per setting. All entries are percentages.}
\label{tab:pooled_validation}
\begin{tabular*}{\textwidth}{@{\extracolsep{\fill}}rrcrrrrr@{}}
\toprule
& & & \multicolumn{3}{c}{Empirical size} & \multicolumn{2}{c}{Failure rate} \\
\cmidrule(lr){4-6}\cmidrule(l){7-8}
Model & $n$ & $\bm{\Sigma}_0$ & F & $F_{20}$ & $\mathrm{Chi}_{20}$ & F & $F_{20}$ \\
\midrule
1 & 20 & I & 5.70 & 6.00 & 6.27 & 9.50 & 0.03 \\
1 & 20 & S & 7.90 & 7.70 & 7.87 & 14.13 & 0.43 \\
1 & 50 & I & 6.33 & 5.60 & 5.63 & 0.83 & 0.00 \\
1 & 50 & S & 8.37 & 6.87 & 6.93 & 7.47 & 0.00 \\
2 & 20 & I & 5.93 & 6.20 & 6.40 & 9.60 & 0.00 \\
2 & 50 & I & 6.20 & 5.70 & 5.73 & 0.93 & 0.00 \\
\bottomrule\end{tabular*}\end{table}

Table~\ref{tab:pooled_validation} shows that pooling reduces the maximum
failure rate from 14.13\% to 0.43\%. For Model~1 with I, the RRMSE of
$\widehat J$ decreases from 0.869 to 0.302 at $n=20$ and from 0.383 to
0.138 at $n=50$. For S, the corresponding reductions are from 2.745 to
1.636 and from 1.642 to 0.819. The empirical sizes of $F_{20}$ range from 5.60\% to 7.70\%. For S, the sizes at $n=20$ and 50 are
7.70\% and 6.87\%, with 95\% Wilson intervals $[6.80,8.71]\%$ and
$[6.02,7.83]\%$. For I at $n=20$, the size is 6.00\%, with interval
$[5.21,6.91]\%$. Pooling reduces the size for S at $n=50$ by 1.50
percentage points, with paired standard error 0.398 percentage points.
The F--chi-square differences are 0.033--0.267 percentage points. These
results show that averaging substantially improves trace estimation,
although the resulting sizes remain liberal.

To examine sensitivity to row order, an exploratory study uses 200
independent datasets per setting and 20 reorderings per dataset. Its 16 settings
combine $p\in\{50,200\}$, $n\in\{20,50\}$, I or S, and either normal
or $(\chi_1^2-1)/\sqrt2$ innovations in both groups. The proportion of
datasets with both rejection and nonrejection among the 20 baseline
orderings ranges from 9.0\% to 35.0\%. The maximum occurs for normal
innovations with I and $(p,n)=(200,20)$, with Monte Carlo standard error
3.38 percentage points. The values of $Q_{\mathrm{full}}$ and
$B_{\mathrm{full}}$ remain invariant to numerical precision, so these
changes arise from trace estimation. Pooling the 80 raw estimates on
these same exploratory datasets gives rejection rates of 4.5--14.5\%.
The maximum is 29 rejections among 200 datasets under skewed innovations,
S, and $(p,n)=(200,20)$, with standard error 2.50 percentage points.

A separate ordering experiment compares F and $F_{20}$ on 50 independent
Gaussian datasets per setting at $p=50$, using ten outer reorderings and
a fixed inner partition seed. For $(n,\bm{\Sigma}_0)=(20,\mathrm I)$,
$(20,\mathrm S)$, $(50,\mathrm I)$, and $(50,\mathrm S)$, the proportions
of datasets with a decision switch decrease from 24, 20, 2, and 16\%
to 12, 4, 0, and 8\%, respectively. The mean pairwise disagreement
probabilities decrease from 7.20, 6.89, 0.40, and 6.22\% to 3.78, 1.64,
0, and 2.27\%. Both studies use datasets as the independent units. The reductions in
decision switches indicate less sensitivity to row order after averaging.

\subsection{Accuracy of the Gaussian reference calibration}

To isolate the moment-matching error, we generate $Q^{o,*}$ and $B^{o,*}$
independently from their weighted chi-square laws under Model~1, with
$p=50$, $n\in\{20,50\}$, and I or S. Each setting uses $100{,}000$ independent
reference ratios, grouping equal feature eigenvalues in the simulation.
Table~\ref{tab:reference_calibration} compares the F critical value with
the simulated 95th reference percentile and its 95\% binomial
order-statistic interval. It also reports the reference rejection
probability at the F critical value.

\begin{table}[!ht]
\centering\small
\caption{Gaussian-reference calibration at $p=50$, based on $100{,}000$
draws per setting. Quantile intervals are 95\% binomial order-statistic
intervals. Tail probabilities and their Monte Carlo standard errors
(in parentheses) are percentages.}
\label{tab:reference_calibration}
\begin{tabular*}{\textwidth}{@{\extracolsep{\fill}}rcrlr@{}}
\toprule
$n$ & $\bm{\Sigma}_0$ & F critical value & Reference quantile [95\% interval] & F tail \\
\midrule
20 & I & 1.06697 & 1.06698 [1.06645, 1.06751] & 5.003 (0.069) \\
20 & S & 1.48094 & 1.52775 [1.51863, 1.53611] & 5.826 (0.074) \\
50 & I & 1.06639 & 1.06692 [1.06639, 1.06744] & 5.138 (0.070) \\
50 & S & 1.47454 & 1.51957 [1.51118, 1.52768] & 5.808 (0.074) \\
\bottomrule\end{tabular*}\end{table}

Under I, all feature eigenvalues are equal and the reference ratio is
exactly F distributed; the reported deviations are consistent with Monte
Carlo error. Under S, the F critical value falls below the reference
percentile at both sample sizes. Its rejection probabilities are about
5.8\%, even though no trace parameters are estimated. This discrepancy
illustrates the matching-error term retained in the theory.

These results are consistent with the persistent matching error under a
concentrated spectrum established in Section~\ref{sec:originalspectralcheck}.

\clearpage
\section{Computational details for the simulation study}
\label{app:simulationcomputation}
\subsection{Simulation settings and numerical summaries}
The main comparisons in Tables~\ref{tab:newsize20}--\ref{tab:newpowerC}
and the extensions in Appendix~\ref{app:v13extensions} use $M=20$ permutation pairs
with fixed partition seed 2026101743. The experiment seed is
2026102141. Separate settings and null-threshold calibration use
independent random-number streams.
The main grid comprises $760{,}000$ datasets: $360{,}000$ for size,
$40{,}000$ for null calibration, and $360{,}000$ for power.
The additional comparisons use $280{,}000$ datasets, giving
$1{,}040{,}000$ in total. Each cell has $10{,}000$ new datasets.
These counts exclude the separate Gaussian-reference draws and
the earlier experiments in Appendix~\ref{app:earlierstudies}.

The code stores replicate-level results and exact random-number-generator
states at checkpoints. The experiment protocol records the implementation
source hashes. An equivalent Gram-based implementation accelerates the
pooled trace calculations; the original implementation and numerical
equivalence checks are also supplied. No statistical formula, partition
choice, or invalid-calibration convention is changed.
NR3 cumulant estimates are used without clipping or regularization.

The numerical summaries include Monte Carlo standard errors, Wilson
intervals, and paired standard errors for differences in rejection
probabilities. Intervals for adjusted power are conditional on the
estimated null thresholds. The size summaries also report bounds
obtained by treating failed calibrations as nonrejections or rejections.

\subsection{NR3 implementation and calibration diagnostic}
The statistic and trace corrections were checked against Sections 2.1 and
2.3 of arXiv:2212.12338v1. We also applied the three-cumulant mean-test
routine \texttt{ZZ2022.TSBF.3cNRT} from
\href{https://github.com/nie23wp8738/HDNRA}{HDNRA} \citep{WangZhuZhangHDNRA2024} to explicit outer-product
feature matrices. Its C++ calculations and the Gram implementation used
here agreed to numerical precision for the statistic, all corrected
traces and the calibration parameters in 18 evaluations of nine datasets,
including diffuse and equicorrelated settings with both known and estimated means. This
checks the shared mean-test component; the covariance paper's complete
simulation software was not available. The fixed source version and
numerical checks are documented with the reproducibility files.

To examine the conservative sizes, we generated 10,000 additional
Gaussian datasets for each combination of $p\in\{50,200\}$ and covariance
I or S, with $(n_1,n_2)=(50,80)$. Each dataset was evaluated using both
sample means and the known population means. Let $K_2^\circ$ denote the
population variance of the known-mean NR3 statistic. Table~\ref{tab:nrcentering}
reports empirical sizes and failures, summaries of
$\widehat K_2/K_2^\circ$, and the empirical mean and variance of the
statistic relative to $K_2^\circ$.

Under I, the empirical variance with known means is close to
$K_2^\circ$, but the mean estimate $\widehat K_2$ is about 2.9 and 8.4
times as large at $p=50$ and 200, respectively. This discrepancy persists
without mean estimation and identifies trace estimation as a substantial
source of conservative calibration in these settings. The discrepancy is
smaller under S. Gaussian observations yield nonnormal quadratic
features, so the Gaussian-unbiased trace corrections need not be unbiased
for their feature covariance at these sample sizes. At I with $p=200$,
the failures arise from nonpositive estimates of the third cumulant.

For known means under Gaussian I, this variance overestimate can also
be calculated directly. Let $\bm{\Omega}$ be the covariance of
$\operatorname{vec}(\bm X\bm X^\top)$, and put $s_p=2p(p+1)$ and
$\kappa_p=4p(2p^2+5p+5)$. If $\widehat s_i$ denotes the corrected
within-group estimate of $\tr(\bm{\Omega}^2)$ used by NR3, then
\[
 \E\widehat s_i=s_p+
 \frac{n_i-1}{n_i(n_i+1)}\kappa_p.
\]
The mixed trace estimate is unbiased, giving
\[
 \frac{\E\widehat K_2}{K_2^\circ}
 =1+\frac{\kappa_p}{s_p}
 \frac{\{n_1^2(n_1+1)\}^{-1}+\{n_2^2(n_2+1)\}^{-1}}
      {\{n_1(n_1-1)\}^{-1}+2/(n_1n_2)+\{n_2(n_2-1)\}^{-1}}.
\]
At $(n_1,n_2)=(50,80)$, this ratio is 2.888 and 8.386 for $p=50$
and 200, respectively, in close agreement with the known-mean estimates
in Table~\ref{tab:nrcentering}. This calculation concerns independent
features formed using the known means.
Condition C4 of the preprint requires $p^2/n\to c\in(0,\infty)$,
where $n=n_1+n_2$. With proportional sample sizes, the relative bias above
is $O(p/n_{\min})$ and vanishes in that regime; the present diagnostic
describes its magnitude at the finite sample sizes considered here.

\begin{table}[!ht]
\centering\small
\setlength{\tabcolsep}{3.5pt}
\caption{NR3 diagnostics under Gaussian sampling, $(n_1,n_2)=(50,80)$
and $N=10{,}000$ per setting. Size and failure frequencies are percentages.
Mean and median refer to $\widehat K_2/K_2^\circ$;
$\overline T$ and $s_T^2$ are the empirical mean and variance of the statistic.}
\label{tab:nrcentering}
\begin{tabular}{cclrrrrrr}
\toprule
Cov. & $p$ & Means & Size & Failure & Mean & Median
 & $\overline T/\sqrt{K_2^\circ}$ & $s_T^2/K_2^\circ$ \\
\midrule
I & 50 & Estimated & 1.11 & 0.00 & 2.695 & 2.645 & 0.422 & 0.959 \\
I & 50 & Known & 0.31 & 0.00 & 2.885 & 2.832 & -0.011 & 0.997 \\
I & 200 & Estimated & 0.24 & 3.39 & 7.834 & 7.738 & 1.712 & 0.992 \\
I & 200 & Known & 0.00 & 3.99 & 8.381 & 8.264 & 0.003 & 1.023 \\
S & 50 & Estimated & 4.89 & 0.00 & 1.245 & 0.995 & 0.070 & 0.987 \\
S & 50 & Known & 4.29 & 0.00 & 1.332 & 1.066 & 0.002 & 1.019 \\
S & 200 & Estimated & 4.39 & 0.00 & 1.267 & 1.000 & 0.066 & 0.939 \\
S & 200 & Known & 3.99 & 0.00 & 1.359 & 1.074 & -0.006 & 0.965 \\
\bottomrule
\end{tabular}
\end{table}

\subsection{Implementation and computation times}
The simulation implementation reuses the observation Gram matrix across
the fixed partitions. Its arithmetic cost is
$O\{(n_1+n_2)^2p+M(n_1+n_2)^2\}$, compared with
$O\{M(n_1+n_2)^2p\}$ for the direct block-vector implementation.
An optional single-thread C routine evaluates the same trace sums;
the supplied Python code falls back to the original calculation when
a compiler is unavailable. Across 278 numerical checks, including
21 invalid calibrations, the rejection decisions and failure flags
agreed. The largest absolute difference in a valid $p$-value was
$1.20\times10^{-14}$.

We timed the implementations on fixed Gaussian datasets using sequential
single-thread evaluations after the simulations had finished. Each
reported time is the median of 100 evaluations after three warm-up
calls, excluding data generation and compilation. Across
$p\in\{50,200\}$ and $(n_1,n_2)\in\{(20,30),(50,80)\}$, the median
times are 0.64--1.58 milliseconds for the accelerated $F_{20}$,
1.57--11.46 milliseconds for the direct implementation, and
0.15--0.38 milliseconds for NR3. At $(p,n_1,n_2)=(200,100,160)$,
the corresponding times are 3.80, 35.61, and 1.31 milliseconds.
Both $F_{20}$ implementations include the additional four-partition
baseline diagnostics. These measurements describe the supplied
implementations and depend on hardware.

\FloatBarrier
\clearpage
\section{Additional power and mixed-innovation comparisons}
\label{app:v13extensions}
All results in this appendix use the fixed $F_{20}$ algorithm of
Section~\ref{sec:releasealgorithm} and $10{,}000$ independent datasets
per setting. The data-generating distributions and sample means are
as in Section~\ref{sec:simulation}. These results supplement the main
comparison; they are separate from the earlier four-partition F
experiments in Appendix~\ref{app:earlierstudies}.

\subsection{Dimension and sample size}
Under covariance I, we take $(p,n_1,n_2)=(200,50,80)$ or $(50,80,120)$
and use either G or K innovations in both groups. The dense, sparse,
and rank-one alternatives have $\delta=0.30$, 1, and 2, respectively.
Table~\ref{tab:v13extendedpower} reports the twelve resulting power
comparisons. Each innovation--dimension--sample-size combination has
a separate set of $10{,}000$ null datasets for calibrating adjusted
power. Null validation for $(200,50,80)$ uses the independent size
datasets in the main study; two additional sets of $10{,}000$ null
datasets validate the adjusted rules at $(50,80,120)$.
Across the original and additional calibrated models, the independent
adjusted null rejection rates are 4.36--5.59\%.

\begin{table}[!ht]
\centering\footnotesize
\setlength{\tabcolsep}{4pt}
\caption{Raw/size-adjusted power (percent) for the dimension and sample-size extensions under covariance I. Each cell and each independent null calibration use $10{,}000$ datasets; the value in parentheses after the alternative is $\delta$.}
\label{tab:v13extendedpower}
\begin{tabular*}{\textwidth}{@{\extracolsep{\fill}}crclrrrr@{}}
\toprule
Dist. & $p$ & $(n_1,n_2)$ & Alternative & $F_{20}$ & $\mathrm{Chi}_{20}$ & LC & NR3 \\
\midrule
G & 200 & $(50,80)$ & dense (0.3) & 30.3/29.4 & 30.5/29.3 & 23.2/23.0 & 5.6/33.2 \\
G & 200 & $(50,80)$ & sparse (1) & 10.6/10.1 & 10.7/10.1 & 9.7/9.7 & 0.7/9.6 \\
G & 200 & $(50,80)$ & rank (2) & 9.4/9.1 & 9.7/9.1 & 9.0/8.9 & 0.7/8.8 \\
K & 200 & $(50,80)$ & dense (0.3) & 23.9/24.2 & 24.0/24.2 & 26.7/18.6 & 0.0/27.6 \\
K & 200 & $(50,80)$ & sparse (1) & 10.1/10.3 & 10.2/10.3 & 15.1/9.9 & 0.0/10.8 \\
K & 200 & $(50,80)$ & rank (2) & 9.4/9.6 & 9.6/9.7 & 13.9/9.5 & 0.0/10.0 \\
G & 50 & $(80,120)$ & dense (0.3) & 49.2/48.6 & 49.4/48.6 & 44.3/42.5 & 34.0/52.2 \\
G & 50 & $(80,120)$ & sparse (1) & 65.7/65.2 & 65.9/65.1 & 63.8/62.3 & 48.9/66.4 \\
G & 50 & $(80,120)$ & rank (2) & 53.8/53.2 & 53.9/53.2 & 52.7/51.2 & 36.4/53.6 \\
K & 50 & $(80,120)$ & dense (0.3) & 28.2/26.9 & 28.3/26.9 & 45.9/22.5 & 1.0/27.9 \\
K & 50 & $(80,120)$ & sparse (1) & 34.3/33.2 & 34.4/33.2 & 53.6/31.8 & 2.3/33.7 \\
K & 50 & $(80,120)$ & rank (2) & 28.8/27.6 & 28.9/27.5 & 46.2/27.8 & 2.2/28.4 \\
\bottomrule
\end{tabular*}
\end{table}

\subsection{Changes in correlations with fixed marginal variances}
Let $\bm{e}_r$ denote the $r$th coordinate vector and set
\[
 \bm{C}_p=\sum_{r=1}^{5}
  \left(\bm{e}_{2r-1}\bm{e}_{2r}^{\top}
       +\bm{e}_{2r}\bm{e}_{2r-1}^{\top}\right).
\]
We take $\bm{\Sigma}_1=\bm{I}_p$ and
$\bm{\Sigma}_2=\bm{I}_p+\delta\bm{C}_p$. Thus, only the five
correlations $(1,2)$, $(3,4)$, $(5,6)$, $(7,8)$, and $(9,10)$ and their symmetric
entries change. Every marginal variance remains one. The eigenvalues
are $1-\delta$ and $1+\delta$, each with multiplicity five, and one
with multiplicity $p-10$, so the matrices are positive definite for
the selected values of $\delta$.

We use $(n_1,n_2)=(50,80)$ and G or K innovations, with
$\delta\in\{0.30,0.60\}$ at $p=50$ and $\delta=0.60$ at $p=200$.
Table~\ref{tab:v13correlationpower} gives raw and size-adjusted powers.
The null calibration samples are the independent samples already
generated for the corresponding I model. The squared Frobenius
difference is $10\delta^2$, independent of $p$; the comparisons across
dimensions therefore retain the same covariance discrepancy.

\begin{table}[!ht]
\centering\small
\setlength{\tabcolsep}{4pt}
\caption{Raw/size-adjusted power (percent) for changes in five correlations with all marginal variances fixed at one; $(n_1,n_2)=(50,80)$ and $N=10{,}000$ per setting. Adjustment uses $10{,}000$ independent null datasets.}
\label{tab:v13correlationpower}
\begin{tabular*}{\textwidth}{@{\extracolsep{\fill}}crrrrrr@{}}
\toprule
Dist. & $p$ & $\delta$ & $F_{20}$ & $\mathrm{Chi}_{20}$ & LC & NR3 \\
\midrule
G & 50 & 0.3 & 8.9/8.6 & 9.0/8.6 & 8.8/8.5 & 2.2/8.4 \\
G & 50 & 0.6 & 29.3/28.7 & 29.6/28.7 & 28.5/27.8 & 11.2/28.5 \\
K & 50 & 0.3 & 7.9/7.6 & 8.1/7.6 & 19.8/7.4 & 0.0/7.6 \\
K & 50 & 0.6 & 17.3/16.8 & 17.6/16.7 & 34.5/15.1 & 0.1/16.6 \\
G & 200 & 0.6 & 8.6/8.1 & 8.7/8.1 & 8.1/8.0 & 0.4/7.7 \\
K & 200 & 0.6 & 7.4/7.5 & 7.5/7.5 & 12.1/7.7 & 0.0/7.9 \\
\bottomrule
\end{tabular*}
\end{table}

\subsection{Different innovation laws under the covariance null}
To allow different fourth-order structures in the two groups, we
generate G innovations in group 1 and K innovations in group 2,
while setting $\bm{\Sigma}_1=\bm{\Sigma}_2=\bm{\Sigma}_0$.
We use $p=200$, covariance I or S, and
$(n_1,n_2)=(20,30)$ or $(50,80)$.
Table~\ref{tab:v13mixedsize} reports empirical sizes and calibration
failure rates for all four methods. This extends the earlier G/$t_{10}$
mixed-innovation checks to unequal sample sizes and skewed observations.

\begin{table}[!ht]
\centering\small
\setlength{\tabcolsep}{4pt}
\caption{Empirical sizes (percent) under equal covariance matrices and different innovation distributions: Gaussian in group 1 and standardized $\chi^2_1$ in group 2. Here $p=200$ and $N=10{,}000$ per setting. The last two columns give calibration failure percentages.}
\label{tab:v13mixedsize}
\begin{tabular*}{\textwidth}{@{\extracolsep{\fill}}ccrrrrrr@{}}
\toprule
Cov. & $(n_1,n_2)$ & $F_{20}$ & $\mathrm{Chi}_{20}$ & LC & NR3 & $F_{20}$ fail & NR3 fail \\
\midrule
I & $(20,30)$ & 5.87 & 6.09 & 5.61 & 0.00 & 0.00 & 0.40 \\
I & $(50,80)$ & 5.30 & 5.39 & 5.83 & 0.00 & 0.00 & 0.01 \\
S & $(20,30)$ & 9.20 & 9.38 & 9.83 & 4.44 & 0.54 & 0.00 \\
S & $(50,80)$ & 6.97 & 7.02 & 9.57 & 4.21 & 0.00 & 0.00 \\
\bottomrule
\end{tabular*}
\end{table}

\FloatBarrier
\clearpage
\section{Numerical checks of local power and the spiked null}
\label{app:theoryexperiments}
We examine the two predictions in Section~\ref{sec:fullpower} and
\eqref{eq:spikedsize} using Gaussian observations. The sample-size paths are
$(p,n_1,n_2)=(50,25,40),(100,50,80),(200,100,160)$.
Each original-data setting uses 10,000 independent replications. $F_{20}$ uses
the same 20 data-independent permutations, with partition seed 2026101743,
in every replication of a given sample-size pair. Calibration failures are
counted as nonrejections. The choices in these additional experiments were
fixed before their simulation. The plotted intervals are 95\% Wilson
intervals for the corresponding binomial probabilities.

For the local alternatives, let $\bm{\Sigma}_1=\bm{I}_p$ and
$\bm{\Sigma}_2=(1+t)\bm{I}_p$. Write $d=p(p+1)/2$ and
$v=2/n_1+2(1+t)^2/n_2$. Then
\[
 \delta_n=pt^2,\quad b^o=dv,\quad J^o=dv^2,\quad
 R^o=\frac{4d}{n_1^2(n_1-1)}+
     \frac{4d(1+t)^4}{n_2^2(n_2-1)},\qquad
 a^o=\{2(J^o+R^o)\}^{1/2}.
\]
We solve $pt^2/a^o=h$ for $h=0,1,2$, evaluating $a^o$ at the
alternative. The linear-term ratio is
$r_{\mathrm{lin}}=\bm{\Delta}^{\top}\bm{V}^o\bm{\Delta}/(a^o)^2
=vpt^2/(a^o)^2$; the variance of the linear term divided by
$(a^o)^2$ is $4r_{\mathrm{lin}}$. The nonzero ratios decrease along the path.
The covariance eigenvalues remain bounded and the feature spectrum is
diffuse, so these alternatives satisfy both local-power conditions.

The finite central-reference prediction in Table~\ref{tab:localtheory}
is $1-G_n^o(t_n^o-h)$, evaluated with 1,000,000 independent draws per
setting. Specifically, independently,
\[
 Q_0^{o,*}=v\chi_d^2,\qquad
 B^{o,*}=\frac{2\chi^2_{d(n_1-1)}}{n_1(n_1-1)}+
         \frac{2(1+t)^2\chi^2_{d(n_2-1)}}{n_2(n_2-1)}.
\]
This is the central reference in the theorem, with the signal entering
through the threshold. It is not an exact finite-sample power formula
for the original observations. The normal-limit prediction is
$1-\Phi(z_{0.95}-h)$.

\begin{table}[!ht]
\centering\small
\caption{Diffuse local alternatives. Rejection probabilities and predictions
are reported as proportions. $r_{\mathrm{lin}}$ is the linear-term ratio.}
\label{tab:localtheory}
\begin{tabular}{rrrrrrrrr}
\toprule
$p$ & $h$ & $t$ & $a^o$ & $r_{\mathrm{lin}}$ & $F_{20}$ & SE & Central ref. & Normal \\
\midrule
50 & 0 & 0.0000 & 6.6286 & 0.00000 & 0.0583 & 0.0023 & 0.0498 & 0.0500 \\
50 & 1 & 0.4309 & 9.2831 & 0.01965 & 0.2647 & 0.0044 & 0.2496 & 0.2595 \\
50 & 2 & 0.6679 & 11.1510 & 0.03929 & 0.6179 & 0.0049 & 0.6259 & 0.6388 \\
100 & 0 & 0.0000 & 6.5637 & 0.00000 & 0.0550 & 0.0023 & 0.0498 & 0.0500 \\
100 & 1 & 0.2866 & 8.2122 & 0.00991 & 0.2653 & 0.0044 & 0.2545 & 0.2595 \\
100 & 2 & 0.4286 & 9.1836 & 0.01982 & 0.6211 & 0.0049 & 0.6318 & 0.6388 \\
200 & 0 & 0.0000 & 6.5317 & 0.00000 & 0.0526 & 0.0022 & 0.0495 & 0.0500 \\
200 & 1 & 0.1950 & 7.6050 & 0.00498 & 0.2591 & 0.0044 & 0.2578 & 0.2595 \\
200 & 2 & 0.2858 & 8.1703 & 0.00995 & 0.6378 & 0.0048 & 0.6361 & 0.6388 \\
\bottomrule
\end{tabular}
\end{table}

\FloatBarrier
\newpage
For the spiked null, both groups have covariance
$\bm{\Sigma}=0.7\bm{I}_p+0.3\bm{1}_p\bm{1}_p^\top$.
Let $\lambda_1=0.3p+0.7$ and $\lambda_0=0.7$.
The three distinct eigenvalues of $\bm{\Gamma}$ are
$2\lambda_1^2$, $2\lambda_1\lambda_0$, and $2\lambda_0^2$, with
multiplicities $1$, $p-1$, and $p(p-1)/2$, respectively.
We compare $F_{20}$ with a population-parameter benchmark that uses the same
original-data statistic $Q/B$ and the critical value
$q_{0.95}(d_1^o,d_2^o)$. Thus the benchmark removes estimation of the
degrees of freedom, while retaining the finite-sample behavior of the
original-data statistic.

The Gaussian-reference column uses 100,000 independent draws per path
point from the exact weighted-chi-square reference distribution, using
the same population critical value. Each multiplicity is incorporated in
the chi-square degrees of freedom. All numerator and denominator
chi-square variables are drawn independently. The limiting rejection
probability from \eqref{eq:spikedsize} is
$\alpha_\infty(0.3)=0.06143$.
The Gaussian reference and the population benchmark distinguish matching
error from finite-sample effects in the original observations. None of
these comparisons requires the spike to satisfy the diffuse-spectrum
condition; its failure is precisely the source of the nonnominal limit.

\begin{table}[!ht]
\centering\small
\caption{Spiked Gaussian null. Parentheses contain Monte Carlo standard
errors. $F_{20}$ and the population benchmark use 10,000 replications; the
Gaussian reference uses 100,000.}
\label{tab:spiketheory}
\begin{tabular}{rrrrrr}
\toprule
$p$ & $(n_1,n_2)$ & $d_1^o$ & $F_{20}$ & Population benchmark & Gaussian reference \\
\midrule
50 & $(25,40)$ & 28.654 & 0.0797 (0.0027) & 0.0488 (0.0022) & 0.0589 (0.0007) \\
100 & $(50,80)$ & 32.293 & 0.0711 (0.0026) & 0.0560 (0.0023) & 0.0600 (0.0008) \\
200 & $(100,160)$ & 34.378 & 0.0667 (0.0025) & 0.0567 (0.0023) & 0.0600 (0.0008) \\
\bottomrule
\end{tabular}
\end{table}
Across these 120,000 original-data replications, $F_{20}$ had 4 calibration failures.
The data-generation seeds are 2026102100--2026102111; reference seeds are
2026103100--2026103111. Per-replication outputs and the exact counts are
provided with the reproduction code.

\FloatBarrier
\end{document}